\documentclass[review,3p,square]{elsarticle}

\usepackage{amssymb}
\usepackage{amsmath}
\usepackage{amsthm}
\usepackage{bbm}
\usepackage{enumitem}

\usepackage{fancyhdr}
\usepackage{hyperref}
\hypersetup{%
  pdftitle={BSDEs},%
  pdfsubject={BSDEs},%
  pdfauthor={ShengJun FAN, Ying Hu},%
  pdfkeywords={BSDEs},%
  pdfstartview=FitH,%
  CJKbookmarks=true,%
  bookmarksnumbered=true,%
  bookmarksopen=true,%
  colorlinks=true, linkcolor=blue, urlcolor=blue, citecolor=blue, %
}

\usepackage{cleveref}
\crefname{thm}{Theorem}{Theorems}
\crefname{pro}{Proposition}{Propositions}
\crefname{lem}{Lemma}{Lemmas}
\crefname{rmk}{Remark}{Remarks}
\crefname{cor}{Corollary}{Corollaries}
\crefname{dfn}{Definition}{Definitions}
\crefname{ex}{Example}{Examples}
\crefname{section}{Section}{Sections}
\crefname{subsection}{Subsection}{Subsections}

\newcommand{\eps}{\varepsilon}

\newcommand{\To}{\rightarrow}
\newcommand{\as}{{\rm d}\mathbb{P}\times{\rm d} t-a.e.}
\newcommand{\ass}{{\rm d}\mathbb{P}\times{\rm d} s-a.e.}
\newcommand{\ps}{\mathbb{P}-a.s.}

\newcommand{\F}{\mathcal{F}}
\newcommand{\E}{\mathbb{E}}

\newcommand{\s}{\mathcal{S}}
\newcommand{\lcal}{\mathcal{L}}
\newcommand{\mcal}{\mathcal{M}}
\newcommand{\ecal}{\mathcal{E}}

\newcommand{\hcal}{\mathcal{H}}

\newcommand{\T}{[0,T]}

\newcommand{\R}{{\mathbb R}}

\newcommand{\RE}{\forall}

\newcommand {\Dis}{\displaystyle}

\newtheorem{thm}{Theorem}[section]
\newtheorem{lem}[thm]{Lemma}
\newtheorem{pro}[thm]{Proposition}
\newtheorem{rmk}[thm]{Remark}
\newtheorem{cor}[thm]{Corollary}

\newtheorem{ex}[thm]{Example}

\journal{arXiv}
\begin{document}
\begin{frontmatter}

\title{{Systems of quadratic BSDEs with unbounded terminal data}\tnoteref{found}}
\tnotetext[found]{This work is supported by National Natural Science Foundation of China (Nos. 12631016, 12171471, 12031009 and 11631004), by Key Laboratory of Mathematics for Nonlinear Sciences (Fudan University), Ministry of Education, Handan Road 220, Shanghai 200433, China, and by Lebesgue Center of Mathematics ``Investissements d'avenir" program-ANR-11-LABX-0020-01.
\vspace{0.2cm}}


\author[Fan]{Shengjun Fan} \ead{shengjunfan@cumt.edu.cn}
\author[Hu]{Ying Hu} \ead{ying.hu@univ-rennes1.fr}
\author[Tang]{Shanjian Tang\corref{cor}} \ead{sjtang@fudan.edu.cn} \vspace{-0.5cm}

\affiliation[Fan]{organization={School of Mathematics, China University of Mining and Technology},
            city={Xuzhou 221116},
            country={China}}

\affiliation[Hu]{organization={Univ. Rennes, CNRS, IRMAR-UMR6625},
            city={F-35000, Rennes},
            country={France}}

\affiliation[Tang]{organization={Department of Finance and Control Sciences, School of Mathematical Sciences, Fudan University},
            city={Shanghai 200433},
            country={China}}

\cortext[cor]{Corresponding author\vspace{0.2cm}}

\begin{abstract}
The paper is devoted to global solution of multidimensional backward stochastic differential equations (BSDEs) with non-convex quadratic generators and unbounded terminal data. The generators of some individual BSDEs of the system are allowed to have interactions among their second unknown variables and can be non-Lipschitz continuous with general growth in the first unknown variable $y$ and non-convex/non-concave in the second unknown variable $z$. We extend the results of the authors' previous paper~[JDE, 2023] by incorporating both the no convexity/concavity of the generators and the interaction among the second unknown variables of individual BSDEs, and update the assumptions of Theorems 2.10 and 2.27 stated in the authors' previous paper~[EJP, 2025] so as to fill in a gap in the \vspace{0.2cm}original proofs.
\end{abstract}

\begin{keyword}
Multidimensional BSDEs \sep Quadratic growth \sep Unbounded terminal data \sep Global solution, \\ \hspace*{1.8cm}  Existence and uniqueness \sep Stability. \vspace{0.2cm}

\MSC[2010] 60H10\vspace{0.2cm}
\end{keyword}

\end{frontmatter}
\vspace{-0.4cm}

\section{Introduction}
\label{sec:1-Introduction}
\setcounter{equation}{0}

Given a terminal time $T\in (0,+\infty)$ and a complete filtered probability space $(\Omega, \F, (\F_t)_{t\in\T}, \mathbb{P})$ which satisfies the usual conditions. Suppose that the filtration $(\F_t)_{t\in\T}$ is generated by a $d$-dimensional Brownian motion $(B_t)_{t\in\T}$ and augmented by the $\mathbb{P}$-null sets. All equalities and inequalities between
random variables and processes should be understood in the senses of $\ps$ and $\as$, respectively. Consider the following backward stochastic differential equation (BSDE in short):
\begin{equation}\label{eq:1.1}
  Y_t=\xi+\int_t^T g(s,Y_s,Z_s){\rm d}s-\int_t^T Z_s {\rm d}B_s, \ \ t\in\T,
\end{equation}
where $\xi\in \F_T$ is an $d$-dimensional random vector called the terminal value, and for each pair $(y,z)$, the random function
$$g(\omega, t, y, z):\Omega\times\T\times\R^n\times\R^{n\times d}\To \R^n$$
is $(\F_t)$-progressively measurable and is called the generator. A solution to such an equation is a pair of $(\F_t)$-progressively measurable processes $(Y_t,Z_t)_{t\in\T}$ taking values in $\R^n\times\R^{n\times d}$ which almost surely satisfies \eqref{eq:1.1}. Linear BSDEs were solved by \citet{Bismut1973JMAA}, and  nonlinear BSDEs were solved in \citet{PardouxPeng1990SCL}, where the generators are uniformly Lipschitz continuous  in both unknown variables $(y,z)$. Subsequently,  extensions, simplifications and applications have received an extensive study on these equations. In particular, BSDEs whose generators have a quadratic growth in $z$, the so-called quadratic BSDEs, have attracted lots of attentions in the last three decades. See \citet{ElKarouiPengQuenez1997MF}, \citet{YongZhou1999}, \citet{HuImkellerMuller2005AAP}, \citet{HuTangWang2022SPA}, \citet{FanWangYong2025JDE}, \citet{HaoHuWen2025AAP} and the authors' paper~\cite{FanHuTang2026MF} for more details among others.

For the one-dimensional case, since the first existence and uniqueness result of \citet{Kobylanski2000AP}  on quadratic BSDEs with bounded terminal values, there have been  numerous subsequent intensive studies on this topic. See, for example,  \citet{Tevzadze2008SPA}, \citet{BriandElie2013SPA}, \citet{Fan2016SPA} and \citet{LuoFan2018SD} for the bounded terminal value case, and \citet{BriandHu2006PTRF,BriandHu2008PTRF},
\citet{DelbaenHuRichou2011AIHPPS}, \citet{BarrieuElKaroui2013AoP} and \citet{DelbaenHuRichou2015DCD} for unbounded terminal value case. In the three authors' previous work \cite{FanHuTang2020CRM}, an existence and uniqueness result was given on solution of a one-dimensional quadratic BSDE, whose terminal value has exponential moments of an arbitrary order and whose generator $g$ may be non-Lipschitz continuous with general growth in $y$ and may be non-convex/non-concave in $z$. For more details, see \cref{lem:A.1} and \cref{rmk:A.2} in the appendix of the paper. When the terminal value is bounded, a solution $(Y, Z)$ of (a system of) quadratic BSDEs can have a finite exponential moment of $(\int_0^T|Z_t|^2\, dt+\max_{t\in [0,T]}|Y_t|)$, due to the bounded mean oscillation (BMO) martingale property of $\int_0^\cdot Z_t\, dW_t$. Once the terminal value is unbounded, the stochastic integral $\int_0^\cdot Z_t\, dW_t$ goes beyond being a BMO martingale,  and an a priori estimate of an  exponential moment of  $(\int_0^T|Z_t|^2\, dt+\max_{t\in [0,T]}|Y_t|)$ is crucial and then is a key point.

The present paper focuses on a system of quadratic BSDEs, which naturally arise in various applications, including financial price-impact models, financial market equilibrium problems for several interacting agents, nonzero-sum risk-sensitive stochastic differential games, stochastic equilibrium problems in incomplete financial markets. See for example the references mentioned in this paragraph and subsequent paragraph for more details. Solvability of a system of quadratic BSDEs is a big challenge, for conventional tools used in the scalar quadratic BSDEs, like Girsanov's transform, comparison theorem and monotone convergence, fail to work. A striking counterexample with a rather simple generator, given by \citet{FreiDosReis2011MFE},  demonstrates that, in contrast to the one-dimensional case, a system of quadratic BSDEs has to face the endogenous unbounded interaction among the second unknown variables, and thus might easily fail to have a global solution even though the terminal value is bounded. In fact, global solvability of multidimensional quadratic BSDEs with unbounded terminal values has been posed as an open problem in~\citet[Section 5, page 270]{Peng1999Book}.

Assuming that the terminal value is sufficiently small in the supremum norm, \citet{Tevzadze2008SPA} gave the first existence and uniqueness result on global solution of multidimensional quadratic BSDEs with the BMO martingale tool and the fixed-point principle, which has inspired some subsequent works under various ``smallness" assumptions either on the generators or on the terminal values, see \citet{Frei2014SPA}, \citet{Kardaras2022Book}, \citet{JamneshanKupperLuo2017ECP}, \citet{KramkovPulido2016AAP}, \citet{KramkovPulido2016SIAM} and \citet{HarterRichou2019EJP} for more details. In the Markovian setting, \citet{cheridito2015Stochastics} and \citet{XingZitkovic2018AoP} addressed the global solvability of multidimensional quadratic BSDEs with bounded terminal values and
weak regularity assumptions on the generator by Girsanov's transform and analytic PDE methods. And, by virtue of the fixed-point argument along with the BMO martingale tool, \citet{cheridito2015Stochastics}, \citet{HuTang2016SPA}, \citet{Luo2020EJP}, \citet{FanHuTang2023JDE,FanHuTang2025EJP} established the existence and uniqueness on local and global bounded solutions of multidimensional BSDEs with bounded terminal values and projectable quadratic generators, sub-quadratic generators, diagonally quadratic generators and triangularly quadratic generators, respectively. In particular, in the three authors' previous work~\cite{FanHuTang2023JDE}, we also proved the global solvability for a class of multidimensional quadratic BSDEs with unbounded terminal values, where the $i$th components $g^i$ of the generators $g(\omega,t,y,z)$ vary only with $(\omega,t,y,z^i)\in \Omega\times [0,T]\times \R^n\times \R^{1\times d}$, where $z^i$ represents the $i$th row of the matrix $z$, and they are Lipschitz continuous in $y$ and convex/concave with a quadratic growth in $z^i$. This seems to be the first general result on the global solvability of multidimensional quadratic BSDEs with general unbounded terminal values.

The present paper is a continuation of the authors' previous work~\cite{FanHuTang2023JDE} on the case of unbounded terminal values, which is restricted to component-wisely convex/concave generators. We first establish (see \cref{thm:2.2} in \cref{sec:2-statement of main result}) an existence and uniqueness of solutions of multidimensional quadratic BSDEs with unbounded terminal values, where the $i$th components $g^i$ of the generators $g$ depend only on the $i$th row $ z^i$ of the matrix $z$ in a quadratic way. The $i$th components $g^i$ may be non-Lipschitz continuous and have a rather general growth in $y$, and neither convex nor concave in $z^i$, extending  Theorem 2.8 on multidimensional quadratic BSDEs with unbounded terminal values and convex/concave generators in~\citet{FanHuTang2023JDE} to the case of multidimensional quadratic BSDEs with some components of the generators being neither convex nor concave. Our \cref{thm:2.2} is proved with the iteration argument and the $\theta$-technique (initially only used for the convex/concave function), and the uniform a priori estimates on solutions of one-dimensional quadratic BSDEs along with Doob's maximal inequality for martingales play natural and crucial roles in its proof. Furthermore, in \cref{thm:2.5} we establish the stability result in the case of multidimensional quadratic BSDEs with unbounded terminal conditions and non-convex/non-concave generators addressed in \cref{thm:2.2}. The $\theta$-technique, the uniform a priori estimates on one-dimensional quadratic BSDEs and Doob's maximal inequality for martingales play an important role again in the proof of  \cref{thm:2.5}. Next, in \cref{thm:2.7} we give a more general solvability of multidimensional quadratic BSDEs with unbounded terminal values, where the $i$-th components $g^i(t,y,z)$ of the generators $g(t,y,z)$ admit interactions among the last argument $z$, and can be non-Lipschitz continuous with general growth in $y$ and non-convex/non-concave in $z$. The proof of \cref{thm:2.7} is an adaptation of  that of \cref{thm:2.2}, and appeals similarly to the iteration argument, the $\theta$-technique and the uniform a priori estimates and Doob's maximal inequality, but adds novel ideas to incorporate the interactions among the second unknown variables of the generators of the systems. Finally, in \cref{thm:2.13} we present an existence and uniqueness result on the unbounded global solution of a system of quadratic BSDEs with a special type of unbounded terminal value, which is expressed as the sum of a bounded random vector and a stochastic integral of a BMO process with respect to the Brownian motion. In this case, more general interactions among the $z$-variables of individual BSDEs are incorporated. Before that, an existence and uniqueness result on the bounded global solution of a more general system of quadratic BSDEs with bounded terminal value is first established in \cref{thm:2.11}. Note that in~\cref{thm:2.11,thm:2.13}, we update the assumptions of \citet[Theorems 2.10 and 2.27]{FanHuTang2025EJP} so as to fill in a gap in the original proofs, and  the whole ideas of the proofs are almost unaltered.

The rest of the paper is organized as follows. In the next section, we state our main results-\cref{thm:2.2,thm:2.5,thm:2.7,thm:2.11,thm:2.13}, and introduce some remarks, examples and corollaries to illustrate these results. Their proofs are provided in Section 3. Finally, in Appendix we present and prove some auxiliary results for one-dimensional quadratic BSDEs with unbounded terminal values.\vspace{0.2cm}

\section{Notations and Statement of main results}
\label{sec:2-statement of main result}
\setcounter{equation}{0}

\subsection{Notations\vspace{0.2cm}}

Denote by ${\bf 1}_A(\cdot)$ the indicator of set $A$, and ${\rm sgn}(x):={\bf 1}_{x>0}-{\bf 1}_{x\leq 0}$. Let $a\wedge b$ be the minimum of $a$ and $b$, $a^-:=-(a\wedge 0)$, $a^+:=(-a)^-$, and for any real $b^j$, let
$$\sum_{j=1}^{0}b^j=\sum_{j=n+1}^{n}b^j=0.$$
The Euclidean norm is always denoted by $|\cdot|$. For each $\F_T$-measurable one-dimensional or multidimensional random variable $\xi$, let $\E_\tau[\xi]$ denote the conditional expectation of $\xi$ with respect to $\F_\tau$ for $\tau\in \mathcal {T}_{[0,T]}$, here and hereafter, $\mathcal {T}_{[a,b]}$ stands for the set of all $(\F_t)$-stopping times valued in $[a,b]\subset \T$. Denote by $L^{\infty}(\R^n)$ the set of $\R^n$-valued $\F_T$-measurable random variables $\xi$ satisfying
$$\|\xi\|_{\infty}:=\mathop{\rm ess\ sup}_{\omega\in \Omega}|\xi(\omega)|<+\infty$$
and by $\s^{\infty}(\R^n)$ the set of $\R^n$-valued continuous $(\F_t)$-adapted processes $(Y_t)_{t\in\T}$ satisfying
$$\|Y\|_{{\s}^{\infty}}:=\mathop{\rm ess\ sup}_{(\omega,t)\in \Omega\times \T}|Y_t(\omega)| =\Big\|\sup_{t\in\T} |Y_t| \Big\|_{\infty}<+\infty.\vspace{0.2cm}$$

For each $p\geq 1$, by $\s^p(\R^n)$ we denote the set of $\R^n$-valued continuous $(\F_t)$-adapted processes $(Y_t)_{t\in\T}$ satisfying
$$\|Y\|_{{\s}^p}:=\left(\E\left[\sup_{t\in\T} |Y_t|^p\right]\right)^{1\over p}<+\infty\vspace{0.1cm}$$
and by $\hcal^p(\R^{n\times d})$ the set of $\R^{n\times d}$-valued $(\F_t)$-progressively measurable processes $(Z_t)_{t\in\T}$ satisfying
$$
\|Z\|_{\hcal^p}:=\left(\E\left[\left(\int_0^T |Z_s|^2{\rm d}s\right)^{p\over 2}\right] \right)^{1\over p}<+\infty.\vspace{0.1cm}
$$
We write $Y\in \mathcal{E}(\R^n)$ if $\exp{|Y|}\in \bigcap_{p\geq 1}\s^p(\R^n)$, and $Z\in \mcal(\R^{n\times d})$ if $Z\in \bigcap_{p\geq 1}\hcal^p(\R^{n\times d})$. Moreover, by ${\rm BMO}(\R^{n\times d})$ we denote the set of $Z\in \hcal^2(\R^{n\times d})$ satisfying
$$
\|Z\|_{\rm BMO}:=\sup_{\tau\in \mathcal {T}_{\T}}\left\|\E_{\tau}\left[\int_{\tau}^T |Z_s|^2 {\rm d}s\right]\right\|_{\infty}^{1\over 2}<+\infty,\vspace{0.2cm}
$$
and by $\mcal^{\infty}$ we denote the set of non-negative real-valued $(\F_t)$-progressively measurable process $(h_t)_{t\in\T}$ satisfying
$$
\|h\|_{{\mcal}^{\infty}}:=\sup_{\tau\in \mathcal {T}_{\T}}\left\|\E_{\tau}\left[\int_{\tau}^T h_s {\rm d}s\right]\right\|_{\infty}<+\infty.\vspace{0.2cm}
$$
For each $r>0$, by $\ecal^{\infty}(r)$ we denote the set of all $(h_t)_{t\in \T}\in \mcal^{\infty}$ satisfying
$$\|h\|_{\ecal^{\infty}(r)}:={1\over r}\ln\left(\sup_{\tau\in \mathcal{T}_{[0,T]}}\left\|\E_{\tau}\left[\exp\left(r\int_\tau^T h_s {\rm d}s\right)\right]\right\|_{\infty}\right)<+\infty,\vspace{0.1cm}$$
and by $\lcal^{\infty}$ we denote the set of all $(h_t)_{t\in \T}\in \cap_{r>0}\ecal^{\infty}(r)$ satisfying
$$
\|h\|_{{\lcal}^{\infty}}:=\left\|\int_0^T h_s {\rm d}s \right\|_{\infty}<+\infty.
$$
For $[a,b]\subset \T$, we denote by $\hcal_{[a,b]}$ the restriction to $[a,b]$ of the space $\hcal$ of processes on $[0,T]$, i.e., the space of all processes which are restrictions to $[a,b]$ of processes in $\hcal$.\vspace{0.2cm}

In the whole paper, we fix an $\R^n$-valued $\F_T$-measurable random vector $\xi$, some non-negative real-valued $(\F_t)$-progressively measurable processes $(\alpha_t)_{t\in\T}$ having a certain integrability to be specified,
$$
(\bar\alpha_t)_{t\in\T}\in \mcal^{\infty},\ \ (\tilde\alpha_t)_{t\in\T}\in \lcal^{\infty}\ \ {\rm and}\ \ (v_t)_{t\in\T}\in {\rm BMO}(\R),
$$
and some real constants $\beta,k,\bar k,\lambda,\bar\lambda,c,\bar c,\theta\geq 0$, $0<\bar\gamma\leq \gamma$ and $\delta\in [0,1)$.\vspace{0.2cm}

Finally, for $i=1,\cdots, n$, let $z^i$, $y^i$ and $g^i$ respectively represent the $i$-th row of matrix $z\in\R^{n\times d}$ and the $i$th component of the vector $y\in \R^n$ and the generator $g$. For any $U\in\R^n$ and $x\in\R$, define by $U(x;i)$ the vector in $\R^n$ whose $i$th component is $x$ and whose $j$th component is $U^j$ for any $j\neq i$. And, for each $H\in \R^{n\times d}$, $w\in \R^{1\times d}$ and $i=1,\cdots,n$, denote by $H(w;i)$ the matrix in $\R^{n\times d}$ whose $i$th row is $w$ and whose $j$th row is $H^j$ for $j\neq i$.\vspace{0.2cm}

\subsection{Statement of main results: the case of general unbounded terminal values\vspace{0.2cm}}

In this subsection, we will state several novel well-posedness results on systems of quadratic BSDEs with the general unbounded terminal values having  exponential moments of any order. The generator $g$ may be neither Lipschitz continuous nor linearly growthing in the first unknown variable $y$, and may be neither convex nor concave in the second unknown variable $z$. \vspace{0.2cm}

Let us start with introducing the following two assumptions \ref{A:H1} and \ref{A:H2} on the generator $g$ defined on $\Omega\times [0,T]\times \R^n\times \R^{n\times d}$.

\begin{enumerate}
\renewcommand{\theenumi}{(H\arabic{enumi})}
\renewcommand{\labelenumi}{\theenumi}
\item\label{A:H1} Almost everywhere  $(\omega, t)\in \Omega\times\T$, $g(\omega,t,\cdot,\cdot)$ is continuous.

\item\label{A:H2} For $i=1,\cdots,n$, the random field $f(\omega,t,y,z):=g^i(\omega,t,y,z)$ or $f(\omega,t,y,z):=-g^i(\omega,t,y(-y^i;i),-z)$ varies only with $(\omega,t,y,z^i)\in \Omega\times [0,T]\times \R^n\times \R^{1\times d}$, and $\as$, for each $(y,z)\in \R^n\times\R^{1\times d}$,
\begin{equation}\label{eq:2.2}
    {\rm sgn}(y^i)f(\omega,t,y,z)\leq \alpha_t(\omega)+\beta|y|+\frac{\gamma}{2} |z|^2
\end{equation}
and
\begin{equation}\label{eq:2.3}
|f(\omega,t,y,z)|\leq \alpha_t(\omega)+k\exp(k|y|)+\frac{\gamma}{2} |z|^2;\vspace{0.1cm}
\end{equation}
and furthermore, 
\begin{equation}\label{eq:2.4}
\hspace*{-0.3cm}\begin{array}{l}
   \Dis {\bf 1}_{\{y^i_1-\theta y^i_2>0\}}\left(f(\omega,t,y_1,z_1)-\theta f(\omega,t,y_2,z_2)\right)\vspace{0.1cm}\\
   \ \ \leq \Dis (1-\theta)\left(\alpha_t(\omega)+ \beta \left|\frac{y_1-\theta y_2}{1-\theta}\right|+\frac{\gamma}{2} \left|\frac{z_1-\theta z_2}{1-\theta}\right|^2+k(|y_1|+|y_2|)+\bar k (|z_1|^{1+\delta}+|z_2|^{1+\delta})\right)\\
   \Dis \hspace{2cm}\RE (y_1,z_1, y_2, z_2)\in \left (\R^n\times\R^{1\times d}\right)^2, \quad  \RE\theta\in (0,1)
   \end{array}
\end{equation}
 with $\bar k=0$  or  the additional inequality
\begin{equation}\label{eq:2.5}
 \Dis f(\omega,t,y,z)\geq \frac{\bar\gamma}{2}|z|^2-\beta|y|-\alpha_t(\omega), \ \ \RE (y,z)\in \R\times\R^{1\times d}.
\end{equation}
\end{enumerate}

\begin{rmk}\label{rmk:2.1}
In assumption \ref{A:H2}, all the integers $1,\cdots,n$ consists of two groups  $J_1$ and $J_2$ such that  for $i\in J_1$, $f=g^i$ satisfies \eqref{eq:2.2}-\eqref{eq:2.3}, and \eqref{eq:2.4} with $\bar k=0$  or  the additional inequality  \eqref{eq:2.5}, and for $i\in J_2$,  $f(\omega,t,y,z)=-g^i(\omega,t,y(-y^i;i),-z)$ satisfies \eqref{eq:2.2}-\eqref{eq:2.3}, and \eqref{eq:2.4} with $\bar k=0$  or  the additional inequality  \eqref{eq:2.5}, and $J_1+J_2=\{1,\cdots,n\}$. Note that either of  $J_1$ and $J_2$ might be $\emptyset$. Now, we define, for $(\omega,t,y,z)\in \Omega\times \T\times \R^n\times\R^{n\times d}$,
$$
\bar g^i(\omega,t,y,z):=\left\{
\begin{array}{l}
\hspace*{0.55cm} g^i(\omega,t,\bar y, z^i),\ \ i\in J_1;\\
-g^i(\omega,t,\bar y,-z^i),\ \ i\in J_2
\end{array}
\right.
\ \ \ {\rm with}\ \ \
\bar y^i=\left\{
\begin{array}{l}
\ \ y^i,\ \ i\in J_1;\\
-y^i,\ \ i\in J_2.
\end{array}
\right. \vspace{0.1cm}
$$
It can be directly verified that $f:=\bar g^i$ satisfies~\eqref{eq:2.2}-\eqref{eq:2.3}, and \eqref{eq:2.4} with $\bar k=0$  or the additional inequality  \eqref{eq:2.5} for all $i=1,\cdots,n$. Therefore, in assumption~\ref{A:H2}, we may assume without loss of generality that $f=g^i$ satisfies \eqref{eq:2.2}-\eqref{eq:2.3}, and \eqref{eq:2.4} with $\bar k=0$  or the additional inequality  \eqref{eq:2.5} for all $i=1,\ldots,n$. \vspace{0.2cm}
\end{rmk}


\cref{thm:2.2} below establishes the existence and uniqueness of solutions to multidimensional quadratic BSDEs with general unbounded terminal values and generators satisfying assumptions \ref{A:H1} and \ref{A:H2}, and constitutes the first main result of the paper.

\begin{thm}\label{thm:2.2}
Assume that for each $p\geq 1$,
\begin{equation}\label{eq:2.1}
\E\left[\exp\left(p|\xi|+p\int_0^T \alpha_s{\rm d}s\right)\right]<+\infty,
\end{equation}
and that the generator $g$ satisfies assumptions \ref{A:H1} and \ref{A:H2}. Then BSDE \eqref{eq:1.1} admits a unique solution $(Y_t,Z_t)_{t\in\T}$ such that $(Y,Z)\in \mathcal{E}(\R^n)\times \mcal(\R^{n\times d})$.

Moreover, for those indexes $i\in \{1,\ldots,n\}$ such that $g^i$ satisfies \eqref{eq:2.5}, and $\eps\in (0,\bar\gamma/9]$, we have
\begin{equation}\label{eq:2.8}
\E\left[\exp\left(\frac{\bar \gamma}{2}\eps\int_0^T |Z^i_t|^2{\rm d}t \right)\right] \leq \E\left[\exp\left(6\eps(1+\beta T+kT)\sup_{t\in [0,T]}|Y_t|+3\eps\int_0^T \alpha_t {\rm d}t\right)\right]. \vspace{0.42cm}
\end{equation}
\end{thm}

\begin{rmk}\label{rmk:2.3}

(i) By a similar argument as that in \citet[Proposition 3.5]{FanHu2021SPA}, it is not difficult to verify that for $i=1,\cdots,n$, \eqref{eq:2.4} holds with $\bar k=0$  for $f$ but with different parameters provided that it satisfies \eqref{eq:2.2} and any one of the following three conditions:
\begin{enumerate}
\item [a)] Almost everywhere in $\Omega\times\T$, $f(\omega,t,\cdot,\cdot)$ is convex;
\item [b)] Almost everywhere in $\Omega\times\T$, for each $(y,z)\in \R^n\times\R^{1\times d}$, $f(\omega,t,\cdot,z)$ is Lipschitz continuous and $f(\omega,t, y,\cdot)$ is convex;
\item [c)] $f(\omega,t,y,z):=l(y)q(z)$ for each $(y,z)\in \R^n\times\R^{1\times d}$, where both functions $l:\R^n\To\R$ and $q:\R^{1\times d}\To\R$ are bounded Lipschitz continuous, and the function $q(z)$ has a bounded support.
\end{enumerate}
Furthermore, in the same way as that in \citet[Proposition 3]{FanHuTang2020CRM} it can also be proved that for $i=1,\cdots,n$, the inequality \eqref{eq:2.4} holds for $f$ but with different parameters provided that it satisfies \eqref{eq:2.2} and either of the following two conditions:
\begin{enumerate}
\item [d)] Almost everywhere in $\Omega\times\T$, for each $(y,z)\in \R^n\times\R^{1\times d}$, $f(\omega,t,\cdot,z)$ is Lipschitz continuous, and $f(\omega,t, y,\cdot)$ is $\delta$-locally Lipschitz continuous,  i.e., $\as$, for each $(y,\bar y,z,\bar z)\in \R^n\times\R^n\times\R^{1\times d}\times\R^{1\times d}$,
\begin{equation}\label{eq:2.9}
|f(\omega,t,y,z)-f(\omega,t,\bar y,\bar z)|\leq \beta |y-\bar y|+\gamma (1+|z|^\delta+|\bar z|^\delta)|z-\bar z|;\vspace{-0.1cm}
\end{equation}
\item [e)] Almost everywhere in $\Omega\times\T$, for each $(y,z)\in \R^n\times\R^{1\times d}$, $f(\omega,t,\cdot,z)$ is convex, and $f(\omega,t, y,\cdot)$ is $\delta$-locally Lipschitz continuous, i.e., \eqref{eq:2.9} holds with $y=\bar y$.\vspace{0.1cm}
\end{enumerate}
In addition, we have \eqref{eq:2.4} with $\bar k=0$ (resp. \eqref{eq:2.4} and \eqref{eq:2.5}) for the sum of finite numbers of generators $f$ satisfying \eqref{eq:2.4} with $\bar k=0$ (resp. \eqref{eq:2.4} and \eqref{eq:2.5}), but with different parameters.\vspace{0.2cm}

(ii) For $i=1,\cdots,n$, letting $y_1=y_2=y$ and $z_1=z_2=z$ in the inequality \eqref{eq:2.4} yields that
$$
{\bf 1}_{\{y^i>0\}} f(\omega,t,y,z)\leq 2\bar k+\alpha_t(\omega)+(\beta+2k)|y|+\left(\frac{\gamma}{2}+2\bar k\right)|z|^{2},\vspace{0.1cm}
$$
which means that a part of the assertion stated in \eqref{eq:2.2} holds with different parameters. In addition, letting first $z_1=z_2=z$ in the inequality~\eqref{eq:2.4} and then sending $\theta\To 1$ yields that for $i=1,\cdots,n$,
$$
{\bf 1}_{\{y_1^i-y_2^i>0\}} \left(f(\omega,t,y_1,z)-f(\omega,t,y_2,z)\right)\leq \beta |y_1-y_2|.\vspace{0.2cm}
$$
\end{rmk}

\begin{ex}\label{ex:2.4}
By virtue of \cref{rmk:2.3}, it is not hard to verify that the following generator $g$ satisfies assumptions \ref{A:H1} and \ref{A:H2}: for $(\omega,t,y,z)\in \Omega\times \T\times \R^n\times \R^{n\times d}$,
$$
g^i(\omega,t,y,z)=(y^i)^2{\bf 1}_{y^i\leq 0}+\sin |y|+l(y)q(z^i)+|z^i|^2,\ \ i\in I_1
$$
and
$$
g^i(\omega,t,y,z)=-(y^i)^3{\bf 1}_{y^i\geq 0}+|y|+\cos |z^i|+|z^i|^{3\over 2}-|z^i|^2,\ \ i\in I_2,\vspace{0.1cm}
$$
where $I_1$ and $I_2$ are two nonempty sets such that $I_1+I_2=\{1,\cdots,n\}$, $l(\cdot)$ and $q(\cdot)$ are defined in (i) of \cref{rmk:2.3}. It is evident that $g$ is non-Lipschitz continuous with general growth in $y$, and $g^i$ is non-convex/non-concave in $z^i$ for each $i=1,\cdots,n$.\vspace{0.2cm}
\end{ex}

The second main result of the paper is \cref{thm:2.5} below, which addresses a stability result on solutions of multidimensional quadratic BSDEs with general unbounded terminal values and generators satisfying assumptions \ref{A:H1} and \ref{A:H2}.

\begin{thm}\label{thm:2.5}
Assume that $g$ satisfies assumptions \ref{A:H1} and \ref{A:H2}, and $(Y_t,Z_t)_{t\in\T}$ is the unique solution to BSDE \eqref{eq:1.1} such that $(Y,Z)\in \mathcal{E}(\R^n)\times \mcal(\R^{n\times d})$, and that for each $m\geq 1$, $\xi^{(m)}$ is an $\F_T$-measurable and $\R^n$-valued random variable, $g^{(m)}$ is a generator which satisfies \ref{A:H1} and \ref{A:H2} with an $(\F_t)$-progressively measurable nonnegative real-valued  process $\alpha^{(m)}$ instead of the process $\alpha$, and $(Y^{(m)}_t,Z^{(m)}_t)_{t\in\T}$ is the unique solution to BSDE \eqref{eq:1.1} associated with $(\xi^{(m)},g^{(m)})$ such that $(Y^{(m)},Z^{(m)})\in \mathcal{E}(\R^n)\times \mcal(\R^{n\times d})$. We also assume that for each $p\geq 1$,\vspace{0.1cm}
\begin{equation}\label{eq:2.12}
\E\left[\exp\left(p|\xi|+p\int_0^T \alpha_t {\rm d}t \right)\right]+\sup_{m\geq 1}\E\left[\exp\left(p|\xi^{(m)}|+p\int_0^T \alpha^{(m)}_t {\rm d}t \right)\right]<+\infty.\vspace{0.2cm}
\end{equation}
If $\ps$, $\lim_{m\To\infty}|\xi^{(m)}-\xi|=0$ and \vspace{0.1cm}
\begin{equation}\label{eq:2.13}
\lim\limits_{m\To\infty} \int_0^T \left| g^{(m)}(s,Y_s,Z_s)-g(s,Y_s,Z_s)\right| {\rm d}s=0,
\end{equation}
then for each $p\geq 1$, we have
\begin{equation}\label{eq:2.14}
\lim\limits_{m\To\infty} \E\left[\exp\left(p \sup_{t\in\T}|Y^{(m)}_t-Y_t|\right)
+\left(\int_0^T|Z^{(m)}_s-Z_s|^2{\rm d}s\right)^{p\over 2}
\right]=1.\vspace{0.4cm}
\end{equation}
\end{thm}

In the sequel, let $n\geq 2$ and $1\leq l\leq n-1$ be given. Denote by $y^{[1:l]}$ the $\R^{l}$-valued vector whose $i$-th component is just the $i$-th component of the vector $y\in\R^n$ for $i=1,\cdots,l$, and by $z^{[1:l]}$ the $\R^{l\times d}$-valued matrix whose $i$-th row is just the $i$-th row of the matrix $z\in \R^{n\times d}$ for $i=1,\cdots,l$. Furthermore, for any $U\in \R^n$ and $y\in \R^n$, denote by $U(y^{[1:l]})$ the vector whose first $l$-th components are the first $l$-th components of $y$, and whose last $(n-l)$-th components are the last $(n-l)$-th components of $U$. By relabelling, the following assumptions on the generator $g$ will be needed in the following theorem.

\begin{enumerate}
\renewcommand{\theenumi}{(H3)}
\renewcommand{\labelenumi}{\theenumi}
\item\label{A:H3} For $i=1,\cdots,l$, $f=g^i$ or $f(\omega,t,y,z)=-g^i(\omega,t,y(-y^i;i),-z)$ varies with $(\omega,t,y,z^{[1;l]})\in \Omega\times\T\times \R^n\times\R^{l\times d}$, and $\as$, for each $(y,z)\in \R^n\times\R^{l\times d}$,
\begin{equation}\label{eq:2.13-1}
{\rm sgn}(y^i)f(\omega,t,y,z)\leq \alpha_t(\omega)+\beta|y|+\gamma |z|
\end{equation}
and
\begin{equation}\label{eq:2.14-1}
\left|f(\omega,t,y,z)\right|\leq \alpha_t(\omega)+k\exp{(k|y|)}+\gamma |z|;\vspace{0.1cm}
\end{equation}
and furthermore, for $i=1,\cdots,l$, $\as$, it satisfies the following inequality
\begin{equation}\label{eq:2.15}
  \begin{array}{l}
   \Dis {\bf 1}_{\{y^i_1-\theta y^i_2>0\}}\left(f(\omega,t,y_1,z_1)-\theta f(\omega,t,y_2,z_2)\right)\vspace{0.2cm}\\
   \ \ \leq \Dis (1-\theta)\left(\alpha_t(\omega)+ \beta \left|\frac{y_1-\theta y_2}{1-\theta}\right|+\gamma \left|\frac{z_1-\theta z_2}{1-\theta}\right|+k(|y_1|+|y_2|)\right)\vspace{0.1cm}\\
   \Dis \hspace{1.5cm}\RE (y_1,y_2,z_1,z_2)\in \R^n\times\R^n\times\R^{l\times d}\times\R^{l\times d},\ \RE\theta\in (0,1).
   \end{array}
\end{equation}
In addition, for $i=l+1,\cdots,n$, $f=g^i$ or $f(\omega,t,y,z)=-g^i(\omega,t,y(-y^i;i),-z)$ varies only with $(\omega,t,y,z^i)\in \Omega\times [0,T]\times \R^n\times \R^{1\times d}$, and satisfies \eqref{eq:2.2}-\eqref{eq:2.3}, and \eqref{eq:2.4} with $\bar k=0$  or along with the additional inequality  \eqref{eq:2.5} in assumption \ref{A:H2}.\vspace{0.2cm}
\end{enumerate}


\begin{rmk}\label{rmk:2.6}
(i) By a similar argument to those in \cref{rmk:2.1}, in assumption~\ref{A:H3}, we may assume without loss of generality that $f=g^i$ satisfies \eqref{eq:2.13}-\eqref{eq:2.15} for all $i=1,\ldots,l$, and \eqref{eq:2.2}-\eqref{eq:2.3} as well as \eqref{eq:2.4} with $\bar k=0$  or along with the additional inequality  \eqref{eq:2.5} for all $i=l+1,\ldots,n$.\vspace{0.2cm}

(ii) For each $i=1,\cdots, l$, letting $y_1=y_2=y$ and $z_1=z_2=z$ in \eqref{eq:2.15} yields that
$$
{\bf 1}_{\{y^i>0\}} f(\omega,t,y,z)\leq \alpha_t(\omega)+(\beta+2k)|y|+\gamma |z|,\vspace{-0.1cm}
$$
which means that a part of the assertion stated in \eqref{eq:2.13} holds with different parameters. Furthermore, assume that for each $i=1,\ldots,l$, $f=g^i$ satisfies \eqref{eq:2.13}-\eqref{eq:2.15}. By sending $\theta\To 1$ in \eqref{eq:2.15} we can deduce that $\as$, for each $i=1,\cdots,l$ and $(y_1,y_2,z_1,z_2)\in \R^n\times \R^n\times \R^{n\times d}\times\R^{n\times d}$,
$$
{\bf 1}_{\{y_1^i-y_2^i>0\}} \left(g^i(\omega,t,y_1,z^{[1:l]}_1)-g^i(\omega,t,y_2,z^{[1:l]}_2)\right)\leq \beta |y_1-y_2|+\gamma |z^{[1:l]}_1-z^{[1:l]}_2|
$$
and then
$$
{\rm sgn}(y_1^i-y_2^i)\left(g^i(\omega,t,y_1,z^{[1:l]}_1)
-g^i(\omega,t,y_2,z^{[1:l]}_2)\right)\leq \beta |y_1-y_2|+\gamma |z^{[1:l]}_1-z^{[1:l]}_2|.
$$
Thus, for any fixed $(\F_t)$-progressively measurable $\R^n$-valued process $(U_t)_{t\in \T}$, by the last inequality we can directly verify that $\as$, for each $(y_1,y_2,z)\in \R^n\times \R^n\times \R^{n\times d}$,
\begin{equation}\label{eq:2.17}
\begin{array}{l}
\Dis\left\langle y^{[1:l]}_1-y^{[1:l]}_2, g^{[1:l]}(t,U_t(y^{[1:l]}_1),z^{[1:l]})-
g^{[1:l]}(t,U_t(y^{[1:l]}_2),z^{[1:l]})\right\rangle\vspace{0.1cm}\\
\ \ \Dis \leq \sum_{i=1}^{l} |y^i_1-y^i_2|\cdot \beta |y^{[1:l]}_1-y^{[1:l]}_2| \leq \beta\sqrt{l}|y^{[1:l]}_1-y^{[1:l]}_2|^2,
\end{array}
\end{equation}
and that for each $(y,z_1,z_2)\in \R^n\times \R^{n\times d}\times \R^{n\times d}$, with the help of the continuity of $g^i$ in $y$,
\begin{equation}\label{eq:2.18}
\left|g^{[1:l]}(t,U_t(y^{[1:l]}),z^{[1:l]}_1)-
g^{[1:l]}(t,U_t(y^{[1:l]}),z^{[1:l]}_2) \right|\leq l\gamma |z^{[1:l]}_1-z^{[1:l]}_2|.
\end{equation}
That is to say, $g^{[1:l]}(t,U_t(y^{[1:l]}),z^{[1:l]})$ for $(y,z)\in \R^n\times\R^{n\times d}$ satisfies the well-known monotonicity condition in $y^{[1:l]}$ and the classical Lipschitz continuity condition in $z^{[1:l]}$. In addition, it follows from \eqref{eq:2.13} and \eqref{eq:2.14} that $\as$, for each $(y,z)\in \R^n\times\R^{n\times d}$,
\begin{equation}\label{eq:2.19}
\begin{array}{l}
\Dis \left\langle \frac{y^{[1:l]}}{|y^{[1:l]}|}{\bf 1}_{|y^{[1:l]}|>0}, g^{[1:i]}(t,U_t(y^{[1:l]}),z^{[1:l]})\right\rangle\leq \frac{{\bf 1}_{|y^{[1:l]}|>0}}{|y^{[1:l]}|} \sum_{i=1}^{l} |y^i| {\rm sgn}(y^i)g^i(t,U_t(y^{[1:l]}),z^{[1:l]})\vspace{0.1cm}\\
\ \ \hspace{6.15cm}\Dis \leq \sqrt{l}\alpha_t(\omega)+\sqrt{l}\beta |U_t|+\sqrt{l}\beta |y^{[1:l]}|+\sqrt{l}\gamma |z^{[1:l]}|
\end{array}
\end{equation}
and
\begin{equation}\label{eq:2.20}
\left|g^{[1:l]}(t,U_t(y^{[1:l]}),z^{[1:l]})\right|\leq l\alpha_t+lk\exp{(2k|U_t|)}+lk\exp{(2k|y^{[1:l])}|}+l\gamma |z^{[1:l]}|.\vspace{0.2cm}
\end{equation}

(iii) In the same way as that in \citet[Proposition 3.5]{FanHu2021SPA}, it is not difficult to verify that for $i=1,\cdots,l$, \eqref{eq:2.15} holds for $f$ but with different parameters provided that $f$ satisfies \eqref{eq:2.13} and any one of the following three conditions:
\begin{enumerate}
\item [a)] Almost everywhere in $\Omega\times\T$, $f(\omega,t,\cdot,\cdot)$ is convex;
\item [b)] Almost everywhere in $\Omega\times\T$, for each $(y,z)\in \R^n\times\R^{l\times d}$, $f(\omega,t,\cdot,z)$ is Lipschitz continuous and $f(\omega,t, y,\cdot)$ is convex;
\item [c)] $f(\omega,t,y,z):= l(y)q(z)$ for $(y,z)\in \R^n\times\R^{l\times d}$, where both functions $l:\R^n\To\R$ and $q:\R^{l\times d}\To\R$ are bounded Lipschitz continuous, and the function $q(\cdot)$ has a bounded support.
\end{enumerate}
And, it is evident that inequality \eqref{eq:2.15} holds also true for the sum of finite numbers of generators $f$ satisfying \eqref{eq:2.15}, but with different parameters.\vspace{0.2cm}
\end{rmk}

\cref{thm:2.7} below constitutes the third main result of the paper. To the best of our knowledge, it seems to be the first result on existence and uniqueness of unbounded global solution of multidimensional quadratic BSDEs with general unbounded terminal values and with interactions in the second unknown variables among individual BSDEs.

\begin{thm}\label{thm:2.7}
Assume that \eqref{eq:2.1} holds and the generator $g$ satisfies assumptions \ref{A:H1} and \ref{A:H3}. Then BSDE \eqref{eq:1.1} admits a unique solution $(Y_t,Z_t)_{t\in\T}$ such that $(Y,Z)\in \mathcal{E}(\R^n)\times \mcal(\R^{n\times d})$.\vspace{0.2cm}
\end{thm}

The following corollary follows immediately from \cref{thm:2.7} and \cref{rmk:2.6}.

\begin{cor}\label{cor:2.8}
Assume that \eqref{eq:2.1} holds and the generator $g$ satisfies assumptions \ref{A:H1} and the following
\begin{enumerate}
\renewcommand{\theenumi}{(H3')}
\renewcommand{\labelenumi}{\theenumi}
\item\label{A:H3'} For $i=1,\cdots,l$, $g^i$ depends only on $(\omega,t,y,z^{[1;l]})$, and $\as$, $g^i$ is Lipschitz continuous in $y$ and convex with linear growth in $z^{[1;l]}$. And, for $i=l+1,\cdots,n$, $g^i$ depends only on $(\omega,t,y,z^i)$, and $\as$, $g^i$ is Lipschitz continuous in $y$ and convex in $z^i$. Furthermore, we assume that
 $$
\as,\ \ \left|g(\omega,t,y,z)\right|\leq \alpha_t(\omega)+\beta |y|+\frac{\gamma}{2} |z|^2,\ \ \RE (y,z)\in \R^n\times\R^{n\times d}.
$$
\end{enumerate}
Then BSDE \eqref{eq:1.1} admits a unique solution $(Y_t,Z_t)_{t\in\T}$ such that $(Y,Z)\in \mathcal{E}(\R^n)\times \mcal(\R^{n\times d})$.\vspace{0.2cm}
\end{cor}

\begin{ex}\label{ex:2.9}
By virtue of \cref{rmk:2.6}, it is not very hard to verify that the following generator $g$ satisfies assumptions \ref{A:H1} and \ref{A:H3}: for $(\omega,t,y,z)\in \Omega\times \T\times \R^n\times \R^{n\times d}$,
$$
g^i(\omega,t,y,z)=-\sqrt[3]{y^i}{\bf 1}_{y^i\leq 0}+ \sin |y|+l(y)q(z^{[1:l]})+|z^{[1:l]}|,\ \ 1\leq i\leq l
$$
and
$$
g^i(\omega,t,y,z)=(e^{-y^i}-1){\bf 1}_{y^i\geq 0}-|y|-|z^i|^{5\over 4}+|z^i|^2,\ \ l+1\leq i\leq n,\vspace{0.1cm}
$$
where $l(\cdot)$ and $q(\cdot)$ are defined in (iii) of \cref{rmk:2.6}. It is evident that $g$ is non-Lipschitz continuous with general growth in $y$ and non-convex (non-concave) with quadratic growth in $z$.\vspace{0.2cm}
\end{ex}

\begin{rmk}\label{rmk:2.10}
By combining the proofs of \cref{thm:2.5} and \cref{thm:2.7} it can be verified that \cref{thm:2.5} is also true when the assumption \ref{A:H2} is replaced with \ref{A:H3}.\vspace{0.2cm}
\end{rmk}

\subsection{Statement of main results: the case of special unbounded terminal values\vspace{0.2cm}}

In this subsection, we will present an existence and uniqueness result on the unbounded solution of a system of quadratic BSDEs with a special type of unbounded terminal value, which is expressed as the sum of a bounded random vector and a stochastic integral of a BMO process with respect to the Brownian motion. In this case, more general interactions among the second unknown variables of individual BSDEs are incorporated. Before that, we need to first establish existence and uniqueness on the bounded solution of a more general system of quadratic BSDEs with bounded terminal value.\vspace{0.2cm}

Let us start with introducing the following two assumptions \ref{A:H1'} and \ref{A:H4} on the generator $g$ defined on $\Omega\times [0,T]\times \R^n\times \R^{n\times d}$. Assumption \ref{A:H1'} is just \citet[assumption (B2)]{FanHuTang2025EJP}. It is evident that \ref{A:H1'} can imply \ref{A:H1}. Assumption \ref{A:H4} has a similar form to \cite[assumption (C1a)]{FanHuTang2025EJP}, but there also exists some crucial differences among them. In assumptions \ref{A:H1'} and \ref{A:H4}, it is assumed that $\phi(\cdot):[0,+\infty)\To (0,+\infty)$ is a deterministic nondecreasing function.

\begin{enumerate}
\renewcommand{\theenumi}{(H1')}
\renewcommand{\labelenumi}{\theenumi}
\item\label{A:H1'} For $i=1,\cdots,n$, $g^i$ satisfies that $\as$, for any $(y,\bar y, z, \bar z)\in \R^n\times\R^n\times\R^{n\times d}\times\R^{n\times d}$,
    $$
    \begin{array}{l}
   |g^i(\omega,t,y,z)-g^i(\omega,t,\bar y,\bar z)|\\
   \ \ \leq \Dis \phi(|y|\vee|\bar y|) \Dis \Big\{\left([v_t(\omega)]^{1+\delta}+|z|^{1+\delta}+|\bar z|^{1+\delta}\right)|y-\bar y|+\left(v_t(\omega)+|z|+|\bar z|\right)\sum_{j=1}^i |z^j-\bar z^j|\\
     \Dis \hspace{2.8cm} +\left[\left([v_t(\omega)]^\delta+|z|^\delta+|\bar z|^\delta\right)+\theta \left(v_t(\omega)+|z|+|\bar z|\right)\right]\sum_{j=i+1}^n |z^j-\bar z^j|\Big\}.
    \end{array}
    $$
\end{enumerate}

\begin{enumerate}
\renewcommand{\theenumi}{(H4)}
\renewcommand{\labelenumi}{\theenumi}
\item\label{A:H4} For each fixed $i=1,\cdots,n$, any one of the following three conditions holds:
\begin{enumerate}
\item [(i)] The random field $f:=g^i$ or $f(\omega,t,y,z):=-g^i(\omega,t,y(-y^i;i),z(-z^i;i))$ satisfies that $\as$, for any $(y,z)\in \R^n\times\R^{n\times d}$,
$$
\begin{array}{l}
\Dis \frac{\bar \gamma}{2} |z^i|^2-\bar\alpha_t(\omega)-\beta |y|-\sum_{j=i+1}^{n} \big(\bar\lambda |z^j|^{1+\delta}+\theta |z^j|^2\big)-\bar c\sum_{j=1}^{i-1}|z^j|^2\leq f(\omega,t,y,z)\\[2mm]
\quad\quad \leq  \Dis \alpha_t(\omega)+\left[\beta |y| \, {\bf 1}_{y^i>0}+\phi(|y|)\, {\bf 1}_{y^i\leq 0}\right]+\lambda \sum_{j\neq i} |z^j|^{\delta \, {\bf 1}_{y^i\leq 0}+1}+\frac{\gamma}{2} |z^i|^2;
\end{array}
$$

\item [(ii)] Almost everywhere in $\Omega\times\T$, for any $(y,z)\in \R^n\times\R^{n\times d}$, we have
$$
\begin{array}{l}
\Dis -\left[\beta |y^i|\,  {\bf 1}_{y^i<0}+\phi(|y|)\, {\bf 1}_{y^i\geq 0}\right]-l_i(\omega,t,y,z)\leq g^i(\omega,t,y,z)\\[3mm]
\hspace{1.5cm}\Dis \leq \left[\beta |y^i|\,  {\bf 1}_{y^i>0}+\phi(|y|)\, {\bf 1}_{y^i\leq 0}\right]+l_i(\omega,t,y,z)
\end{array}
$$
with
$$
l_i(\omega,t,y,z):=\tilde\alpha_t(\omega)+\frac{\gamma}{2} |z^i|^2+\Big\{v_t(\omega)+\beta \sqrt{|y|}+c\sum_{j=1}^{i-1}|z^j|+\sum_{j=i+1}^{n} \big(\lambda |z^j|^\delta+\theta |z^j|\big)\Big\}|z^i|;
$$

\item [(iii)] Almost everywhere in $\Omega\times\T$, for any $(y,z)\in \R^n\times\R^{n\times d}$, we have
$$
\begin{array}{l}
\Dis -\big[\beta (\,|y^i|+\sqrt{|y|}\,)\,{\bf 1}_{y^i<0}+\phi(|y|)\, {\bf 1}_{y^i\geq 0}\big]-\bar l_i(\omega,t,z)\leq g^i(\omega,t,y,z)\\[3mm]
\hspace{1.5cm}\Dis \leq \big[\beta (\,|y^i|+\sqrt{|y|}\,)\,{\bf 1}_{y^i>0}+\phi(|y|)\, {\bf 1}_{y^i\leq 0}\big]+\bar l_i(\omega,t,z)
\end{array}
$$
with
$$
\bar l_i(\omega,t,z):=\bar\alpha_t(\omega)+\bar\lambda |z|.\vspace{0.2cm}
$$
\end{enumerate}
\end{enumerate}

\cref{thm:2.11} below is the fourth main result of the paper, which establishes existence and uniqueness on the bounded global solution of a general system of quadratic BSDEs with bounded terminal value. Notably, \cref{thm:2.11} updates the assumptions of the three authors' previous \cite[Theorem 2.10]{FanHuTang2025EJP} so as to fill in a gap in the original proof. \footnote{The authors are greatly grateful to Professor Huilin Zhang from Shandong University for indicating a gap in the original proof of the authors' previous \cite[Theorem 2.10]{FanHuTang2025EJP} and for his helpful discussions.}

\begin{thm}\label{thm:2.11}
Let $\xi\in L^{\infty}(\R^n)$, $\alpha\in \ecal^{\infty}(p\gamma\exp(\beta T))$ for some real $p>1$ and the generator $g$ satisfy assumptions \ref{A:H1'} and \ref{A:H4}. Then, there exists a constant $\theta_0>0$ depending only on parameters $(\|\xi\|_{\infty}, \|\alpha\|_{\ecal^{\infty}(p\gamma\exp(\beta T))}, \|\bar\alpha\|_{\mcal^{\infty}}, \|\tilde\alpha\|_{\lcal^{\infty}}, \|v\|_{{\rm BMO}}, n,\beta,\gamma,\bar\gamma,\lambda,\bar\lambda, c,\bar c, \delta, T, p)$ such that when $\theta\in [0,\theta_0]$, BSDE \eqref{eq:1.1} has a unique solution $(Y,Z)\in \s^\infty(\R^n)\times {\rm BMO}(\R^{n\times d})$. Moreover, the above assertion remains true for $p=1$ when $\lambda=0$ and $\theta_0=0$.\vspace{0.2cm}
\end{thm}

\begin{ex}\label{ex:2.12}
Let us consider the following system of two BSDEs:
\begin{equation}\label{eq:2.4b}
\left\{
\begin{array}{l}
\Dis Y^1_t=\xi^1+\int_t^T \left(\theta_1|Z^1_s|^2+\vartheta_1|Z^2_s|^2\right){\rm d}s-\int_t^T Z^1_s{\rm d}B_s, \vspace{0.1cm}\\
\Dis Y^2_t=\xi^2+\int_t^T \left(\vartheta_2|Z^1_s|^2+\theta_2|Z^2_s|^2\right){\rm d}s-\int_t^T Z^2_s{\rm d}B_s,\quad t\in \T,
\end{array}
\right.
\end{equation}
where $\xi^i\in L^{\infty}(\R)$, $(\theta_i,\vartheta_i)^T\in \R^2$ and $\theta_i\vartheta_i<0$ for $i=1,2$. It is easy to verify that  $g^1(z):=\theta_1|z^1|^2+\vartheta_1|z^2|^2$ satisfies assumption \ref{A:H4}(i) with $\gamma=\bar\gamma=2|\theta_1|$ and $\theta=|\vartheta_1|$, and that $g^2(z):=\vartheta_2|z^1|^2+\theta_2|z^2|^2$ satisfies assumption \ref{A:H4}(i) with $\gamma=\bar\gamma=2|\theta_2|$ and $\bar c=|\vartheta_2|$. Then, according to \cref{thm:2.5}, we see that for each $\xi\in L^\infty(\R^2)$, there exists a positive constant $\vartheta_0>0$ depending only on parameters $(\theta_1,\theta_2,\vartheta_2,\|\xi\|_\infty)$ such that for each $\vartheta_1$ satisfying $|\vartheta_1|\leq \vartheta_0$, the system \eqref{eq:2.4b} of BSDEs admits a unique solution $(Y,Z)\in \s^\infty(\R^2)\times {\rm BMO}(\R^{2\times d})$. We note that this conclusion is different from that in \citet[Theorem 2.3]{JamneshanKupperLuo2017ECP}, where $\|\xi\|_{\infty}$ and $|\vartheta_2|$ are required to be small enough.\vspace{0.1cm}
\end{ex}

Let us further introduce the following two assumptions \ref{A:H1"} and \ref{A:H4'} on the generator $g$ defined on $\Omega\times [0,T]\times \R^n\times \R^{n\times d}$. It is evident that \ref{A:H1"} can imply \ref{A:H1'}, and \ref{A:H4'} can imply \ref{A:H4}. Notably, assumption \ref{A:H4'} has a similar form to \citet[assumption (D1)]{FanHuTang2025EJP}, but there also exists some crucial differences among them.

\begin{enumerate}
\renewcommand{\theenumi}{(H1")}
\renewcommand{\labelenumi}{\theenumi}
\item\label{A:H1"} For $i=1,\cdots,n$, $g^i$ satisfies that $\as$, for each $(y,\bar y, z, \bar z)\in \R^n\times\R^n\times\R^{n\times d}\times\R^{n\times d}$,
$$
\hspace*{-1cm}
\begin{array}{l}
\Dis |g^i(\omega,t,y,z)-g^i(\omega,t,\bar y,\bar z)|\leq \Dis \Dis  \gamma\left(v_t(\omega)+|z|+|\bar z|\right)\Big(|y-\bar y|+\sum_{j=1}^i |z^j-\bar z^j|\Big)\\
\Dis \hspace{2.6cm}+\left[\gamma \left([v_t(\omega)]^\delta+|z|^\delta+|\bar z|^\delta\right)+\theta \left(v_t(\omega)+|z|+|\bar z|\right)\right]\sum_{j=i+1}^n |z^j-\bar z^j|.
\end{array}
$$
\end{enumerate}

\begin{enumerate}
\renewcommand{\theenumi}{(H4')}
\renewcommand{\labelenumi}{\theenumi}
\item\label{A:H4'} For each fixed $i=1,\cdots,n$, any one of the following three conditions holds:
\begin{enumerate}
\item [(i)] The random field $f:=g^i$ or $f(\omega,t,y,z):=-g^i(\omega,t,y(-y^i;i),z(-z^i;i))$ satisfies that $\as$, for any $(y,z)\in \R^n\times\R^{n\times d}$,
$$
\hspace*{-1cm}\begin{array}{l}
\Dis \frac{\bar \gamma}{2} |z^i|^2-\bar\alpha_t(\omega)-\sum_{j=i+1}^{n} \big(\bar\lambda |z^j|^{1+\delta}+\theta |z^j|^2\big)-\bar c\sum_{j=1}^{i-1}|z^j|^2 \leq f(\omega,t,y,z)\\[5mm]
\quad\quad\quad \leq  \Dis \alpha_t(\omega)+\lambda\sum_{j\neq i}|z^j|+\frac{\gamma}{2} |z^i|^2;
\end{array}
$$

\item [(ii)] Almost everywhere in $\Omega\times\T$, for any $(y,z)\in \R^n\times\R^{n\times d}$, we have
$$
\hspace*{-0.5cm}|g^i(\omega,t,y,z)|\leq \tilde\alpha_t(\omega)+\frac{\gamma}{2} |z^i|^2+\Big\{v_t(\omega)
+c\sum_{j=1}^{i-1}|z^j|+\sum_{j=i+1}^{n} \big(\lambda |z^j|^\delta+\theta |z^j|\big)\Big\}|z^i|;
$$

\item [(iii)] Almost everywhere in $\Omega\times\T$, for any $(y,z)\in \R^n\times\R^{n\times d}$, we have
$$
\hspace*{-2cm}|g^i(\omega,t,y,z)|\leq \bar\alpha_t(\omega)+\bar\lambda |z|.\vspace{0.2cm}
$$
\end{enumerate}
\end{enumerate}

\cref{thm:2.13} below is the last main result of the paper, which establishes existence and uniqueness on the unbounded global solution of a general system of quadratic BSDEs with a special type of unbounded terminal value. Notably, \cref{thm:2.13} updates the assumptions of the three authors' previous \cite[Theorem 2.27]{FanHuTang2025EJP} so as to fill in a gap in the original proof.

\begin{thm}\label{thm:2.13}
Assume that $\bar\xi\in L^{\infty}(\R^n)$, $H\in {\rm BMO}(\R^{n\times d})$ and
$$\xi:=\bar\xi+\int_0^T H_s {\rm d}B_s,$$
and that the generator $g$ satisfy \ref{A:H1"} and \ref{A:H4'}. When $g^i$ satisfies \ref{A:H4'}(ii) and violates both \ref{A:H4'}(i) and \ref{A:H4'}(iii), and $-g^i$ violates \ref{A:H4'}(i), we further assume that $H^i\equiv 0$. If $\alpha\in \ecal^{\infty}(p\gamma)$ for some $p>1$, and $|H|^2\in \ecal^{\infty}(2\bar p (q\gamma)^2)$ for some $\bar p>1$ with $q=p/(p-1)$ such that $1/p+1/q=1$, then there exists a constant $\theta_0>0$ depending only on $(\|\bar\xi\|_{\infty}, \|\alpha\|_{\ecal^{\infty}(p\gamma)}, \|\bar\alpha\|_{\mcal^{\infty}}, \|\tilde\alpha\|_{\lcal^{\infty}}, \|v\|_{{\rm BMO}}, \||H|^2\|_{\ecal^{\infty}(2\bar p (q\gamma)^2)}, n,\gamma, \bar\gamma,\\ \lambda,\bar\lambda,c,\bar c,\delta, T, p,\bar p)$ such that when $\theta\in [0,\theta_0]$, BSDE \eqref{eq:2.1} has a unique solution $(Y_t,Z_t)_{t\in\T}$ such that
$$
\Big(Y-\int_0^\cdot H_s {\rm d}B_s, Z\Big)\in \s^\infty(\R^n)\times {\rm BMO}(\R^{n\times d}).\vspace{0.2cm}
$$
\end{thm}

\begin{rmk}\label{rmk:2.14}
It follow from the John-Nirenberg inequality that if $\|H\|_{\rm BMO}<1/(\sqrt{2\bar p}q\gamma)$, then $|H|^2\in \ecal^{\infty}(2\bar p (q\gamma)^2)$. Thus, according to \cref{thm:2.13}, we have a unique solution for the system of BSDEs, as stated in \citet[Theorem 3.1]{KramkovPulido2016AAP}.\vspace{0.2cm}
\end{rmk}

\begin{ex}\label{ex;2.15}
Let $A:=(a_{ij})_{n\times n}$ and $C:=(c_{ij})_{n\times n}$ be real-valued matrixes, and $h^i:\R^n\To \R$ and $\bar h^i:\R^{n\times d}\To \R$ be Lipschitz continuous functions for each $i=1,\cdots,n$. Assume that $a_{ij}\geq 0$ for each $i,j=1,\cdots,n$, and $h^i(\cdot)$ is bounded for each $i=1,\cdots,n$. Define generators $g:=(g^1,\cdots,g^n)^\top$ and $\bar g:=(\bar g^1,\cdots,\bar g^n)^\top$ as follows: for each $i=1,\cdots,n$, and $(\omega,t,y,z)\in \Omega\times \T\times \R^n\times\R^{n\times d}$,
$$
g^i(\omega,t,y,z):=\tilde\alpha_t(\omega)+h^i(y)+\bar h^i(z)+(-1)^i
\Big(\sum_{j=1}^{i-1}a_{ij}|z^j|^2-a_{ii}|z^i|^2
+\sum_{j=i+1}^{n}a_{ij}|z^j|^{1+\delta}\Big)
$$
and
$$
\bar g^i(\omega,t,y,z):=\tilde\alpha_t(\omega)+h^i(y)+z^i\sum_{j=1}^{i}c_{ij} (z^j)^\top+\left(v_t(\omega)+\sum_{j=i+1}^{n}c_{ij} |z^j|^\delta\right)|z^i|.\vspace{0.2cm}
$$
It is not very difficult to verify that both $f:=g$ and $f:=\bar g$ satisfy assumption \ref{A:H1"} with $\theta=0$, and that for each $i=1,\cdots,n$, $f:=g^i$ or $f(\omega,t,y,z):=-g^i(\omega,t,y(-y^i;i),z(-z^i;i))$ satisfies assumption \ref{A:H4'}(i) with $\theta=0$, and $f:=\bar g^i$ satisfies assumption \ref{A:H4'}(ii) with $\theta=0$.\vspace{0.2cm}
\end{ex}

\section{The proof of main results\vspace{0.2cm}}
\label{sec:3-The proof of main results}
\setcounter{equation}{0}

\subsection{The proof of \cref{thm:2.2}\vspace{0.2cm}}

Let assumptions \ref{A:H1} and \ref{A:H2} be satisfied by the generator $g$ and inequality \eqref{eq:2.1} hold. By \cref{rmk:2.1}, we assume without loss of generality that $f=g^i$ satisfies \eqref{eq:2.2}-\eqref{eq:2.3} for all $i=1,\ldots,n$, and that the integers $1,\cdots,n$ are divided into two parts $I_1$ and $I_2$ such that $f=g^i$ satisfies the inequality~\eqref{eq:2.4} with $\bar k=0$ for $i\in I_1$ and the inequalities~\eqref{eq:2.4}  and \eqref{eq:2.5} for $i\in I_2$, and $I_1+I_2=\{1,\cdots,n\}$.\vspace{0.2cm}

For a pair of processes $(U,V)\in \mathcal{E}(\R^n)\times \mcal (\R^{n\times d})$, consider the following decoupled system of quadratic BSDEs: for $i=1,\cdots, n$,\vspace{0.1cm}
\begin{equation}\label{eq:3.2}
  Y_t^i=\xi^i+\int_t^T g^i(s,U_s(Y_s^i;i),Z_s^i){\rm d}s-\int_t^T Z_s^i {\rm d}B_s, \ \ t\in\T.\vspace{0.1cm}
\end{equation}
In light of assumptions \ref{A:H1} and \ref{A:H2}, it is not difficult to verify that for each fixed $i=1,\cdots,n$, the generator $f^i(t,y,z):=g^i(t,U_t(y;i),z)$ satisfies assumptions \ref{A:A1} and \ref{A:A2} in Appendix with $\hat\alpha_t:=\alpha_t+(\beta+2k)|U_t|$ and $\varphi_t(r):=k\exp\{k(|U_t|+r)\}$. Furthermore, in light of \eqref{eq:2.1} and the fact that $U\in \mathcal{E}(\R^n)$, by H\"{o}lder's inequality we have, for each $q>1$,
$$
\E\left[\exp\left(q|\xi^i|+q\int_0^T \left[\alpha_t+(\beta+2k)|U_t|\right]{\rm d}t\right)\right]<+\infty,\ \ \ \ i=1,\cdots,n,\vspace{0.2cm}
$$
which means that \eqref{eq:A.2} holds. It then follows from \cref{lem:A.1} in Appendix that for each $i=1,\cdots,n$, the following BSDE
$$
Y_t^i=\xi^i+\int_t^T g^i(s,U_s(Y_s^i;i),Z_s^i){\rm d}s-\int_t^T Z_s^i {\rm d}B_s, \ \ t\in\T
$$
admits a unique adapted solution $(Y^i,Z^i)$ such that for each $q> 1$,
$$\E\left[\exp\left(q\sup_{t\in \T}|Y_t^i|\right)+\left(\int_0^T |Z_t^i|^2{\rm d}t\right)^{q\over 2}\right]<+\infty.$$
Consequently, by H\"{o}lder's inequality we conclude that multidimensional BSDE \eqref{eq:3.2} has a unique solution $(Y,Z)$ in the space of processes $\mathcal{E}(\R^n)\times \mcal (\R^{n\times d})$.\vspace{0.2cm}

Based on the above analysis, we set $(Y^{(0)},Z^{(0)})=(0,0)$ and define, recursively, the sequence of processes $\{(Y^{(m)},Z^{(m)})\}_{m=1}^\infty$ in the space of processes $\mathcal{E}(\R^n)\times \mcal (\R^{n\times d})$ by the unique adapted solution of multidimensional quadratic BSDEs: for $i=1,\cdots, n$,
\begin{equation}\label{eq:3.3}
  Y_t^{(m+1);i}=\xi^i+\int_t^T g^i(s,Y_s^{(m)}(Y_s^{(m+1);i};i),Z_s^{(m+1);i}){\rm d}s-\int_t^T Z_s^{(m+1);i}{\rm d}B_s, \ \ t\in\T,
\end{equation}
where and hereafter, we denote by $Y^{(m);i}$ and $Z^{(m);i}$, respectively, the $i$th component of $Y^{(m)}$ and the $i$th row of $Z^{(m)}$. In the sequel, we will show that $\{(Y^{(m)},Z^{(m)})\}_{m=1}^\infty$ is a Cauchy sequence in the space of processes $\s^q(\R^n)\times \hcal^q(\R^{n\times d})$ for each $q>1$, and then converges to a pair of processes $(Y,Z)$ in $\mathcal{E}(\R^n)\times \mcal(\R^{n\times d})$, which is the desired solution of BSDE \eqref{eq:1.1}.\vspace{0.2cm}

We first show that for each $q>1$,
\begin{equation}\label{eq:3.4}
\sup_{m\geq 0}\E\left[\exp\left(q \sup_{t\in \T}(\gamma e^{c_0 t}|Y_t^{(m)}|)\right)\right]\leq C(q)
\end{equation}
with $c_0:=8n(\beta+\gamma)e^{\beta T}$ and
\begin{equation}\label{eq:3.5}
C(q):=\left(\frac{q}{q-1}\right)^{2q}\E\left[\exp\left\{qc_0 e^{c_0 T} \left(|\xi|+\int_0^T \alpha_t{\rm d}t\right)\right\}\right]<+\infty,
\end{equation}
and that
\begin{equation}\label{eq:3.6}
\sup_{m\geq 0}\E\left[\left(\int_0^T |Z_s^{(m)}|^2 {\rm d}s\right)^{q\over 2}\right]<+\infty. \vspace{0.2cm}
\end{equation}
In fact, it follows from \eqref{eq:2.2} that for each $i=1,\cdots,n$ and $m\geq 0$, $\ass$,
\begin{equation}\label{eq:3.7}
{\rm sgn}(Y^{(m+1);i}_s) g^i(s,Y^{(m)}_s(Y^{(m+1);i}_s;i),Z^{(m+1);i}_s)\leq \alpha_s+\beta |Y^{(m)}_s|+\beta |Y^{(m+1);i}_s|+\frac{\gamma}{2} |Z^{(m+1);i}_s|^2.
\end{equation}
Furthermore, in light of \eqref{eq:2.1} and the fact that $Y^{(m)}\in \mathcal{E}(\R^n)$, by H\"{o}lder's inequality we have
\begin{equation}\label{eq:3.8}
\RE\ q>1,\ \  \E\left[\exp\left\{q\left(|\xi^i|+\int_0^T \left(\alpha_s+\beta |Y_s^{(m)}|\right){\rm d}s\right)\right\}\right]<+\infty.
\end{equation}
Then, applying (i) of \cref{lem:A.3} in Appendix yields that for $m\geq 0$, $i=1,\cdots,n$, and $q>1$,
$$
\exp\left(\gamma |Y_t^{(m+1);i}|\right)\leq \E_t\left[\exp\left(\gamma e^{\beta T}|\xi^i|+\gamma e^{\beta T} \int_t^T \left(\alpha_s+\beta |Y_s^{(m)}|\right) {\rm d}s\right)\right],\ \ t\in\T
$$
and
$$
\E\left[\left(\int_0^T |Z_s^{(m+1);i}|^2 {\rm d}s\right)^{q\over 2}\right]\leq K_q \left\{\E\left[\left(\int_0^T \left(\alpha_s+\beta |Y_s^{(m)}| \right){\rm d}s\right)^q+\exp\left(q\gamma \sup_{s\in\T}|Y_s^{(m+1);i}|\right)\right]\right\},\vspace{0.2cm}
$$
where $K_q$ is a positive constant depending only on $(q,\beta,\gamma,T)$.
Consequently, taking $e^{c_0 t}$-order power in both sides of the last second inequality and noticing that
$$
|Y_t^{(m+1)}|\leq \sum_{i=1}^{n}|Y_t^{(m+1);i}|\ \ \ {\rm and}\ \ \ |Z_t^{(m+1)}|^2= \sum_{i=1}^{n}|Z_t^{(m+1);i}|^2,
$$
by virtue of Jensen's inequality and the basic inequality $(\sum_{i=1}^{n} |a_i|)^{q/2}\leq n^q \sum_{i=1}^{n} |a_i|^{q/2}$ we can obtain that for each $m\geq 0$ and $q>1$, in light of the definition of the constant $c_0$,
\begin{equation}\label{eq:3.9}
\begin{array}{l}
\Dis \exp\left(\gamma e^{c_0 t}|Y_t^{(m+1)}|\right)\vspace{0.1cm}\\
\ \ \leq \Dis \E_t\left[\exp\left(n\gamma e^{\beta T}e^{c_0 T} |\xi|+n\gamma e^{\beta T}e^{c_0 t} \int_t^T \left(\alpha_s+\beta |Y_s^{(m)}|\right) {\rm d}s\right)\right]\vspace{0.1cm}\\
\ \ \leq \Dis \E_t\left[\exp\left(\frac{c_0 e^{c_0 T}}{2}\left(|\xi|+\int_0^T \alpha_s {\rm d}s\right)+ \frac{c_0 e^{c_0 t}}{2} \int_t^T \left(e^{-c_0 s}\cdot \gamma e^{c_0 s}|Y_s^{(m)}|\right) {\rm d}s\right)\right],\vspace{0.3cm}\\
\ \ \leq \Dis \E_t\left[\exp\left(\frac{c_0 e^{c_0 T}}{2}\left(|\xi|+\int_0^T \alpha_s {\rm d}s\right)+ \frac{1}{2} \sup_{s\in\T} \left(\gamma e^{c_0 s}|Y_s^{(m)}|\right)\right)\right],\ \ \ \ t\in\T
\end{array}
\end{equation}
and
\begin{equation}\label{eq:3.10}
\begin{array}{l}
\Dis \E\left[\left(\int_0^T |Z_s^{(m+1)}|^2 {\rm d}s\right)^{q\over 2}\right]\vspace{0.1cm}\\
\ \ \Dis \leq K_q n^{\color{red}q+1} \E\left[\left(\int_0^T \left(\alpha_s+\beta |Y_s^{(m)}| \right){\rm d}s\right)^q+\exp\left(q\gamma \sup_{s\in\T}|Y_s^{(m+1)}|\right)\right].
\end{array}\vspace{0.1cm}
\end{equation}
In light of \eqref{eq:3.9}, using Doob's maximal inequality for martingales and H\"{o}lder's inequality we can deduce that for each $m\geq 0$ and $q>1$,
\begin{equation}\label{eq:3.11}
\begin{array}{l}
\Dis \E\left[\exp\left(q\sup_{s\in [0,T]}(\gamma e^{c_0 s} |Y_s^{(m+1)}|)\right)\right]\vspace{0.2cm}\\
\ \ \leq \Dis \left(\frac{q}{q-1}\right)^{q}\E\left[\exp\left\{\frac{q c_0 e^{c_0 T}}{2} \left(|\xi|+\int_0^T \alpha_t{\rm d}t\right)+\frac{q}{2}\sup_{s\in [0,T]}\left(\gamma e^{c_0 s} |Y_s^{(m)}|\right)\right\}\right]\vspace{0.2cm}\\
\ \ \leq \Dis (C(q))^{1\over 2}\left(\E\left[\exp\left(q \sup_{s\in [0,T]}\left(\gamma e^{c_0 s} |Y_s^{(m)}|\right)\right)\right]\right)^{1\over 2},\vspace{0.1cm}
\end{array}
\end{equation}
where $C(q)$ is defined in \eqref{eq:3.5}. Thus, in light of \eqref{eq:3.11}, by induction we can get
$$
\begin{array}{l}
\Dis \E\left[\exp\left(q\sup_{s\in [0,T]}\left(\gamma e^{c_0 s} |Y_s^{(m+1)}|\right)\right)\right]\vspace{0.1cm}\\
\ \ \Dis \leq \sqrt{C(q)}^{1+{1\over 2}+\cdots +{1\over 2^m}} \left(\E\left[\exp\left(q\sup_{s\in [0,T]}\left(\gamma e^{c_0 s} |Y_s^{(0)}|\right)\right)\right]\right)^{1\over 2^{m+1}}\leq C(q),
\end{array}\vspace{0.1cm}
$$
from which inequality \eqref{eq:3.4} follows immediately. And, \eqref{eq:3.6} follows from \eqref{eq:3.4} and \eqref{eq:3.10}.\vspace{0.2cm}

Next, in light of \ref{A:H1}-\ref{A:H2}, it is not difficult to verify that for each fixed $i=1,\cdots,n$ and $m\geq 0$, the generator $f^i(t,y,z):=g^i(t,Y_t^{(m)}(y;i),z)$ of BSDE \eqref{eq:3.3} satisfies \ref{A:A1}-\ref{A:A2} in Appendix with $\hat\alpha_t:=\alpha_t+(\beta+2k) |Y_t^{(m)}|$ and $\varphi_t(r):=k\exp\{k(|Y_t^{(m)}|+r)\}$. Then, note that$f^i$ satisfies \eqref{eq:A.6} for $i\in I_2$ due to the assumption that $g^i$ satisfies \eqref{eq:2.5} for $i\in I_2$. Since $(Y^{(m+1);i},Z^{(m+1);i})$ is the unique solution of BSDE \eqref{eq:3.3}, from \eqref{eq:A.11} in \cref{lem:A.1} we can conclude that for each $m\geq 0$, $i\in I_2$ and $\eps\in (0,\bar\gamma/9]$,
\begin{equation}\label{eq:3.12}
\begin{array}{l}
\Dis \E\left[\exp\left(\frac{\bar \gamma}{2}\eps\int_0^T |Z_t^{(m+1);i}|^2{\rm d}t \right)\right]\vspace{0.1cm}\\
\ \ \leq \Dis \E\left[\exp\left(3\eps(2+\beta T)\sup_{t\in [0,T]}|Y_t^{(m+1);i}|+3\eps\int_0^T \left(\alpha_t+(\beta+2k)|Y_t^{(m)}| \right){\rm d}t\right)\right]\vspace{0.2cm}\\
\ \ \leq \Dis \E\left[\exp\left(6\eps(1+\beta T+kT)\sup_{t\in [0,T]}(|Y_t^{(m+1)}|+|Y_t^{(m)}|)+3\eps\int_0^T \alpha_t {\rm d}t\right)\right].
\end{array}\vspace{0.2cm}
\end{equation}

Now, we define, for each fixed $m,p\geq 1$ and $\theta\in (0,1)$,
$$
\delta_{\theta}Y^{(m,p)}:=\frac{Y^{(m+p)}-\theta Y^{(m)}}{1-\theta}\ \ \ {\rm and}\ \ \ \delta_{\theta}Z^{(m,p)}:=\frac{Z^{(m+p)}-\theta Z^{(m)}}{1-\theta}.
$$
Then, in light of \eqref{eq:3.3}, $(\delta_{\theta}Y^{(m,p)},\delta_{\theta}Z^{(m,p)})\in \mathcal{E}(\R^n)\times \mcal(\R^{n\times d})$ solves the following system of quadratic BSDEs: for $i=1,\cdots, n$,\vspace{0.1cm}
\begin{equation}\label{eq:3.13}
  \delta_{\theta} Y_t^{(m,p);i}=\xi^i+\int_t^T \delta_{\theta} g^{(m,p);i}(s){\rm d}s-\int_t^T \delta_{\theta} Z_s^{(m,p);i}{\rm d}B_s, \ \ t\in\T,\vspace{0.1cm}
\end{equation}
where $\ass$,
\begin{equation}\label{eq:3.14}
\Dis\delta_{\theta} g^{(m,p);i}(s)
:=\Dis\frac{g^{i}\left(s, Y_s^{(m+p-1)}(Y_s^{(m+p);i};i), Z_s^{(m+p);i}\right)-\theta g^{i}\left(s,Y_s^{(m-1)}(Y_s^{(m);i};i),Z_s^{(m);i}\right)}{1-\theta}.
\end{equation}
Since \eqref{eq:2.4} holds with $\bar k=0$  for $i\in I_1$ and \eqref{eq:2.4} and \eqref{eq:2.5} hold for $i\in I_2$, by the above definition we can deduce that $\ass$, for each $i=1,\cdots,n$,
\begin{equation}\label{eq:3.15}
\begin{array}{l}
\Dis {\bf 1}_{\{\delta_{\theta} Y_s^{(m,p);i}>0\}}\delta_{\theta} g^{(m,p);i}(s)\leq \Dis \alpha_s+\beta |\delta_{\theta}Y_s^{(m-1,p)}(\delta_{\theta} Y_s^{(m,p);i};i)|+\frac{\gamma}{2}|\delta_{\theta} Z_s^{(m,p);i}|^2\vspace{0.1cm}\\
\hspace{4.6cm}\Dis +k\left(|Y_s^{(m-1)}(Y_s^{(m);i};i)|+
|Y_s^{(m+p-1)}(Y_s^{(m+p);i};i)|\right)\vspace{0.1cm}\\
\hspace{4.6cm}\Dis +\bar k{\bf 1}_{\{i\in I_2\}}\left(|Z_s^{(m);i}|^{1+\delta}+ |Z_s^{(m+p);i}|^{1+\delta}\right)\vspace{0.1cm}\\
\hspace{4.05cm}\leq \Dis \alpha_s+\beta |\delta_{\theta}Y_s^{(m-1,p)}|+H^{(m,p)}(s)+\beta |\delta_{\theta} Y_s^{(m,p);i}|+\frac{\gamma}{2}|\delta_{\theta} Z_s^{(m,p);i}|^2
\end{array}
\end{equation}
with
$$
H^{(m,p)}(s):=k\left(|Y_s^{(m-1)}|+|Y_s^{(m)}|+|Y_s^{(m+p-1)}|+|Y_s^{(m+p)}|\right)
+\bar k \sum_{i\in I_2}\left(|Z_s^{(m);i}|^{1+\delta}+|Z_s^{(m+p);i}|^{1+\delta}\right).
$$
Furthermore, by virtue of \eqref{eq:A.27}, \eqref{eq:2.1}, \eqref{eq:3.4} and \eqref{eq:3.12} along with H\"{o}lder's inequality we can obtain the existence of a constant $K(q)>0$ depending only on $(q,n,\beta,\gamma,\bar\gamma,\delta,k,\bar k,T)$ such that for each $q>1$,
\begin{equation}\label{eq:3.16}
\begin{array}{l}
\Dis\left(\frac{q}{q-1}\right)^{2q}\sup_{m,p\geq 1}\E\left[\exp\left\{qc_0 e^{c_0 T} \left(|\xi|+\sup_{s\in\T}\left(|Y_s^{(m)}|+|Y_s^{(m+p)}|\right)
\right.\right.\right.\vspace{0.1cm}\\
\hspace{4.55cm} \Dis +\left.\left.\left.\int_0^T \left(\alpha_s+H^{(m,p)}(s)\right) {\rm d}s\right)\right\}\right]\leq K(q).
\end{array}
\end{equation}
Consequently, in light of \eqref{eq:3.13}, \eqref{eq:3.15} and \eqref{eq:3.16}, it follows from (ii) of \cref{lem:A.3} with $p=1$, $\bar\alpha_t\equiv 0$ and $\hat\alpha_t=\alpha_t+\beta |\delta_{\theta}Y_t^{(m-1,p)}|+H^{(m,p)}(t)$ that for $i=1,\cdots,n$ and $t\in\T$,
\begin{equation}\label{eq:3.17}
\begin{array}{l}
\Dis \exp\left(\gamma \left(\delta_{\theta} Y_t^{(m,p);i}\right)^+\right)\vspace{0.1cm}\\
\ \ \leq \Dis \E_t\left[\exp\left(\gamma e^{\beta T} (\xi^i)^+ +\gamma e^{\beta T} \int_t^T \left(\alpha_s+\beta |\delta_{\theta}Y_s^{(m-1,p)}|+H^{(m,p)}(s)\right) {\rm d}s\right)\right].\vspace{0.2cm}
\end{array}
\end{equation}

On the other hand, define
$$
\delta_{\theta}\widetilde Y^{(m,p)}:=\frac{Y^{(m)}-\theta Y^{(m+p)}}{1-\theta}\ \ \ {\rm and}\ \ \ \delta_{\theta}\widetilde Z^{(m,p)}:=\frac{Z^{(m)}-\theta Z^{(m+p)}}{1-\theta}.
$$
The same computation as the last paragraph leads to that for $i=1,\cdots,n$,
\begin{equation}\label{eq:3.18}
\begin{array}{l}
\Dis \exp\left(\gamma \left(\delta_{\theta} \widetilde Y_t^{(m,p);i}\right)^+\right)\vspace{0.1cm}\\
\ \ \leq \Dis \E_t\left[\exp\left(\gamma e^{\beta T} (\xi^i)^+ +\gamma e^{\beta T} \int_t^T \left(\alpha_s+\beta |\delta_{\theta}\widetilde Y_s^{(m-1,p)}|+H^{(m,p)}(s)\right) {\rm d}s\right)\right],\ \ t\in\T.
\end{array}
\end{equation}
Furthermore, observe that for $i=1,\cdots,n$ and $t\in\T$, we have\vspace{0.1cm}
$$
\begin{array}{l}
\Dis \left(\delta_{\theta} Y_t^{(m,p);i}\right)^- = \Dis \frac{\left(Y_t^{(m+p);i}-\theta Y_t^{(m);i}\right)^-} {1-\theta}= \Dis \frac{\left(\theta Y_t^{(m);i}-Y_t^{(m+p);i}\right)^+} {1-\theta} \vspace{0.1cm}\\
\hspace{2.25cm}\leq \Dis \frac{\theta\left( Y_t^{(m);i}-\theta Y_t^{(m+p);i}\right)^+ +(1-\theta^2)|Y_t^{(m+p);i}|} {1-\theta}\vspace{0.1cm}\\
\hspace{2.25cm}\leq \Dis \left(\delta_{\theta} \widetilde Y_t^{(m,p);i}\right)^+ +2|Y_t^{(m+p)}|
\end{array}
$$
and, similarly,
$$
\left(\delta_{\theta} \widetilde Y_t^{(m,p);i}\right)^- \leq \left(\delta_{\theta} Y_t^{(m,p);i}\right)^+ +2|Y_t^{(m)}|.\vspace{0.1cm}
$$
It follows from \eqref{eq:3.17} and \eqref{eq:3.18} along with Jensen's inequality that for each $i=1,\cdots, n$,\vspace{0.1cm}
$$
\begin{array}{l}
\Dis \exp\left(\gamma |\delta_{\theta} Y_t^{(m,p);i}|\right) = \exp\left(\gamma \left(\delta_{\theta} Y_t^{(m,p);i}\right)^+\right)\cdot \exp\left(\gamma \left(\delta_{\theta} Y_t^{(m,p);i}\right)^-\right) \vspace{0.1cm}\\
\ \ \leq \Dis
\E_t\left[\exp\left(2\gamma e^{\beta T} |\xi|+2\gamma |Y_t^{(m+p)}| +2\gamma e^{\beta T} \int_0^T \left(\alpha_s+H^{(m,p)}(s)\right) {\rm d}s \right.\right. \vspace{0.1cm}\\
\ \ \Dis \hspace{1.8cm}\left.\left. +2\gamma\beta e^{\beta T} \int_t^T \left( |\delta_{\theta}Y_s^{(m-1,p)}|+|\delta_{\theta}\widetilde Y_s^{(m-1,p)}| \right) {\rm d}s\right)\right],\ \ t\in\T
\end{array}
$$
and
$$
\begin{array}{l}
\Dis \exp\left(\gamma |\delta_{\theta} \widetilde Y_t^{(m,p);i}|\right) = \exp\left(\gamma \left(\delta_{\theta} \widetilde Y_t^{(m,p);i}\right)^+\right)\cdot \exp\left(\gamma \left(\delta_{\theta} \widetilde Y_t^{(m,p);i}\right)^-\right) \vspace{0.1cm}\\
\ \ \leq \Dis
\E_t\left[\exp\left(2\gamma e^{\beta T} |\xi|+2\gamma |Y_t^{(m)}| +2\gamma e^{\beta T} \int_0^T \left(\alpha_s+H^{(m,p)}(s)\right) {\rm d}s \right.\right. \vspace{0.1cm}\\
\ \ \Dis \hspace{1.8cm}\left.\left. +2\gamma\beta e^{\beta T} \int_t^T \left( |\delta_{\theta}Y_s^{(m-1,p)}|+|\delta_{\theta}\widetilde Y_s^{(m-1,p)}| \right) {\rm d}s\right)\right],\ \ t\in\T.\vspace{0.1cm}
\end{array}
$$
Consequently, by Jensen's inequality again we have, for each $t\in \T$, \vspace{0.1cm}
$$
\begin{array}{l}
\Dis \exp\left(\gamma |\delta_{\theta} Y_t^{(m,p);i}|+\gamma |\delta_{\theta} \widetilde Y_t^{(m,p);i}| \right) \vspace{0.1cm}\\
\ \ \leq \Dis
\E_t\left[\exp\left(4\gamma e^{\beta T} \left(|\xi|+|Y_t^{(m)}| + |Y_t^{(m+p)}|+\int_0^T \left(\alpha_s+H^{(m,p)}(s)\right) {\rm d}s\right) \right.\right. \vspace{0.1cm}\\
\ \ \Dis \hspace{1.8cm}\left.\left. +4\gamma\beta e^{\beta T} \int_t^T \left( |\delta_{\theta}Y_s^{(m-1,p)}|+|\delta_{\theta}\widetilde Y_s^{(m-1,p)}| \right) {\rm d}s\right)\right],\ \ \ \ i=1,\cdots, n,
\end{array}
$$
and then, in light of the definition of the constant $c_0$,
$$
\begin{array}{l}
\Dis \exp\left(\gamma |\delta_{\theta} Y_t^{(m,p)}|+\gamma |\delta_{\theta} \widetilde Y_t^{(m,p)}|\right) \vspace{0.1cm}\\
\ \ \leq \Dis
\E_t\left[\exp\left\{\frac{c_0}{2}\left(|\xi|+|Y_t^{(m)}| + |Y_t^{(m+p)}|+\int_0^T \left(\alpha_s+H^{(m,p)}(s)\right) {\rm d}s\right) \right.\right. \vspace{0.1cm}\\
\ \ \Dis \hspace{1.8cm}\left.\left. +\frac{c_0\gamma}{2} \int_t^T \left( |\delta_{\theta}Y_s^{(m-1,p)}|+|\delta_{\theta}\widetilde Y_s^{(m-1,p)}| \right) {\rm d}s\right\}\right],
\end{array}
$$
which further leads to, similar to \eqref{eq:3.9},
\begin{equation}\label{eq:3.19}
\begin{array}{l}
\Dis \exp\left(\gamma e^{c_0 t}\left(|\delta_{\theta} Y_t^{(m,p)}|+ |\delta_{\theta} \widetilde Y_t^{(m,p)}|\right)\right) \vspace{0.1cm}\\
\ \ \leq \Dis
\E_t\left[\exp\left\{\frac{c_0e^{c_0 T}}{2}\left(|\xi|+\sup_{s\in \T}(|Y_s^{(m)}| + |Y_s^{(m+p)}|)+\int_0^T \left(\alpha_s+H^{(m,p)}(s)\right) {\rm d}s\right) \right.\right. \vspace{0.1cm}\\
\ \ \Dis \hspace{1.8cm}\left.\left. +\frac{1}{2}\sup_{s\in \T}\left[\gamma e^{c_0 s}\left(|\delta_{\theta}Y_s^{(m-1,p)}|+|\delta_{\theta}\widetilde Y_s^{(m-1,p)}|\right)\right]\right\}\right].\vspace{0.2cm}
\end{array}
\end{equation}

Next, in light of \eqref{eq:3.19}, \eqref{eq:3.16}, Doob's maximal inequality for martingales along with H\"{o}lder's inequality yields that for each $q>1$,
\begin{equation}\label{eq:3.20}
\hspace*{-0.2cm}\begin{array}{l}
\Dis \E\left[\exp\left\{q \sup_{t\in [0,T]}\left[\gamma e^{c_0 t}\left(|\delta_{\theta} Y_t^{(m,p)}|+|\delta_{\theta} \widetilde Y_t^{(m,p)}|\right)\right]\right\}\right]\vspace{0.2cm}\\
\ \ \leq \Dis \left(\frac{q}{q-1}\right)^{q}\E\left[\exp\left\{\frac{qc_0e^{c_0 T}}{2} \left(|\xi|+\sup_{s\in \T}(|Y_s^{(m)}|+|Y_s^{(m+p)}|)+\int_0^T \left(\alpha_s+H^{(m,p)}(s)\right) {\rm d}s\right)\right.\right.\vspace{0.1cm}\\
\ \ \Dis \hspace{2.5cm} \left.\left.+\frac{q}{2}\sup_{s\in \T}\left[\gamma e^{c_0 s}\left(|\delta_{\theta}Y_s^{(m-1,p)}|+|\delta_{\theta}\widetilde Y_s^{(m-1,p)}|\right)\right]\right\}\right]\vspace{0.2cm}\\
\ \ \leq \Dis (K(q))^{1\over 2}\left(\E\left[\exp\left\{q\sup_{t\in \T}\left[\gamma e^{c_0 t}\left(|\delta_{\theta}Y_t^{(m-1,p)}|+|\delta_{\theta}\widetilde Y_t^{(m-1,p)}| \right)\right]\right\}\right]\right)^{1\over 2}. \vspace{0.1cm}
\end{array}
\end{equation}
Based on the last inequality, by induction we can derive that
$$
\begin{array}{l}
\Dis \E\left[\exp\left\{q \sup_{t\in [0,T]}\left[\gamma e^{c_0 t}\left(|\delta_{\theta} Y_t^{(m,p)}|+|\delta_{\theta} \widetilde Y_t^{(m,p)}|\right)\right]\right\}\right]\vspace{0.1cm}\\
\ \ \leq \Dis K(q) \left(\E\left[\exp\left\{q \sup_{t\in [0,T]}\left[\gamma e^{c_0 t}\left(|\delta_{\theta} Y_t^{(1,p)}|+|\delta_{\theta} \widetilde Y_t^{(1,p)}|\right)\right]\right\}\right]\right)^{\frac{1}{2^m}},
\end{array}
$$
from which along with \eqref{eq:3.4} it follows that for each $q>1$ and $\theta\in (0,1)$,
\begin{equation}\label{eq:3.21}
\limsup_{m\To\infty}\sup_{p\geq 1}\E\left[\exp\left(q\gamma \sup_{t\in \T}\left(|\delta_{\theta} Y_t^{(m,p)}|+|\delta_{\theta} \widetilde Y_t^{(m,p)}|\right)\right)\right]\leq K(q).
\end{equation}
Thus, in light of $x\leq e^x$ for $x\geq 0$, we have for each $\theta\in (0,1)$,
$$
\limsup_{m\To\infty}\sup_{p\geq 1}\E\left[\sup_{t\in \T}|Y_t^{(m+p)}-\theta Y_t^{(m)}|\right]\leq (1-\theta)\frac{K(2)}{2\gamma},
$$
and then, in light of \eqref{eq:3.4},
$$
\limsup_{m\To\infty}\sup_{p\geq 1}\E\left[\sup_{t\in \T}|Y_t^{(m+p)}- Y_t^{(m)}|\right]\leq (1-\theta)\left(\frac{K(2)}{2\gamma}+\sup_{m\geq 1}\E\left[\sup_{t\in \T}|Y_t^{(m)}|\right]\right)<+\infty.
$$
Sending $\theta\To 1$ in the last inequality and using \eqref{eq:3.4} again yields the existence of a process $Y\in \mathcal{E}(\R^n)$ such that for each $q>1$,
\begin{equation}\label{eq:3.22}
\lim_{m\To\infty} \E\left[\sup_{t\in \T}|Y_t^{(m)}-Y_t|^q\right]=0
\end{equation}
and then
\begin{equation}\label{eq:3.23}
\lim_{m\To\infty} \E\left[\exp\left(q\sup_{t\in \T}|Y_t^{(m)}-Y_t|\right)\right]=1.\vspace{0.2cm}
\end{equation}
Furthermore, it follows from It\^{o}'s formula that for each $m,p\geq 1$,
\begin{equation}\label{eq:3.24}
\E\left[\int_0^T |Z^{(m+p)}_s-Z^{(m)}_s|^2 {\rm d}s\right]\leq 2\E\left[\sup_{t\in\T}\left|Y^{(m+p)}_t-Y^{(m)}_t\right| \int_0^T \Delta^{m,p}_s\ {\rm d}s \right]
\end{equation}
with
$$
\Delta^{m,p}_s:=\sum_{i=1}^n\left|g^i(s,Y^{(m+p-1)}_s(Y^{(m+p);i}_s;i), Z^{(m+p);i}_s)-g^i(s,Y^{(m-1)}_s(Y^{(m);i}_s;i), Z^{(m);i}_s)\right|.\vspace{0.1cm}
$$
And, by virtue of \eqref{eq:2.3}, \eqref{eq:3.4} and \eqref{eq:3.6} we can deduce that
\begin{equation}\label{eq:3.25}
\sup_{m,p\geq 1}\E\left[\left(\int_0^T \Delta^{m,p}_s\ {\rm d}s \right)^2\right]<+\infty.
\end{equation}
Then, applying H\"{o}lder's inequality to \eqref{eq:3.24} and using \eqref{eq:3.22} and \eqref{eq:3.25} yields that
$$
\lim_{m\To\infty}\sup_{p\geq 1}\E\left[\int_0^T |Z_s^{(m+p)}-Z_s^{(m)}|^2{\rm d}s\right]=0,
$$
from which along with \eqref{eq:3.6} we obtain the existence of a process $Z\in \mcal(\R^{n\times d})$ such that
\begin{equation}\label{eq:3.26}
\RE\ q>1,\ \ \lim_{m\To\infty} \E\left[\left(\int_0^T |Z_s^{(m)}-Z_s|^2{\rm d}s\right)^{q\over 2}\right]=0.
\end{equation}
Finally, in light of equalities \eqref{eq:3.22}, \eqref{eq:3.23} and \eqref{eq:3.26}, sending $m\To \infty$ in BSDE \eqref{eq:3.3} yields that $(Y,Z)$ is a desired solution of BSDE \eqref{eq:1.1}. And, in light of assumptions \ref{A:H1}-\ref{A:H2}, it is easy to verify that for each fixed $i=1,\cdots,n$, the generator $f^i(t,y,z):=g^i(t, Y_t(y;i),z)$ satisfies assumptions \ref{A:A1}-\ref{A:A2} in Appendix with $\hat\alpha_t:=\alpha_t+(\beta+2k)|Y_t|$ and $\varphi_t(r):=k\exp\{k(|Y_t|+r)\}$. Then, note that $f^i$ satisfies \eqref{eq:A.6} for $i\in I_2$ due to the assumption that $g^i$ satisfies \eqref{eq:2.5}. Since $(Y^i, Z^i)$ is the unique solution of $i$th equation in BSDE \eqref{eq:1.1}, whose generator can be regarded as $f^i$, and inequality \eqref{eq:2.1} holds, it follows from \eqref{eq:A.11} in \cref{lem:A.1} that for each $i\in I_2$ and $\eps\in (0,\bar\gamma/9]$,\vspace{0.1cm}
$$
\begin{array}{l}
\Dis \E\left[\exp\left(\frac{\bar \gamma}{2}\eps\int_0^T |Z_t^{i}|^2{\rm d}t \right)\right]\leq \Dis \E\left[\exp\left(3\eps(2+\beta T)\sup_{t\in [0,T]}|Y_t^{i}|+3\eps\int_0^T \big(\alpha_t+(\beta+2k)|Y_t| \big){\rm d}t\right)\right]\vspace{0.2cm}\\
\hspace{4.15cm} \leq \Dis \E\left[\exp\left(6\eps(1+\beta T+kT)\sup_{t\in [0,T]}|Y_t|+3\eps\int_0^T \alpha_t {\rm d}t\right)\right].
\end{array}
$$
Therefore, we have \eqref{eq:2.8}.\vspace{0.2cm}

It remains to prove the uniqueness part of \cref{thm:2.2}. For this, let $(\widetilde Y,\widetilde Z)$ be another solution of system of BSDEs \eqref{eq:1.1} in the space $\mathcal{E}(\R^n)\times\mcal(\R^{n\times d})$, and define for $\theta\in (0,1)$,
$$
\delta_{\theta}U:=\frac{Y-\theta \widetilde Y}{1-\theta},\ \ \ \ \ \ \delta_{\theta}V:=\frac{Z-\theta \widetilde Z}{1-\theta}
$$
and
$$
\delta_{\theta}\widetilde U:=\frac{\widetilde Y-\theta Y}{1-\theta},\ \ \ \ \ \ \delta_{\theta}\widetilde V:=\frac{\widetilde Z-\theta Z}{1-\theta}.\vspace{0.2cm}
$$
First of all, similar to the proof of \eqref{eq:2.8}, we can conclude that for each $i\in I_2$ and $\eps\in (0,\bar\gamma/9]$, \vspace{0.1cm}
\begin{equation}\label{eq:3.27}
\Dis \E\left[\exp\left(\frac{\bar \gamma}{2}\eps\int_0^T |\widetilde Z_t^{i}|^2{\rm d}t \right)\right]\leq \Dis \E\left[\exp\left(6\eps(1+\beta T+kT)\sup_{t\in [0,T]}|\widetilde Y_t|+3\eps\int_0^T \alpha_t {\rm d}t\right)\right].
\vspace{0.1cm}
\end{equation}
Moreover, by virtue of \eqref{eq:A.27}, \eqref{eq:3.27}, \eqref{eq:2.8}, \eqref{eq:2.1} and H\"{o}lder's inequality we can obtain the existence of a positive constant $\bar K(q)$ depending only on $(q,n,\beta,\gamma,\bar\gamma,\delta,k,\bar k,T)$ such that for each $q>1$,\vspace{0.1cm}
\begin{equation}\label{eq:3.28}
\begin{array}{l}
\Dis \left(\frac{q}{q-1}\right)^{2q}\E\left[\exp\left\{qc_0 e^{c_0 T} \left(|\xi|+\sup_{s\in\T}\left(|Y_s|+|\widetilde Y_s| \right) \right)\right.\right.\vspace{0.1cm}\\
\hspace{1.4cm} \Dis +\left. \left. qc_0 e^{c_0 T}\int_0^T \left(\alpha_s+k(|Y_s|+|\widetilde Y_s|)+\bar k\sum_{i\in I_2}(|Z^i_s|^{1+\delta}+|\widetilde Z^i_t|^{1+\delta})\right){\rm d}s\right\}\right]\leq \bar K(q).
\end{array}
\end{equation}
Next, in light of \eqref{eq:3.28}, using a similar argument as that from \eqref{eq:3.13} to \eqref{eq:3.20} yields that for each $q>1$,
$$
\begin{array}{l}
\Dis \E\left[\exp\left\{q\sup_{t\in [0,T]}\left[\gamma e^{c_0 t}\left(|\delta_{\theta} U_t|+|\delta_{\theta} \widetilde U_t| \right)\right]\right\}\right]\vspace{0.2cm}\\
\ \ \leq \Dis \left(\frac{q}{q-1}\right)^{q} \E\left[\exp\left\{\frac{qc_0 e^{c_0 T}}{2} \left(|\xi|+\sup_{t\in\T}(|Y_t|+\widetilde Y_t|)\right)+\frac{q}{2}\sup_{t\in [0,T]}\left[\gamma e^{c_0 t} \left(|\delta_{\theta} U_t|+|\delta_{\theta} \widetilde U_t| \right)\right]\right.\right.\vspace{0.1cm}\\
\Dis \hspace{2.8cm} \left.\left.+\frac{qc_0 e^{c_0 T}}{2}\int_0^T \left(\alpha_t+k(|Y_t|+|\widetilde Y_t|)+\bar k \sum_{i\in I_2}(|Z^i_t|^{1+\delta}+ |\widetilde Z^i_t|^{1+\delta})\right){\rm d}t\right\}\right]\vspace{0.1cm}\\
\ \ \leq \Dis \left(\bar K(q)\right)^{1\over 2}\left(\E\left[\exp\left\{q\sup_{t\in [0,T]}\left[\gamma e^{c_0 t}\left(|\delta_{\theta} U_t|+|\delta_{\theta} \widetilde U_t| \right)\right]\right\}\right]\right)^{1\over 2},
\end{array}
$$
and then
\begin{equation}\label{eq:3.29}
\E\left[\exp\left\{q\sup_{t\in [0,T]}\left[\gamma e^{c_0 t}\left(|\delta_{\theta} U_t|+|\delta_{\theta} \widetilde U_t| \right)\right]\right\}\right]\leq \bar K(q).\vspace{0.2cm}
\end{equation}
Finally, by picking $q=2$ in \eqref{eq:3.29} we can derive that
$$
\E\left[2\gamma \sup_{t\in [0,T]}|\delta_{\theta} U_t|\right] \leq \E\left[\exp\left(2\gamma \sup_{t\in [0,T]}|\delta_{\theta} U_t|\right)\right]\leq \bar K(2),\vspace{0.2cm}
$$
and then
$$
\begin{array}{l}
\Dis \E\left[\sup_{t\in [0,T]}|Y_t-\widetilde Y_t|\right]\leq \Dis \E\left[\sup_{t\in [0,T]}|Y_t-\theta \widetilde Y_t|\right]+(1-\theta)\E\left[\sup_{t\in [0,T]}|\widetilde Y_t|\right]\vspace{0.2cm}\\
\hspace{3cm}\leq \Dis (1-\theta)\left(\frac{\bar K(2)}{2\gamma}+\E\left[\exp\left(\sup_{t\in [0,T]}|\widetilde Y_t|\right)\right]\right),\vspace{0.2cm}
\end{array}
$$
in which sending $\theta\To 1$ leads to that $Y=\widetilde Y$ and then $Z=\widetilde Z$. The proof is then complete.\vspace{0.2cm}

\subsection{The proof of \cref{thm:2.5}\vspace{0.2cm}}
\label{sec:4-Stability}

Throughout this section, let all assumptions in \cref{thm:2.5} holds. For each $m\geq 1$ and $i=1,\cdots,n$, define by $g^{(m);i}$ the $i$th component of $g^{(m)}$, and by \cref{rmk:2.1}, with lost of generality, let both $g^i$ and $g^{(m);i}$ satisfy \eqref{eq:2.2}-\eqref{eq:2.3} for all $i=1,\cdots,n$,
\eqref{eq:2.4} with $\bar k=0$ for $i\in I_1$, and \eqref{eq:2.4} and \eqref{eq:2.5} for $i\in I_2$, and $I_1+I_2=\{1,\cdots,n\}$. Furthermore, we define
\begin{equation}\label{eq:4.1}
\bar c:= 64n(\beta+\gamma)e^{\beta T}.\vspace{0.2cm}
\end{equation}

First of all, it follows from assumption \eqref{eq:2.2} of the generator $g$ that for each $i=1,\cdots,n$,
$$
\as,\ \ {\rm sgn}(Y^i_t) g^i(t,Y_t,Z^i_t)\leq \alpha_t+\beta |Y_t|+\frac{\gamma}{2} |Z^i_t|^2.
$$
In light of \eqref{eq:2.12} and the fact that $Y\in \mathcal{E}(\R^n)$, by H\"{o}lder's inequality we have
$$
\RE\ p>1,\ \  \E\left[\exp\left\{p\left(|\xi^i|+\int_0^T \left(\alpha_t+\beta |Y_t|\right){\rm d}t\right)\right\}\right]<+\infty.
$$
Then, applying (i) of \cref{lem:A.3} in Appendix to the $i$th equation of BSDE \eqref{eq:1.1} yields that for each $i=1,\cdots,n$ and each $p>1$,
$$
\exp\left(\gamma |Y_t^i|\right)\leq \E_t\left[\exp\left(\gamma e^{\beta T}|\xi^i|+\gamma e^{\beta T} \int_t^T \left(\alpha_s+\beta |Y_s|\right) {\rm d}s\right)\right],\ \ t\in\T
$$
and
$$
\Dis \E\left[\left(\int_0^T |Z_s^i|^2 {\rm d}s\right)^{p\over 2}\right]\leq \Dis K_p \left\{\E\left[\left(\int_0^T \left(\alpha_s+\beta |Y_s| \right){\rm d}s\right)^p+\exp\left(p\gamma \sup_{s\in\T}|Y_s^i|\right)\right]\right\},\vspace{0.1cm}
$$
where $K_p$ is a positive constant depending only on $(p,\beta,\gamma,T)$. Consequently, by Jensen's inequality we obtain that for each $p>1$,
\begin{equation}\label{eq:4.2}
\exp\left(\gamma |Y_t|\right)\leq \E_t\left[\exp\left(n\gamma     e^{\beta T} |\xi|+n\gamma e^{\beta T} \int_t^T \left(\alpha_s+\beta |Y_s|\right) {\rm d}s\right)\right],\ \ t\in\T
\end{equation}
and\vspace{0.1cm}
\begin{equation}\label{eq:4.3}
\Dis \E\left[\left(\int_0^T |Z_s|^2 {\rm d}s\right)^{p\over 2}\right]\leq \Dis K_p n^{p+1} \left\{\E\left[\left(\int_0^T \left(\alpha_s+\beta |Y_s| \right){\rm d}s\right)^p+\exp\left(p\gamma \sup_{s\in\T}|Y_s|\right)\right]\right\}.\vspace{0.2cm}
\end{equation}
Taking $e^{\bar c t}$-order power in both sides of \eqref{eq:4.2} and using Jensen's inequality, by an identical argument as \eqref{eq:3.9} we can deduce that for each $p>1$,
$$
\exp\left(\gamma e^{\bar c t} |Y_t|\right)\leq \E_t\left[
\frac{\bar c e^{\bar c T}}{64}\left(|\xi|+\int_0^T \alpha_s {\rm }ds\right)+\frac{1}{64}\sup_{s\in \T}\left(\gamma e^{\bar c s} |Y_s|\right)\right].
$$
In light of the last inequality and \eqref{eq:4.1}, using Doob's maximal inequality for martingales and H\"{o}lder's inequality we derive that for each $p>1$,
$$
\begin{array}{ll}
\Dis \E\left[\exp\left(p \sup_{t\in [0,T]}(\gamma e^{\bar c t}|Y_t|)\right)\right]\vspace{0.2cm}\\
\ \ \leq \Dis \left(\frac{p}{p-1}\right)^{p}\E\left[\exp\left\{\frac{p\bar c e^{\bar c T}}{64} \left(|\xi|+\int_0^T \alpha_t{\rm d}t\right)+\frac{p}{64} \sup_{t\in [0,T]}(\gamma e^{\bar c t}|Y_t|)\right\}\right]\vspace{0.2cm}\\
\ \ \leq \Dis \left(\frac{p}{p-1}\right)^{p}\left(\E\left[\exp\left\{\frac{p\bar c e^{\bar c T}}{32} \left(|\xi|+\int_0^T \alpha_t{\rm d}t\right)\right\}\right]\right)^{1\over 2}\left(\E\left[\exp\left(\frac{p}{32} \sup_{t\in [0,T]}(\gamma e^{\bar c t}|Y_t|)\right)\right]\right)^{1\over 2},
\end{array}
$$
and then\vspace{0.1cm}
\begin{equation}\label{eq:4.4}
\E\left[\exp\left(p\sup_{t\in \T}(\gamma e^{\bar c t}|Y_t|)\right)\right]\leq \left(\frac{p}{p-1}\right)^{2p} \E\left[\exp\left\{\frac{p\bar c e^{\bar c T}}{32}\left(|\xi|+\int_0^T \alpha_t{\rm d}t\right)\right\}\right].\vspace{0.2cm}
\end{equation}

Furthermore, a similar argument as the above paragraph on the BSDE associated with $(\xi^{(m)},g^{(m)})$ gives that for each $m\geq 1$ and $p>1$,
\begin{equation}\label{eq:4.5}
\E\left[\exp\left(p\sup_{t\in \T}(\gamma e^{\bar c t}|Y^{(m)}_t|)\right)\right]\leq \left(\frac{p}{p-1}\right)^{2p}\E\left[\exp\left\{\frac{p\bar c e^{\bar c T}}{32}\left(|\xi|+\int_0^T \alpha^{(m)}_t{\rm d}t\right)\right\}\right]
\end{equation}
and
\begin{equation}\label{eq:4.6}
\hspace*{-0.2cm}\Dis \E\left[\left(\int_0^T |Z^{(m)}_s|^2 {\rm d}s\right)^{p\over 2}\right]\leq \Dis K_p n^{p+1} \left\{\E\left[\left(\int_0^T \left(\alpha^{(m)}_s+\beta |Y^{(m)}_s| \right){\rm d}s\right)^p+\exp\left(p\gamma \sup_{s\in\T}|Y^{(m)}_s|\right)\right]\right\}.
\end{equation}
And, since $g^{(m);i}$ satisfies \eqref{eq:2.4} and \eqref{eq:2.5} for $i\in I_2$, by (ii) of \cref{thm:2.2} we can obtain that for each $m\geq 1$, $i\in I_2$ and $\eps\in (0,\bar\gamma/9]$, \vspace{0.1cm}
\begin{equation}\label{eq:4.7}
\E\left[\exp\left(\frac{\bar \gamma}{2}\eps\int_0^T |Z_t^{(m);i}|^2{\rm d}t \right)\right] \leq \E\left[\exp\left(6\eps(1+\beta T+kT)\sup_{t\in [0,T]}|Y^{(m)}_t|+3\eps\int_0^T \alpha^{(m)}_t {\rm d}t\right)\right].\vspace{0.1cm}
\end{equation}
Then, in light of \eqref{eq:2.12}, by combining \eqref{eq:4.5}-\eqref{eq:4.7} and H\"{o}lder's inequality we can verify that the sequence of processes
$(Y^{(m)}_t,Z^{(m)}_t)_{t\in\T}$ satisfies that for each $p>1$,
\begin{equation}\label{eq:4.8}
\sup_{m\geq 1} \E\left[\exp\left(p \sup_{t\in\T}|Y^{(m)}_t|\right)
+\left(\int_0^T|Z^{(m)}_s|^2{\rm d}s\right)^{p\over 2}+\exp\left(\frac{\bar \gamma^2}{18}\int_0^T (\sum_{i\in I_2}|Z_s^{(m);i}|^2){\rm d}s \right)
\right]<+\infty.\vspace{0.2cm}
\end{equation}
Based on \eqref{eq:4.3}, \eqref{eq:4.4}, \eqref{eq:4.8}, \eqref{eq:2.8} and \eqref{eq:A.27}, from H\"{o}lder's inequality we can derive the existence of a positive constant $C(p)$ depending only on $(p,\beta,\gamma,\bar\gamma, \delta, k,\bar k, T)$ such that for each $p>0$,
\begin{equation}\label{eq:4.9}
4\sup_{m\geq 1}\E\left[\exp\left\{p\int_0^T \left(\alpha_t^{(m)}+k(|Y_t^{(m)}|+|Y_t|)+\bar k\sum_{i\in I_2}(|Z_t^{(m);i}|^{1+\delta}+|Z^i_t|^{1+\delta})\right){\rm d}t \right\}\right]\leq C(p). \vspace{0.2cm}
\end{equation}

In the sequel, for each fixed $\theta\in (0,1)$, define
$$
\delta_\theta^m  U:=\frac{Y^{(m)}-\theta Y}{1-\theta}\ \  {\rm and} \ \ \delta_\theta^m  V:=\frac{Z^{(m)}-\theta Z}{1-\theta}.
$$
Then, the pair $(\delta_\theta^m  U,\delta_\theta^m  V)$ verifies the following system of BSDEs:
$$
\delta_\theta^m  U_t^i=\delta_\theta^m  U_T^i +\int_t^T \delta_\theta^m g^i (s) {\rm d}s-\int_t^T \delta_\theta^m  V_s^i {\rm d}B_s, \ \ \ \ t\in\T;\ \ \ i=1,\cdots,n,
$$
where $\ass$,
\begin{equation}\label{eq:4.10}
\Dis \delta_\theta^m g^i(s) := \Dis \frac{g^{(m);i}(s, Y^{(m)}_s, Z^{(m);i}_s)-\theta g^{(m);i}(s, Y_s, Z^i_s)}{1-\theta}
\end{equation}
with
\begin{equation}\label{eq:4.11}
\delta_\theta^m D^i_s:=\frac{\theta\left( g^{(m);i}(s,Y_s, Z^i_s)-g^i(s,Y_s, Z^i_s) \right)}{1-\theta}.
\vspace{0.1cm}
\end{equation}
Since $g^{(m);i}$ satisfies \eqref{eq:2.4} with $\bar k=0$ for $i\in I_1$, and \eqref{eq:2.4} and \eqref{eq:2.5} for $i\in I_2$, it follows from \eqref{eq:4.10} that
\begin{equation}\label{eq:4.12}
\begin{array}{l}
\ass,\ \ \Dis {\bf 1}_{\{\delta^m_{\theta} U^i_s>0\}}\delta^m_{\theta} g^i(s)\leq \Dis \alpha^{(m)}_s+\beta |\delta^m_{\theta}U_s|+\frac{\gamma}{2}|\delta^m_{\theta} V^i_s|^2+k(|Y^{(m)}_s|+|Y_s|)\vspace{0.1cm}\\
\hspace{5.8cm}\Dis +\bar k{\bf 1}_{\{i\in I_2\}}(|Z_s^{(m);i}|^{1+\delta}+|Z^i_s|^{1+\delta})+|\delta_\theta^m D_s^i|.
\end{array}
\end{equation}
In light of \eqref{eq:4.12}, \eqref{eq:4.1}, \eqref{eq:4.3}, \eqref{eq:4.4}, \eqref{eq:4.8}, \eqref{eq:4.9}, \eqref{eq:2.12} and \eqref{eq:2.13}, by (ii) of \cref{lem:A.3} with
$$
p:=e^{\bar c}t,\ \ \hat\alpha_s:=\alpha^{(m)}_s+ k(|Y^{(m)}_s|+|Y_s|)+\bar k \sum_{i\in I_2}(|Z_s^{(m);i}|^{1+\delta}+|Z^i_s|^{1+\delta})+\beta |\delta^m_{\theta}U_s|\ \ {\rm and}\ \ \check\alpha_s:=|\delta_\theta^m D_s^i|,
$$
from Doob's maximal inequality and H\"{o}lder's inequality we can deduce that, similar to \eqref{eq:3.9},
\begin{equation}\label{eq:4.13}
\begin{array}{l}
\Dis \E\left[\exp\left(2\sup_{t\in\T}\left[\gamma e^{\bar c t}\left(\delta_\theta^m  U_t^i\right)^+\right]\right)\right] \vspace{0.2cm}\\
\ \ \leq \Dis C\left(4\gamma e^{(\beta+\bar c) T}\right)\cdot \E\left[\exp\left(4\gamma e^{(\beta+\bar c) T}\left(\delta_\theta^m  U_T^i\right)^+ +\frac{1}{16n}\sup_{t\in\T} \left[\gamma e^{\bar c t} |\delta_\theta^m  U_t|\right]\right)\right]\vspace{0.2cm}\\
\ \ \ \ \ \ \Dis +2\gamma e^{(\beta+\bar c)T} \sup_{m\geq 1}\E\left[\exp\left(8\gamma e^{(\beta+\bar c)T}\sup_{t\in\T} \left(\delta_\theta^m  U_t^i\right)^+ +\frac{1}{8n}\sup_{t\in\T} \left[\gamma e^{\bar c t} |\delta_\theta^m  U_t|\right]\right)\right]\vspace{0.2cm}\\
\ \ \Dis \ \ \ \ \ \cdot C\left(8\gamma e^{(\beta+\bar c)T}\right)\cdot \left(\E\left[\left(\int_0^T |\delta_\theta^m D_t^i| {\rm d}t\right)^2\right]\right)^{1\over 2},\ \ i=1,\cdots,n.
\end{array}
\end{equation}
In addition, in light of \eqref{eq:4.11}, \eqref{eq:2.13}, \eqref{eq:2.12} and \eqref{eq:4.8} as well as
\eqref{eq:2.3} for $g^{(m)}$ and $g$, we can conclude that
\begin{equation}\label{eq:4.14}
\lim\limits_{m\To\infty} \E\left[\left(\int_0^T \left| \delta_\theta^m D_t^i \right| {\rm d}s\right)^2\right]=0,\ \ i=1,\cdots, n.
\end{equation}
Thus, in light of \eqref{eq:4.14} and the assumption that $\ps$, $\lim_{m\To\infty}|\xi^{(m)}-\xi|=0$, sending $m$ to infinity in \eqref{eq:4.13} and using H\"{o}lder's inequality yields that for each $i=1,\cdots,n$,
\begin{equation}\label{eq:4.15}
\limsup_{m\To \infty} \E\left[\exp\left(2\sup_{t\in\T}\left[\gamma e^{\bar c t}\left(\delta_\theta^m  U_t^i\right)^+\right]\right)\right] \leq  \Dis \widetilde C \limsup_{m\To \infty} \E\left[\exp\left(\frac{1}{8n}\sup_{t\in\T} \left[\gamma e^{\bar c t} |\delta_\theta^m  U_t|\right]\right)\right],
\end{equation}
where
$$
\widetilde C: =C\left(4\gamma e^{(\beta+\bar c)T}\right)\E\left[\exp\left(8\gamma e^{(\beta+\bar c)T}|\xi|\right)\right].\vspace{0.2cm}
$$

On the other hand, for each fixed $\theta\in (0,1)$, define
$$
\delta_\theta^m \widetilde U:=\frac{Y-\theta Y^{(m)}}{1-\theta}\ \  {\rm and} \ \ \delta_\theta^m  \widetilde V:=\frac{Z-\theta Z^{(m)}}{1-\theta}.
$$
Then the pair $(\delta_\theta^m  \widetilde U,\delta_\theta^m \widetilde V)$ solves the following system of BSDEs: \vspace{0.1cm}
$$
\delta_\theta^m \widetilde U_t^i=\delta_\theta^m \widetilde U_T^i +\int_t^T \delta_\theta^m \widetilde g^i (s) {\rm d}s-\int_t^T \delta_\theta^m  \widetilde V_s^i {\rm d}B_s, \ \ \ \ t\in\T;\ \ \ i=1,\cdots,n, \vspace{0.1cm}
$$
where $\ass$,
\begin{equation}\label{eq:4.16}
\Dis \delta_\theta^m \widetilde g^i(s):= \Dis \delta_\theta^m \widetilde D^i_s+\frac{g^{(m);i}(s, Y_s, Z^i_s)-\theta g^{(m);i}(s, Y^{(m)}_s, Z^{(m);i}_s)  }{1-\theta}
\end{equation}
with
$$
\delta_\theta^m \widetilde D^i_s:=\frac{g^i(s,Y_s, Z^i_s)-g^{(m);i}(s,Y_s, Z^i_s)}{1-\theta}=-\frac{\delta_\theta^m D^i_s}{\theta}.
\vspace{0.1cm}
$$
Since $g^{(m);i}$ satisfies \eqref{eq:2.4} with $\bar k=0$ for $i\in I_1$, and \eqref{eq:2.4} and \eqref{eq:2.5} for $i\in I_2$, it follows from \eqref{eq:4.10} that $\ass$,
$$
\begin{array}{l}
\Dis {\bf 1}_{\{\delta^m_{\theta} \widetilde U^i_s>0\}}\delta^m_{\theta} \widetilde g^i(s)\leq \Dis \alpha^{(m)}_s+\beta |\delta^m_{\theta} \widetilde U_s|+\frac{\gamma}{2}|\delta^m_{\theta} \widetilde V^i_s|^2+k(|Y^{(m)}_s|+|Y_s|)\vspace{0.1cm}\\
\hspace{3.3cm}\Dis +\bar k{\bf 1}_{\{i\in I_2\}}(|Z_s^{(m);i}|^{1+\delta}+|Z^i_s|^{1+\delta})+|\delta_\theta^m \widetilde D_s^i|.
\end{array}
$$
Thus, we can use the same argument as those from \eqref{eq:4.12} to \eqref{eq:4.15} to deduce that that for each $i=1,\cdots,n$,
\begin{equation}\label{eq:4.17}
\limsup_{m\To \infty} \E\left[\exp\left(2\sup_{t\in\T}\left[\gamma e^{\bar c t}\left(\delta_\theta^m \widetilde U_t^i\right)^+\right]\right)\right] \leq  \Dis \widetilde C \limsup_{m\To \infty} \E\left[\exp\left(\frac{1}{8n}\sup_{t\in\T} \left[\gamma e^{\bar c t} |\delta_\theta^m \widetilde U_t|\right]\right)\right].\vspace{0.3cm}
\end{equation}

Next, observe that for $i=1,\cdots,n$ and $t\in\T$,\vspace{0.1cm}
$$
\begin{array}{l}
\Dis \left(\delta_{\theta}^m U_t^i\right)^- = \Dis \frac{\left( Y_t^{(m);i}-\theta Y_t^i\right)^-} {1-\theta}= \Dis \frac{\left(\theta Y_t^i- Y_t^{(m);i}\right)^+} {1-\theta}\vspace{0.2cm}\\
\hspace{1.5cm}\leq \Dis\frac{\theta\left( Y_t^i-\theta Y_t^{(m);i}\right)^+ +(1-\theta^2)|Y^{(m);i}_t|} {1-\theta}\vspace{0.2cm}\\
\hspace{1.5cm} \leq \Dis \left(\delta_{\theta}^m \widetilde U_t^i\right)^+ +2|Y^{(m)}_t|
\end{array}
$$
and, similarly,
$$
\left(\delta_{\theta}^m \widetilde U_t^i \right)^- \leq \left(\delta_{\theta}^m U_t^i\right)^+ +2|Y_t|.
$$
From \eqref{eq:4.15}, \eqref{eq:4.17} and H\"{o}lder's inequality, it follows that for each $i=1,\cdots, n$,
$$
\begin{array}{l}
\Dis \limsup_{m\To \infty} \E\left[\exp\left(\frac{1}{2}\sup_{t\in\T}\left[\gamma e^{\bar c t} |\delta_\theta^m  U_t^i|\right]\right)\right]\vspace{0.2cm}\\
\ \ =\Dis \limsup_{m\To \infty} \E\left[\exp\left(\frac{1}{2} \sup_{t\in\T}\left[\gamma e^{\bar c t} \left(\delta_\theta^m  U_t^i\right)^+\right]\right)\cdot \exp\left(\frac{1}{2} \sup_{t\in\T}\left[\gamma e^{\bar c t} \left(\delta_\theta^m  U_t^i\right)^-\right]\right)\right]\vspace{0.2cm}\\
\ \ \leq \Dis \widetilde C \E\left[\exp\left(2\sup_{t\in\T} \left[\gamma e^{\bar c t} |Y^{(m)}_t|\right]\right)\right]\cdot \limsup_{m\To \infty} \E\left[\exp\left(\frac{1}{8n}\sup_{t\in\T}\left[\gamma e^{\bar c t}  \left(|\delta_\theta^m  U_t|+ |\delta_\theta^m \widetilde U_t|\right)\right] \right)\right]
\end{array}
$$
and
$$
\begin{array}{l}
\Dis \limsup_{m\To \infty} \E\left[\exp\left(\frac{1}{2}\sup_{t\in\T}\left[\gamma e^{\bar c t} |\delta_\theta^m \widetilde U_t^i|\right]\right)\right]\vspace{0.2cm}\\
\ \ =\Dis \limsup_{m\To \infty} \E\left[\exp\left(\frac{1}{2} \sup_{t\in\T}\left[\gamma e^{\bar c t} \left(\delta_\theta^m  \widetilde U_t^i\right)^+\right]\right)\cdot \exp\left(\frac{1}{2} \sup_{t\in\T}\left[\gamma e^{\bar c t} \left(\delta_\theta^m  \widetilde U_t^i\right)^-\right]\right)\right]\vspace{0.2cm}\\
\ \ \leq \Dis \widetilde C \E\left[\exp\left(2\sup_{t\in\T} \left[\gamma e^{\bar c t} |Y_t|\right]\right)\right]\cdot \limsup_{m\To \infty} \E\left[\exp\left(\frac{1}{8n}\sup_{t\in\T}\left[\gamma e^{\bar c t}  \left(|\delta_\theta^m  U_t|+ |\delta_\theta^m \widetilde U_t|\right)\right] \right)\right].\vspace{0.1cm}
\end{array}
$$
Consequently, by H\"{o}lder's inequality again we deduce that for each $i=1,\cdots, n$,
$$
\begin{array}{l}
\Dis \limsup_{m\To \infty} \E\left[\exp\left(\frac{1}{4}\sup_{t\in\T}\left[\gamma e^{\bar c t}  \left(|\delta_\theta^m  U^i_t|+ |\delta_\theta^m \widetilde U^i_t|\right)\right] \right)\right]\vspace{0.2cm}\\
\ \ \leq \Dis \widetilde K \left(\limsup_{m\To \infty} \E\left[\exp\left(\frac{1}{4n}\sup_{t\in\T}\left[\gamma e^{\bar c t}  \left(|\delta_\theta^m  U_t|+ |\delta_\theta^m \widetilde U_t|\right)\right] \right)\right]\right)^{1\over 2},\vspace{0.1cm}
\end{array}
$$
and then, in light of \eqref{eq:4.4} and \eqref{eq:4.8},
$$
\begin{array}{l}
\Dis \limsup_{m\To \infty} \E\left[\exp\left(\frac{1}{4n}\sup_{t\in\T}\left[\gamma e^{\bar c t}  \left(|\delta_\theta^m  U_t|+ |\delta_\theta^m \widetilde U_t|\right)\right] \right)\right]\vspace{0.2cm}\\
\ \ \leq \Dis  \widetilde K \left(\limsup_{m\To \infty} \E\left[\exp\left(\frac{1}{4n}\sup_{t\in\T}\left[\gamma e^{\bar c t}  \left(|\delta_\theta^m  U_t|+ |\delta_\theta^m \widetilde U_t|\right)\right] \right)\right]\right)^{1\over 2}<+\infty
\end{array}
$$
with
$$
\widetilde K :=\widetilde C \E\left[\exp\left(2\sup_{t\in\T} \left[\gamma e^{\bar c t} |Y_t|\right]\right)\right]\cdot \sup_{m\geq 1}\E\left[\exp\left(2\sup_{t\in\T} \left[\gamma e^{\bar c t} |Y^{(m)}_t|\right]\right)\right].\vspace{0.2cm}
$$
Thus,
\begin{equation}\label{eq:4.18}
\limsup_{m\To \infty} \E\left[\exp\left(\frac{1}{4n}\sup_{t\in\T}\left[\gamma e^{\bar c t}  \left(|\delta_\theta^m  U_t|+ |\delta_\theta^m \widetilde U_t|\right)\right] \right)\right]\leq \widetilde K^2. \vspace{0.4cm}
\end{equation}

Finally, it follows from the last inequality that for each $\theta\in (0,1)$,
$$
\limsup_{m\To \infty}\E\left[\sup_{t\in\T}|Y_t^{(m)}-\theta Y_t|\right] \leq (1-\theta){4n\over \gamma} \widetilde K^2.
$$
Therefore,
\begin{equation}\label{eq:4.19}
\limsup_{m\To \infty}\E\left[\sup_{t\in\T}|Y_t^{(m)}-Y_t|\right] \leq (1-\theta)\left( {4n\over \gamma} \widetilde K^2 +\E\left[\sup_{t\in\T}|Y_t|\right]\right). \vspace{0.2cm}
\end{equation}
In light of \eqref{eq:4.4}, sending $\theta\To 1$ in \eqref{eq:4.19} we  can derive that for each $p\geq 1$,
\begin{equation}\label{eq:4.20}
\lim\limits_{m\To \infty} \E\left[\sup_{t\in\T}|Y_t^{(m)}-Y_t|^p\right]=0.
\end{equation}
Furthermore, it follows from It\^{o}'s formula that for each $m\geq 1$,
\begin{equation}\label{eq:4.21}
\begin{array}{l}
\Dis \E\left[\int_0^T |Z^{(m)}_s-Z_s|^2 {\rm d}s\right]\vspace{0.3cm}\\
\ \ \leq \Dis \E\left[|\xi^{(m)}-\xi|^2+2\sup_{t\in\T}|Y^{(m)}_t-Y_t| \int_0^T \left|g^{(m)}(s,Y^{(m)}_s, Z^{(m)}_s)-g(s,Y_s, Z_s)\right| {\rm d}s \right].
\end{array}
\end{equation}
And, by virtue of \eqref{eq:2.3} for $g^{(m)}$ and $g$ as well as \eqref{eq:2.12}, \eqref{eq:4.4} and \eqref{eq:4.8} we obtain that
\begin{equation}\label{eq:4.22}
\sup_{m\geq 1}\E\left[\left(\int_0^T \left|g^{(m)}(s,Y^{(m)}_s, Z^{(m)}_s)-g(s,Y_s, Z_s)\right| {\rm d}s\right)^2\right]<+\infty.
\end{equation}
Then, by using H\"{o}lder's inequality, the desired conclusion \eqref{eq:2.14} follows immediately from \eqref{eq:4.20}-\eqref{eq:4.22} along with \eqref{eq:4.8}. The proof of \cref{thm:2.5} is then complete.\vspace{0.2cm}

\subsection{The proof of \cref{thm:2.7}\vspace{0.2cm}}
\label{sec:5-Mixedcase}

The proof of \cref{thm:2.7} is a modification one of \cref{thm:2.2}. For the sake of simplicity, by (i) of \cref{rmk:2.6} we assume without loss of generality that $f=g^i$ satisfies \eqref{eq:2.13-1}-\eqref{eq:2.15} in assumption \ref{A:H3} for all $i=1,\ldots,l$, and \eqref{eq:2.2}-\eqref{eq:2.4}  with $\bar k=0$ in assumption \ref{A:H2} for all $i=l+1,\cdots,n$.\vspace{0.2cm}

For a pair of processes $(U,V)\in \mathcal{E}(\R^n)\times \mcal (\R^{n\times d})$, consider the following partly-decoupled system of quadratic BSDEs: \vspace{0.1cm}
\begin{equation}\label{eq:5.1}
\Dis Y^i_t=\xi^i+\int_t^T g^i(s,U_s(Y^{[1:l]}_s),Z^{[1:l]}_s){\rm d}s-\int_t^T Z^i_s {\rm d}B_s, \ \ t\in\T,\ \ 1\leq i\leq l
\end{equation}
and
\begin{equation}\label{eq:5.2}
\Dis Y_t^i=\xi^i+\int_t^T g^i(s,U_s(Y_s^i;i),Z_s^i){\rm d}s-\int_t^T Z_s^i {\rm d}B_s, \ \ t\in\T,\ \ l+1\leq i\leq n.\vspace{0.1cm}
\end{equation}
In light of assumptions \ref{A:H1} and \ref{A:H3}, it follows from \eqref{eq:2.20}, \eqref{eq:2.17} and \eqref{eq:2.18} that $\as$, for each $(y,z)\in\R^n\times\R^{n\times d}$ and $(y_1,y_2,z_1,z_2)\in \R^n\times \R^n\times \R^{n\times d}\times \R^{n\times d}$, we have
\begin{equation}\label{eq:5.3}
\left|g^{[1:l]}(t,U_t(0^{[1:l]}),0)\right|\leq l\alpha_t+lk\exp{(2k|U_t|)}+lk,
\end{equation}
\begin{equation}\label{eq:5.4}
\Dis\left\langle y^{[1:l]}_1-y^{[1:l]}_2, g^{[1:l]}(t,U_t(y^{[1:l]}_1),z^{[1:l]})-
g^{[1:l]}(t,U_t(y^{[1:l]}_2),z^{[1:l]})\right\rangle
\Dis \leq \sqrt{l}\beta |y^{[1:l]}_1-y^{[1:l]}_2|^2
\end{equation}
and
\begin{equation}\label{eq:5.5}
\left|g^{[1:l]}(t,U_t(y^{[1:l]}),z^{[1:l]}_1)-
g^{[1:l]}(t,U_t(y^{[1:l]}),z^{[1:l]}_2) \right|\leq l\gamma |z^{[1:l]}_1-z^{[1:l]}_2|.\vspace{0.1cm}
\end{equation}
Then, in light of the last three inequalities \eqref{eq:5.3}-\eqref{eq:5.5} along with \eqref{eq:2.1} and the fact of $U\in \mathcal{E} (\R^n)$, by \citet[Theorem 4.1]{Briand-etal2003SPA} we can conclude that for each $p>1$, BSDE \eqref{eq:5.1} admits a unique solution $(Y^i,Z^i)$ with $1\leq i\leq l$ in the space of $\s^p(\R^l)\times \hcal^p(\R^{l\times d})$. Moreover, it follow from \eqref{eq:2.19} that $\as$, for each $(y,z)\in \R^n\times\R^{n\times d}$,
\begin{equation}\label{eq:5.6}
\Dis \left\langle \frac{y^{[1:l]}}{|y^{[1:l]}|}{\bf 1}_{|y^{[1:l]}|>0}, g^{[1:i]}(t,U_t(y^{[1:l]}),z^{[1:l]})\right\rangle\leq \sqrt{l}\alpha_t+\sqrt{l}\beta |U_t|+\sqrt{l}\beta |y^{[1:l]}|+\sqrt{l}\gamma |z^{[1:l]}|.
\end{equation}
Then, in light of \eqref{eq:5.1} and \eqref{eq:5.6}, by \cite[Proposition 3.2]{Briand-etal2003SPA} we can deduce that for each $q>1$, there exists a positive constant $\bar C_q$ depending only on $(q,l,\beta,\gamma,T)$ such that
\begin{equation}\label{eq:5.7}
\begin{array}{l}
\Dis \E\left[\sup_{t\leq s\leq T}\left(\sum_{i=1}^{l}|Y^i_s|^2\right)^{q\over 2}+\left(\int_t^T \sum_{i=1}^{l}|z^i_s|^2{\rm d}s\right)^{q\over 2}\bigg|\F_t\right]\vspace{0.1cm}\\
\ \ \Dis \leq \bar C_q \E\left[\sum_{i=1}^{l}|\xi^i|^q+\left(\int_t^T \left(\alpha_s+|U_s|\right){\rm d}s\right)^{q}\bigg|\F_t\right],\ \ t\in \T.
\end{array}
\end{equation}
In particular, by letting $q=2$ we know that there exists a constant $C\geq 1$ such that
$$
\sum_{i=1}^{l}|Y^i_t|^2\leq C \E\left[\sum_{i=1}^{l}|\xi^i|^2+\left(\int_t^T \left(\alpha_s+|U_s|\right){\rm d}s\right)^2\bigg|\F_t\right],\ \ t\in \T,
$$
and then, with $c_0$ being a positive constant, which will be determined later,
$$
\begin{array}{l}
\Dis \gamma^2 e^{2c_0 t}\sum_{i=1}^{l}|Y^i_t|^2\leq \gamma^2 C \E\left[e^{2c_0 T}\left\{|\xi|^2+\left(\int_0^T \alpha_s{\rm d}s\right)^2\right\}+\left(e^{c_0 t} \int_t^T \left(e^{c_0 s}|U_s|\cdot e^{-c_0 s}\right){\rm d}s\right)^2\bigg|\F_t\right]\vspace{0.1cm}\\
\hspace{2.1cm}\Dis \leq \gamma^2 C \E\left[e^{2c_0 T}\left\{|\xi|^2+\left(\int_0^T \alpha_s{\rm d}s\right)^2\right\}+\frac{1}{c^2_0}\sup_{s\in \T}\left(e^{2c_0 s}|U_s|^2\right)\bigg|\F_t\right],\ \ t\in \T.
\end{array}
$$
Since the function $f(x)=\exp{(\sqrt{x+1}\,)}, x\ge 0,$ is  convex, from the last inequality and Jensen's inequality,  we  deduce that for each $t\in \T$,
$$
\Dis \exp\left(\gamma e^{c_0 t}\sqrt{\sum_{i=1}^{l}|Y^i_t|^2}\right)\leq \E\left[\exp \left\{\gamma C e^{c_0 T}\left(|\xi|+\int_0^T \alpha_s{\rm d}s\right)+\frac{C}{c_0}\sup_{s\in \T}(\gamma e^{c_0 s}|U_s|)+1\right\}\bigg|\F_t\right].
$$
Furthermore, by applying Doob's maximal inequality for martingales to the last inequality we can obtain that for each $q>1$,
\begin{equation}\label{eq:5.8}
\begin{array}{l}
\Dis \E\left[\exp\left\{q \sup_{t\in\T}\left(\gamma e^{c_0 t}\sqrt{\sum_{i=1}^{l}|Y^i_t|^2}\right)\right\}\right]\vspace{0.1cm}\\
\ \ \Dis \leq \left(\frac{eq}{q-1}\right)^{q} \E\left[\exp\left\{q\gamma Ce^{c_0 T}\left(|\xi|+\int_0^T \alpha_s{\rm d}s\right)+\frac{qC}{c_0}\sup_{t\in \T}(\gamma e^{c_0 t}|U_t|)\right\}\right].\vspace{0.2cm}
\end{array}
\end{equation}

On the other hand, the proof of \cref{thm:2.2} shows that BSDE \eqref{eq:5.2} admits a unique adapted solution $(Y^i,Z^i)$ with $l+1\leq i\leq n$ such that
\begin{equation}\label{eq:5.9}
\E\left[\exp\left(q\sup_{t\in \T}\sqrt{\sum_{i=l+1}^{n}|Y^i_t|^2}\right)+\left(\int_0^T \sum_{i=l+1}^{n}|Z_t^i|^2{\rm d}t\right)^{q\over 2}\right]<+\infty, \quad \forall q>1.
\end{equation}
Thus, in light of \eqref{eq:5.7}-\eqref{eq:5.9}, we conclude that the system of multidimensional BSDEs \eqref{eq:5.1} and \eqref{eq:5.2} admits a unique solution $(Y,Z):=(Y^1,\cdots,Y^n,Z^1,\cdots,Z^n)$ in the space $\mathcal{E}(\R^n)\times \mcal (\R^{n\times d})$.\vspace{0.2cm}

Based on the analysis in the preceding two paragraphs, we set $(Y^{(0)},Z^{(0)})=(0,0)$ and define, recursively, the sequence of processes $\{(Y^{(m)},Z^{(m)})\}_{m=1}^\infty$ in the space of processes $\mathcal{E}(\R^n)\times \mcal (\R^{n\times d})$ as the unique adapted solution of multidimensional quadratic BSDEs: for $1\leq i\leq l$,
\begin{equation}\label{eq:5.10}
  Y_t^{(m+1);i}=\xi^i+\int_t^T g^i(s,Y_s^{(m)}(Y_s^{(m+1);[1:l]}),Z_s^{(m+1);[1:l]}){\rm d}s-\int_t^T Z_s^{(m+1);i}{\rm d}B_s, \ \ t\in\T
\end{equation}
and for $l+1\leq i\leq n$,
\begin{equation}\label{eq:5.11}
  Y_t^{(m+1);i}=\xi^i+\int_t^T g^i(s,Y_s^{(m)}(Y_s^{(m+1);i};i),Z_s^{(m+1);i}){\rm d}s-\int_t^T Z_s^{(m+1);i}{\rm d}B_s, \ \ t\in\T.\vspace{0.2cm}
\end{equation}

We now show that $\{(Y^{(m)},Z^{(m)})\}_{m=1}^\infty$ is a Cauchy sequence in the space $\s^q(\R^n)\times \hcal^q(\R^{n\times d})$ for each $q>1$, and then converges to a pair of processes $(Y,Z)$ in $\mathcal{E}(\R^n)\times \mcal(\R^{n\times d})$, which is the desired solution of BSDE \eqref{eq:1.1}.\vspace{0.2cm}

We first show that for each $q>1$,
\begin{equation}\label{eq:5.12}
\sup_{m\geq 0}\E\left[\exp\left(q \sup_{t\in \T}(\gamma e^{c_0 t}|Y_t^{(m)}|)\right)\right]\leq C(q)
\end{equation}
with $c_0:=32n(\beta+\gamma)e^{\beta T}+8(1+\gamma)C$ and
\begin{equation}\label{eq:5.13}
C(q):=\left(\frac{eq}{q-1}\right)^{4q} \E\left[\exp\left\{qc_0 e^{c_0 T} \left(|\xi|+\int_0^T \alpha_t{\rm d}t\right)\right\}\right]<+\infty,
\end{equation}
and that
\begin{equation}\label{eq:5.14}
\sup_{m\geq 0}\E\left[\left(\int_0^T |Z_s^{(m)}|^2 {\rm d}s\right)^{q\over 2}\right]<+\infty. \vspace{0.2cm}
\end{equation}
In fact, in light of \eqref{eq:5.10} and the definition of $c_0$, it follows from \eqref{eq:5.8} that for each $m\geq 0$ and $q>1$,
\begin{equation}\label{eq:5.15}
\begin{array}{l}
\Dis \E\left[\exp\left\{q \sup_{t\in\T}\left(\gamma e^{c_0 t}\sqrt{\sum_{i=1}^{l}|Y^{(m+1);i}_t|^2}\right)\right\}\right]\vspace{0.1cm}\\
\ \ \Dis \leq \left(\frac{eq}{q-1}\right)^{q} \E\left[\exp \left\{\frac{q c_0 e^{c_0 T}}{8}\left(|\xi|+\int_0^T \alpha_s{\rm d}s\right)+\frac{q}{8}\sup_{t\in \T}(\gamma e^{c_0 t}|Y^{(m)}_t|)\right\}\right].
\end{array}
\end{equation}
On the other hand, in light of \eqref{eq:5.11}, by a similar proof of \eqref{eq:3.11} we can verify that for $m\geq 0$ and $q>1$,
\begin{equation}\label{eq:5.16}
\begin{array}{l}
\Dis \E\left[\exp\left\{q \sup_{t\in\T}\left(\gamma e^{c_0 t}\sqrt{\sum_{i=l+1}^{n}|Y^{(m+1);i}_t|^2}\right)\right\}\right]\vspace{0.1cm}\\
\ \ \leq \Dis \left(\frac{q}{q-1}\right)^{q} \E\left[\exp\left\{\frac{q c_0 e^{c_0 T}}{8} \left(|\xi|+\int_0^T \alpha_t{\rm d}t\right)+\frac{q}{8}\sup_{s\in [0,T]}\left(\gamma e^{c_0 s} |Y_s^{(m)}|\right)\right\}\right].\vspace{0.1cm}
\end{array}
\end{equation}
By combining \eqref{eq:5.15} and \eqref{eq:5.16} and using H\"{o}lder's inequality we obtain that for $m\geq 0$ and $q>1$,
\begin{equation}\label{eq:5.17}
\begin{array}{l}
\Dis \E\left[\exp\left\{q \sup_{t\in\T}\left(\gamma e^{c_0 t}|Y^{(m+1)}_t|\right)\right\}\right]\vspace{0.1cm}\\
\ \ \leq \Dis \left(\frac{eq}{q-1}\right)^{2q}  \E\left[\exp\left\{\frac{q c_0 e^{c_0 T}}{4} \left(|\xi|+\int_0^T \alpha_t{\rm d}t\right)+\frac{q}{4}\sup_{s\in [0,T]}\left(\gamma e^{c_0 s} |Y_s^{(m)}|\right)\right\}\right]\vspace{0.1cm}\\
\ \ \leq \Dis (C(q))^{1\over 2} \left\{\E\left[\exp\left\{q \sup_{t\in\T}\left(\gamma e^{c_0 t}|Y^{(m)}_t|\right)\right\}\right]\right\}^{1\over 2},\vspace{0.1cm}
\end{array}
\end{equation}
where $C(q)$ is defined in \eqref{eq:5.13}. Then, the desired assertion \eqref{eq:5.12} follows by induction for \eqref{eq:5.17}. Furthermore, by \eqref{eq:5.7} and \eqref{eq:3.10} we can derive that there exists a positive constant $K_q$ depending only on $(q,l,n,\beta,\gamma,T)$ such that for each $m\geq 0$ and $q>1$,
\begin{equation}\label{eq:5.18}
\Dis \E\left[\left(\int_0^T \sum_{i=1}^l |Z_s^{(m+1):i}|^2 {\rm d}s\right)^{q\over 2}\right]\leq K_q \E\left[\sum_{i=l}^l |\xi^i|^q +\left(\int_0^T \left(\alpha_s+\beta |Y_s^{(m)}| \right){\rm d}s\right)^q\right]
\end{equation}
and
\begin{equation}\label{eq:5.19}
\begin{array}{l}
\Dis \E\left[\left(\int_0^T \sum_{i=l+1}^n |Z_s^{(m+1);i}|^2 {\rm d}s\right)^{q\over 2}\right]\vspace{0.1cm}\\
\ \ \Dis \leq K_q \E\left[\left(\int_0^T \left(\alpha_s+\beta |Y_s^{(m)}| \right){\rm d}s\right)^q+\exp\left(q\gamma \sup_{s\in\T}|Y_s^{(m+1)}|\right)\right].
\end{array}\vspace{0.1cm}
\end{equation}
Then, the desired inequality \eqref{eq:5.14} follows \eqref{eq:5.18}, \eqref{eq:5.19} and \eqref{eq:5.12}.\vspace{0.2cm}

In the sequel, we define, for each fixed $m,p\geq 1$ and $\theta\in (0,1)$,
$$
\delta_{\theta}Y^{(m,p)}:=\frac{Y^{(m+p)}-\theta Y^{(m)}}{1-\theta}\ \ \ {\rm and}\ \ \ \delta_{\theta}Z^{(m,p)}:=\frac{Z^{(m+p)}-\theta Z^{(m)}}{1-\theta}.
$$
Then, in light of \eqref{eq:5.10} and \eqref{eq:5.11}, $(\delta_{\theta}Y^{(m,p)},\delta_{\theta}Z^{(m,p)})\in \mathcal{E}(\R^n)\times \mcal(\R^{n\times d})$ solves the following system of quadratic BSDEs: for $i=1,\cdots, n$,\vspace{0.1cm}
\begin{equation}\label{eq:5.20}
  \delta_{\theta} Y_t^{(m,p);i}=\xi^i+\int_t^T \delta_{\theta} g^{(m,p);i}(s){\rm d}s-\int_t^T \delta_{\theta} Z_s^{(m,p);i}{\rm d}B_s, \ \ t\in\T,\vspace{0.1cm}
\end{equation}
where $\ass$, for $1\leq i\leq l$,
\begin{equation}\label{eq:5.21}
\hspace*{-0.2cm}\Dis\delta_{\theta} g^{(m,p);i}(s)
:=\Dis\frac{g^{i}\left(s, Y_s^{(m+p-1)}(Y_s^{(m+p);[1;l]}), Z_s^{(m+p);[1:l]}\right)-\theta g^{i}\left(s,Y_s^{(m-1)}(Y_s^{(m);[1:l]}),Z_s^{(m);[1:l]}\right)}{1-\theta},
\end{equation}
and for $l+1\leq i\leq n$,
\begin{equation}\label{eq:5.22}
\Dis\delta_{\theta} g^{(m,p);i}(s)
:=\Dis\frac{g^{i}\left(s, Y_s^{(m+p-1)}(Y_s^{(m+p);i};i), Z_s^{(m+p);i}\right)-\theta g^{i}\left(s,Y_s^{(m-1)}(Y_s^{(m);i};i),Z_s^{(m);i}\right)}{1-\theta}.
\end{equation}
In light of \eqref{eq:5.21}, it follows from \eqref{eq:2.15} with $\bar k=0$ that $\ass$, for each $i=1,\cdots,l$,
$$
\begin{array}{l}
\Dis {\bf 1}_{\{\delta_{\theta} Y_s^{(m,p);i}>0\}}\delta_{\theta} g^{(m,p);i}(s)\vspace{0.1cm}\\
\ \ \leq \Dis \alpha_s+\beta |\delta_{\theta}Y_s^{(m-1,p)}(\delta_{\theta} Y_s^{(m,p);[1:l]})|+\gamma |\delta_{\theta} Z_s^{(m,p);[1:l]}|\vspace{0.1cm}\\
\ \ \ \ \ \Dis +k\left(|Y_s^{(m-1)}(Y_s^{(m);[1:l]})|+
|Y_s^{(m+p-1)}(Y_s^{(m+p);[1:l]})|\right)\vspace{0.1cm}\\
\ \ \leq \Dis \alpha_s+\beta |\delta_{\theta}Y_s^{(m-1,p)}|+H^{(m,p)}(s)+\beta |\delta_{\theta} Y_s^{(m,p);[1:l]}|+\gamma |\delta_{\theta} Z_s^{(m,p);[1:l]}|,
\end{array}
$$
and then
\begin{equation}\label{eq:5.23}
\begin{array}{l}
\Dis \left\langle \frac{\delta_{\theta} Y_s^{(m,p);[1:l]}}{|\delta_{\theta} Y_s^{(m,p);[1:l]}|}{\bf 1}_{|\delta_{\theta} Y_s^{(m,p);[1:l]}|>0}, \delta_{\theta} g^{(m,p);[1:l]}(s)\right\rangle\vspace{0.1cm}\\
\ \ \leq \Dis \sqrt{l}\alpha_s+\sqrt{l}\beta |\delta_{\theta}Y_s^{(m-1,p)}|+H^{(m,p)}(s)+\sqrt{l}\beta |\delta_{\theta} Y_s^{(m,p);[1:l]}|+\sqrt{l}\gamma |\delta_{\theta} Z_s^{(m,p);[1:l]}|
\end{array}
\end{equation}
with
$$
H^{(m,p)}(s):=\sqrt{l}k\left(|Y_s^{(m-1)}|+|Y_s^{(m)}|+|Y_s^{(m+p-1)}|+|Y_s^{(m+p)}|\right).
$$
Furthermore, by \eqref{eq:2.1} and \eqref{eq:5.12} along with H\"{o}lder's inequality we can obtain the existence of a constant $K(q)>0$ depending only on $(q,n,\beta,\gamma,\bar\gamma,\delta,k,\bar k,T)$ such that for each $q>1$,
\begin{equation}\label{eq:5.24}
\begin{array}{l}
\Dis \left(\frac{4eq}{4q-1}\right)^{8q} \sup_{m,p\geq 1}\E\left[\exp\left\{qc_0 e^{c_0 T} \left(|\xi|+\sup_{s\in\T}\left(|Y_s^{(m)}|+|Y_s^{(m+p)}|\right)
\right.\right.\right.\vspace{0.1cm}\\
\hspace{5.5cm} \Dis +\left.\left.\left.\int_0^T \left(\alpha_s+H^{(m,p)}(s)\right) {\rm d}s\right)\right\}\right]\leq K(q).
\end{array}
\end{equation}
Thus, in light of \eqref{eq:5.20}, \eqref{eq:5.21} and \eqref{eq:5.23}, by an identical argument to obtain \eqref{eq:5.15} we can deduce that for each $q>1$,
\begin{equation}\label{eq:5.25}
\begin{array}{l}
\Dis \E\left[\exp\left\{q \sup_{t\in\T}\left(\gamma e^{c_0 t}\sqrt{\sum_{i=1}^{l}|\delta_{\theta} Y_t^{(m,p);i}|^2}\right)\right\}\right]\vspace{0.2cm}\\
\ \ \Dis \leq \left(\frac{eq}{q-1}\right)^{q} \E\left[\exp \left\{\frac{q c_0 e^{c_0 T}}{8}\left(|\xi|+\int_0^T (\alpha_s+H^{(m,p)}(s)){\rm d}s\right)+\frac{q}{8}\sup_{t\in \T}(\gamma e^{c_0 t}|\delta_{\theta} Y_t^{(m-1,p)}|)\right\}\right].
\end{array}
\end{equation}
Similarly, by letting
$$
\delta_{\theta}\widetilde Y^{(m,p)}:=\frac{Y^{(m)}-\theta Y^{(m+p)}}{1-\theta}\ \ \ {\rm and}\ \ \ \delta_{\theta}\widetilde Z^{(m,p)}:=\frac{Z^{(m)}-\theta Z^{(m+p)}}{1-\theta},
$$
it follows that for each $q>1$,
\begin{equation}\label{eq:5.26}
\begin{array}{l}
\Dis \E\left[\exp\left\{q \sup_{t\in\T}\left(\gamma e^{c_0 t}\sqrt{\sum_{i=1}^{l}|\delta_{\theta} \widetilde Y_t^{(m,p);i}|^2}\right)\right\}\right]\vspace{0.2cm}\\
\ \ \Dis \leq \left(\frac{eq}{q-1}\right)^{q} \E\left[\exp \left\{\frac{q c_0 e^{c_0 T}}{8}\left(|\xi|+\int_0^T (\alpha_s+H^{(m,p)}(s)){\rm d}s\right)+\frac{q}{8}\sup_{t\in \T}(\gamma e^{c_0 t}|\delta_{\theta} \widetilde Y_t^{(m-1,p)}|)\right\}\right].\vspace{0.1cm}
\end{array}
\end{equation}
Combining \eqref{eq:5.25}-\eqref{eq:5.26}, by H\"{o}lder's inequality we derive that for each $q>1$,
$$
\begin{array}{l}
\Dis \E\left[\exp\left\{q \sup_{t\in\T}\left(\gamma e^{c_0 t}\sqrt{\sum_{i=1}^{l}|\delta_{\theta} Y_t^{(m,p);i}|^2} +\gamma e^{c_0 t}\sqrt{\sum_{i=1}^{l}|\delta_{\theta} \widetilde Y_t^{(m,p);i}|^2} \right)\right\}\right]\vspace{0.2cm}\\
\ \ \Dis \leq \left(\frac{2eq}{2q-1}\right)^{2q} \E\left[\exp \left\{\frac{q c_0 e^{c_0 T}}{4}\left(|\xi|+\int_0^T (\alpha_s+H^{(m,p)}(s)){\rm d}s\right)\right.\right.\vspace{0.2cm}\\
\hspace*{4.2cm} \Dis \left.\left.+\frac{q}{4}\sup_{t\in \T}\left[\gamma e^{c_0 t}\left(|\delta_{\theta} Y_t^{(m-1,p)}|+|\delta_{\theta} \widetilde Y_t^{(m-1,p)}|\right)\right]\right\}\right]\vspace{0.2cm}\\
\ \ \leq \Dis \left[K\left({q\over 2}\right)\right]^{1\over 2}\left(\E\left[\exp\left\{{q\over 2}\sup_{t\in \T}\left[\gamma e^{c_0 t}\left(|\delta_{\theta}Y_t^{(m-1,p)}|+|\delta_{\theta}\widetilde Y_t^{(m-1,p)}| \right)\right]\right\}\right]\right)^{1\over 2}. \vspace{0.1cm}
\end{array}
$$
On the other hand, in light of \eqref{eq:5.20} and \eqref{eq:5.22}, by an identical argument to obtain \eqref{eq:3.20} we can deduce that for each $q>1$,
$$
\begin{array}{l}
\Dis \E\left[\exp\left\{q \sup_{t\in\T}\left(\gamma e^{c_0 t}\sqrt{\sum_{i=l+1}^{n}|\delta_{\theta} Y_t^{(m,p);i}|^2} +\gamma e^{c_0 t}\sqrt{\sum_{i=l+1}^{n}|\delta_{\theta} \widetilde Y_t^{(m,p);i}|^2} \right)\right\}\right]\vspace{0.1cm}\\
\ \ \leq \Dis \left[K\left({q\over 2}\right)\right]^{1\over 2} \left(\E\left[\exp\left\{{q\over 2}\sup_{t\in \T}\left[\gamma e^{c_0 t}\left(|\delta_{\theta}Y_t^{(m-1,p)}|+|\delta_{\theta}\widetilde Y_t^{(m-1,p)}| \right)\right]\right\}\right]\right)^{1\over 2}.\vspace{0.1cm}
\end{array}
$$
By combining the last two inequalities and H\"{o}lder's inequality we deduce that for each $q>1$,
\begin{equation}\label{eq:5.27}
\begin{array}{l}
\Dis \E\left[\exp\left\{q\sup_{t\in \T}\left[\gamma e^{c_0 t}\left(|\delta_{\theta}Y_t^{(m-1,p)}|+|\delta_{\theta}\widetilde Y_t^{(m-1,p)}| \right)\right]\right\}\right]\\
\ \ \leq \Dis [K(q)]^{1\over 2}\left(\E\left[\exp\left\{q\sup_{t\in \T}\left[\gamma e^{c_0 t}\left(|\delta_{\theta}Y_t^{(m-1,p)}|+|\delta_{\theta}\widetilde Y_t^{(m-1,p)}| \right)\right]\right\}\right]\right)^{1\over 2}. \vspace{0.1cm}
\end{array}
\end{equation}
Based on the last inequality, by induction we get that
$$
\begin{array}{l}
\Dis \E\left[\exp\left\{q \sup_{t\in [0,T]}\left[\gamma e^{c_0 t}\left(|\delta_{\theta} Y_t^{(m,p)}|+|\delta_{\theta} \widetilde Y_t^{(m,p)}|\right)\right]\right\}\right]\\
\ \ \leq \Dis K(q) \left(\E\left[\exp\left\{q \sup_{t\in [0,T]}\left[\gamma e^{c_0 t}\left(|\delta_{\theta} Y_t^{(1,p)}|+|\delta_{\theta} \widetilde Y_t^{(1,p)}|\right)\right]\right\}\right]\right)^{\frac{1}{2^m}},
\end{array}
$$
from which along with \eqref{eq:5.12} it follows that for each $q>1$ and $\theta\in (0,1)$,
$$
\limsup_{m\To\infty}\sup_{p\geq 1}\E\left[\exp\left(q\gamma \sup_{t\in \T}\left(|\delta_{\theta} Y_t^{(m,p)}|+|\delta_{\theta} \widetilde Y_t^{(m,p)}|\right)\right)\right]\leq K(q).
$$
With the last inequality on the hand along with \eqref{eq:2.1}, \eqref{eq:2.3}, \eqref{eq:2.14-1}, \eqref{eq:5.12} and \eqref{eq:5.14}, by using a similar analysis as those from \eqref{eq:3.21} to \eqref{eq:3.26} we can conclude that $\{(Y^{(m)},Z^{(m)})\}_{m=1}^\infty$ is a Cauchy sequence in the space $\s^q(\R^n)\times \hcal^q(\R^{n\times d})$ for each $q>1$, and then converges to a pair of processes $(Y,Z)$ in $\mathcal{E}(\R^n)\times \mcal(\R^{n\times d})$, which is a desired solution of BSDE \eqref{eq:1.1}.\vspace{0.2cm}

It remains to prove the uniqueness assertion of \cref{thm:2.7}. For this, let $(\widetilde Y,\widetilde Z)$ be another solution of system of BSDEs \eqref{eq:1.1} in the space $\mathcal{E}(\R^n)\times\mcal(\R^{n\times d})$, and define for $\theta\in (0,1)$,
$$
\delta_{\theta}U:=\frac{Y-\theta \widetilde Y}{1-\theta},\ \ \ \ \ \ \delta_{\theta}V:=\frac{Z-\theta \widetilde Z}{1-\theta}
$$
and
$$
\delta_{\theta}\widetilde U:=\frac{\widetilde Y-\theta Y}{1-\theta},\ \ \ \ \ \ \delta_{\theta}\widetilde V:=\frac{\widetilde Z-\theta Z}{1-\theta}.\vspace{0.2cm}
$$
First, by virtue of \eqref{eq:2.1} and H\"{o}lder's inequality we can obtain the existence of a positive constant $\bar K(q)$ depending only on $(q,n,\beta,\gamma,\bar\gamma,\delta,k,\bar k,T)$ such that for each $q>1$,\vspace{0.1cm}
$$
\Dis \left(\frac{4eq}{4q-1}\right)^{8q} \E\left[\exp\left\{qc_0 e^{c_0 T} \left[|\xi|+\sup_{s\in\T}\left(|Y_s|+|\widetilde Y_s|\right) +\int_0^T \left(\alpha_s+|Y_s|+|\widetilde Y_s|\right){\rm d}s \right]\right\}\right]\leq \bar K(q).
$$
Next, in light of the last inequality, by using a similar argument as that from \eqref{eq:5.20} to \eqref{eq:5.27} we can derive that for each $q>1$,
$$
\begin{array}{l}
\Dis \E\left[\exp\left\{q\sup_{t\in [0,T]}\left[\gamma e^{c_0 t}\left(|\delta_{\theta} U_t|+|\delta_{\theta} \widetilde U_t| \right)\right]\right\}\right]\vspace{0.1cm}\\
\ \ \leq \Dis \left(\bar K(q)\right)^{1\over 2}\left(\E\left[\exp\left\{q\sup_{t\in [0,T]}\left[\gamma e^{c_0 t}\left(|\delta_{\theta} U_t|+|\delta_{\theta} \widetilde U_t| \right)\right]\right\}\right]\right)^{1\over 2},
\end{array}
$$
and then
$$
\E\left[\exp\left\{q\sup_{t\in [0,T]}\left[\gamma e^{c_0 t}\left(|\delta_{\theta} U_t|+|\delta_{\theta} \widetilde U_t| \right)\right]\right\}\right]\leq \bar K(q).\vspace{0.2cm}
$$
Finally, by picking $q=2$ in the last inequality we have
$$
\E\left[2\gamma \sup_{t\in [0,T]}|\delta_{\theta} U_t|\right] \leq \E\left[\exp\left(2\gamma \sup_{t\in [0,T]}|\delta_{\theta} U_t|\right)\right]\leq \bar K(2),\vspace{0.2cm}
$$
and then
$$
\begin{array}{l}
\Dis \E\left[\sup_{t\in [0,T]}|Y_t-\widetilde Y_t|\right]\leq \Dis \E\left[\sup_{t\in [0,T]}|Y_t-\theta \widetilde Y_t|\right]+(1-\theta)\E\left[\sup_{t\in [0,T]}|\widetilde Y_t|\right]\vspace{0.2cm}\\
\hspace{3cm}\leq \Dis (1-\theta)\left(\frac{\bar K(2)}{2\gamma}+\E\left[\exp\left(\sup_{t\in [0,T]}|\widetilde Y_t|\right)\right]\right),\vspace{0.2cm}
\end{array}
$$
in which sending $\theta\To 1$ leads to that $Y=\widetilde Y$ and then $Z=\widetilde Z$. The proof of \cref{thm:2.7} is then complete.\vspace{0.2cm}

\subsection{The proof of \cref{thm:2.11}\vspace{0.2cm}}
\label{subsec3.4-proof-Th2.11}

To prove \cref{thm:2.11}, we need \cref{pro:3.1} below, which proposes a uniform estimate of the solution. Notably, the whole idea of \cref{pro:3.1} is mainly enlightened by the three authors' previous \cite[Proposition 4.1]{FanHuTang2025EJP}, but a gap in the proof of the latter one is found and filled in here.

\begin{pro}\label{pro:3.1}
Let $\xi\in L^{\infty}(\R^n)$, $\alpha\in \ecal^{\infty}(p\gamma\exp(\beta T))$ for some real $p>1$ and the generator $g$ satisfy assumption \ref{A:H4}. Assume that for some real $h\in (0,T]$, BSDE \eqref{eq:1.1} has a solution $(Y,Z)\in \s^{\infty}_{[T-h,T]}(\R^n)\times {\rm BMO}_{[T-h,T]}(\R^{n\times d})$ on the time interval $[T-h,T]$. Then, there exists constants $\widetilde K>0$ and $\theta_0>0$ depending only on parameters $(\|\xi\|_{\infty}, \|\alpha\|_{\ecal^{\infty}(p\gamma\exp(\beta T))}, \|\bar\alpha\|_{\mcal^{\infty}}, \|\tilde\alpha\|_{\lcal^{\infty}}, \|v\|_{{\rm BMO}}, n,\beta,\gamma,\bar\gamma$, $\lambda,\bar\lambda,c,\bar c,\delta,T,p)$ and being independent of $h$ such that when $\theta\in [0,\theta_0]$, we have
$$
\|Y\|_{\s^{\infty}_{[T-h,T]}}+\|Z\|_{{\rm BMO}_{[T-h,T]}}^2 \leq \widetilde K.\vspace{0.1cm}
$$
Moreover, the preceding assertion remains true for $p=1$ when $\lambda=0$ and $\theta=0$.\vspace{0.1cm}
\end{pro}

\begin{proof}
We only prove the case of $p>1$. The other case can be proved identically. \vspace{0.2cm}

First, we define an increasing non-negative real-valued function $\Phi:[0,+\infty)\To [0,+\infty)$ as follows:
$$
\begin{array}{l}
\Phi(x):=\Dis {2(1+\beta T)+\bar\gamma\over \bar\gamma}\exp\left(2\beta T\right)\vspace{0.2cm}\\
\hspace{2cm} \Dis \times \bigg[2n\bar\lambda T^{1-\delta\over 2}+{\ln 2\over \gamma\exp(\beta T) }+\|\alpha\|_{\ecal^{\infty}(p\gamma \exp(\beta T))}+2\left\|\bar\alpha\right\|_{\mcal^{\infty}}+3x\bigg]\vspace{0.2cm}\\
\hspace{1.5cm} \Dis +{4(\gamma +1) \over \gamma^2}\exp\left\{4\gamma \exp(\beta T)\left(\|\xi\|_{\infty}+\|\tilde\alpha\|_{\lcal^{\infty}}\right)\right\}\vspace{0.2cm}\\
\hspace{2cm} \Dis \times \bigg[1+\beta T \exp(\beta T)\left(\|\xi\|_{\infty}+
\|\tilde\alpha\|_{\lcal^{\infty}}\right)+5n\lambda^2T^{1-\delta}
+\|\tilde\alpha\|_{\lcal^{\infty}}+5\|v\|^2_{{\rm BMO}}\bigg]\vspace{0.2cm}\\
\hspace{1.5cm} +1+c_0 \exp\left(2\beta T+ 2\bar\lambda^2 T\right)\left(x^2+6\|\bar\alpha\|_{\mcal^\infty}^2\right),
\end{array}
$$
and define the following uniform constants: $\bar C_1:= \Phi(\|\xi\|_{\infty})$,
$$
\begin{array}{l}
\bar C_2 := \Dis {6\beta(1+\beta T)+3\beta\bar\gamma\over \bar\gamma} \exp\left(2\beta T\right)\\[4mm]
\hspace{1cm}\Dis +{20\beta^2(\gamma +1) \over \gamma^2}\exp\left\{4\gamma \exp(\beta T)\left(\|\xi\|_{\infty}+
\|\tilde\alpha\|_{\lcal^{\infty}}\right)\right\}+6\beta^2c_0 \exp\left(2\beta T+2\bar\lambda^2 T\right),
\end{array}
$$
$$
\begin{array}{l}
\bar C_3 :=\Dis {2(1+\beta T)+\bar\gamma\over \bar\gamma}\exp\left(2\beta T\right)\left(C_{p,n,\beta,\gamma,\lambda,T}+2n\bar\lambda \right)\\[4mm]
\hspace{1cm}\Dis +{20n\lambda^2(\gamma +1) \over \gamma^2}\exp\left\{4\gamma \exp(\beta T)\left(\|\xi\|_{\infty}+
\|\tilde\alpha\|_{\lcal^{\infty}}\right)\right\}+6n\bar\lambda^2c_0 \exp\left(2\beta T+2\bar\lambda^2 T\right),
\end{array}
$$
$$
\bar C_4 :=\Dis {4\bar c(1+\beta T)+2\bar c\bar\gamma\over \bar\gamma}\exp\left(2\beta T\right)+{20nc^2(\gamma +1) \over \gamma^2}\exp\left\{4\gamma \exp(\beta T)\left(\|\xi\|_{\infty}+
\|\tilde\alpha\|_{\lcal^{\infty}}\right)\right\}
$$
and
$$
\bar C_5 :=\Dis {4(1+\beta T)+2\bar\gamma\over \bar\gamma}\exp\left(2\beta T\right)+{20n(\gamma +1) \over \gamma^2}\exp\left\{4\gamma \exp(\beta T)\left(\|\xi\|_{\infty}+
\|\tilde\alpha\|_{\lcal^{\infty}}\right)\right\},
$$
where $c_0>0$ is a uniform constant defined in \cite[(i) of Proposition A.3]{FanHuTang2025EJP}, and
\begin{equation}\label{eq:6.1}
C_{p,n,\beta,\gamma,\lambda,T}:=
{n\gamma p\lambda^2\over 2(p-1)}\exp(\beta T).
\end{equation}

Since $(Y,Z)\in \s^{\infty}_{[T-h,T]}\times {\rm BMO}_{[T-h,T]}$ is a solution of BSDE \eqref{eq:1.1} on the time interval $[T-h,T]$, we know that for each $i=1,\cdots,n$, $(Y_t,Z_t)_{t\in [T-h,T]}$ also solves the following BSDE:
\begin{equation}\label{eq:6.2}
Y_t^i=\xi^i+\int_t^T f^i(s,Y_s^i,Z_s^i) {\rm d}s-\int_t^T Z_s^i {\rm d}B_s, \ \ t\in [T-h,T],
\end{equation}
where for each $(\omega,t,y,z)\in \Omega\times [T-h,T]\times\R\times\R^{1\times d}$,
\begin{equation}\label{eq:6.3}
f^i(t,y,z):=g^i(t,Y_t(y;i),Z_t(z;i)).
\end{equation}
In what follows, we will prove that for each $i=1,\cdots,n$ and $t\in [T-h,T]$, if $T-t\leq 1$ and $\theta\leq 1$, then
\begin{equation}\label{eq:6.5}
\begin{array}{l}
\Dis \|Y^i\|_{\s^{\infty}_{[T-t,T]}}+\|Z^i\|_{{\rm BMO}_{[T-t,T]}}^2\vspace{0.2cm}\\
\ \ \leq \Dis \bar C_1+\bar C_2 \|Y\|_{\s^{\infty}_{[T-t,T]}} (T-t)+\bar C_3 \|Z\|_{{\rm BMO}_{[T-t,T]}}^2 (T-t)^{1-\delta\over 2}\vspace{0.1cm}\\
\ \ \ \ \ \ \Dis +\bar C_4 \sum_{j=1}^{i-1} \|Z^j\|_{{\rm BMO}_{[T-t,T]}}^{2}+\theta\bar C_5 \|Z\|_{{\rm BMO}_{[T-t,T]}}^2.\vspace{-0.1cm}
\end{array}
\end{equation}
Indeed, for any fixed $i=1,\cdots,n$, since $g$ satisfies assumption \ref{A:H4}, in light of \cite[Remark 2.9]{FanHuTang2025EJP}, we only need to consider the following three cases.\vspace{0.2cm}

(1) $g^i$ satisfies \ref{A:H4}(i). In this case, it follows from \eqref{eq:6.3} that almost everywhere on $\Omega\times [T-h,T]$, for each $(y,z)\in \R\times \R^{1\times d}$, we have
\begin{equation}\label{eq:6.6}
f^i(t,y,z)\, {\bf 1}_{y>0} \leq  \Dis \check\alpha_t+\beta |y|+\frac{\gamma}{2} |z|^2\quad {\rm and} \quad
f^i(t,y,z)\geq \Dis \frac{\bar \gamma}{2} |z|^2-\dot\alpha_t-\beta |y|,
\end{equation}
where
$$
\check\alpha_t:=\alpha_t+\beta |Y_t| +\lambda \sum_{j\neq i} |Z^j_t|
$$
and
$$
\dot\alpha_t:=\bar\alpha_t+\beta |Y_t|+\bar c\sum_{j=1}^{i-1}|Z^j_t|^2+\sum_{j=i+1}^{n} \left(\bar\lambda |Z^j_t|^{1+\delta}+\theta |Z^j_t|^2\right).\vspace{0.2cm}
$$
Since $\alpha\in \ecal^{\infty}(p\gamma \exp(\beta T))$, $\bar\alpha\in \mcal^{\infty}$ and $|Y|\in \lcal^\infty_{[T-h,T]}$, applying \cite[(i) of Lemma 3.1]{FanHuTang2025EJP} with
$$
\theta=0,\ r=\gamma \exp(\beta T),\ \bar\delta=0,\ \bar t=T-h,\ V=\bar V=Z,\ \breve\alpha=\check\alpha\ \ {\rm and}\ \ \hat\alpha=\alpha+\beta |Y|
$$
and \cite[(ii) of  Lemma 3.1]{FanHuTang2025EJP} with
$$\bar\delta=\delta,\ \bar t=T-h,\ V=\bar V=Z,\ \breve\alpha=\dot\alpha\ \ {\rm and}\ \ \hat\alpha=\bar\alpha+\beta |Y|\vspace{-0.1cm}
$$
we can deduce that
\begin{equation}\label{eq:6.7}
\begin{array}{lll}
\Dis \|\check\alpha\|_{\ecal^{\infty}_{[t,T]}(\gamma \exp(\beta T))} &\leq &\Dis \|\alpha\|_{\ecal^{\infty}_{[t,T]}(p\gamma \exp(\beta T))}+\beta \|Y\|_{\s^\infty_{[t,T]}}(T-t)+{\ln 2\over \gamma \exp(\beta T)}\vspace{0.2cm}\\
&&\Dis +C_{p,n,\beta,\gamma,\lambda,T}\|Z\|_{{\rm BMO}_{[t,T]}}^2 (T-t)<+\infty,\ \ t\in [T-h,T]
\end{array}
\end{equation}
and
\begin{equation}\label{eq:6.8}
\begin{array}{lll}
\Dis \left\|\dot\alpha\right\|_{\mcal^{\infty}_{[t,T]}} &\leq &\Dis \left\|\bar\alpha\right\|_{\mcal^{\infty}_{[t,T]}}+\beta \|Y\|_{\s^\infty_{[t,T]}} (T-t)+n\bar\lambda \|Z\|_{{\rm BMO}_{[t,T]}}^{1+\delta}(T-t)^{1-\delta\over 2}\vspace{0.2cm}\\
&&\Dis +\bar c\sum_{j=1}^{i-1}\|Z^j\|^2_{{\rm BMO}_{[t,T]}}+\theta \|Z\|^2_{{\rm BMO}_{[t,T]}}
<+\infty,\ \ t\in [T-h,T],
\end{array}
\end{equation}
where the constant $C_{p,n,\beta,\gamma,\lambda,T}$ is defined in \eqref{eq:6.1}. Combining \eqref{eq:6.6}, \eqref{eq:6.7} and \eqref{eq:6.8} implies that on $\Omega\times [T-h,T]$, the generator $f^i$ satisfies the second inequality for the case of $y>0$ and the first inequality in \cite[assumption (A1)]{FanHuTang2025EJP} with $\bar\beta=\beta$ and $\varphi(\cdot)\equiv 0$. It then follows from \cite[(i) of Proposition A.1]{FanHuTang2025EJP} that for each $t\in [T-h,T]$ such that $T-t\leq 1$, we have\vspace{0.2cm}
\begin{equation}\label{eq:6.9}
\begin{array}{l}
\Dis \|Y^i\|_{\s^{\infty}_{[t,T]}}+\|Z^i\|_{{\rm BMO}_{[t,T]}}^2 \leq  \Dis {2(1+\beta T)+\bar\gamma\over \bar\gamma}\exp\left(2\beta T\right)\bigg\{2n\bar\lambda T^{1-\delta\over 2}+{\ln 2\over \gamma\exp(\beta T) } \vspace{0.2cm}\\
\Dis \hspace{0.5cm} +3\|\xi^i\|_{\infty}+\|\alpha\|_{\ecal^{\infty}(p\gamma \exp(\beta T))}+2\left\|\bar\alpha\right\|_{\mcal^{\infty}}+3\beta \|Y\|_{\s^\infty_{[t,T]}} (T-t)+2\theta \|Z\|^2_{{\rm BMO}_{[t,T]}}\vspace{0.2cm}\\
\Dis \hspace{0.5cm} +\left(C_{p,n,\beta,\gamma,\lambda,T}+2n\bar\lambda \right) \|Z\|_{{\rm BMO}_{[t,T]}}^2 (T-t)^{1-\delta\over 2}+2\bar c\sum_{j=1}^{i-1}\|Z^j\|^2_{{\rm BMO}_{[t,T]}}\Big\}.
\end{array}
\end{equation}
Then, it follows from \eqref{eq:6.9} along with the definitions of constants $\bar C_1$, $\bar C_2$, $\bar C_3$, $\bar C_4$ and $\bar C_5$ that the desired inequality \eqref{eq:6.5} is fulfilled in this case when $T-t\leq 1$.\vspace{0.2cm}

(2) $g^i$ satisfies \ref{A:H4}(ii). For this case, it follows from \eqref{eq:6.3} that almost everywhere on $\Omega\times [T-h,T]$, for each $(y,z)\in \R\times \R^{1\times d}$, we have
$$
f^i(t,y,z){\rm sgn}(y) \leq \tilde\alpha_t+\beta |y|+\bar u_t |z|+\frac{\gamma}{2} |z|^2,
$$
where
$$
\bar u:=v+\beta \sqrt{|Y|}+c\sum_{j=1}^{i-1}|Z^j| + \sum_{j=i+1}^{n}\big(\lambda  |Z^j|^\delta+\theta |Z^j|\big)
\in {\rm BMO}_{[T-h,T]}(\R).\vspace{0.1cm}
$$
This means that on $\Omega\times [T-h,T]$, the generator $f^i$ satisfies the first inequality for the case of $y<0$ and the second inequality for the case of $y>0$ in \cite[assumption (A3)]{FanHuTang2025EJP} with $\ddot\alpha=\tilde\alpha$ and $\bar\beta=\beta$. Then, in light of the fact that for each $t\in [T-h,T]$, $\tau\in \mathcal{T}_{[t,T]}$ and $j=i+1,\cdots,n$, it holds that
$$
\E_\tau\left[\int_\tau^T|Z^j_s|^{2\delta}{\rm d}s\right]\leq \left(\E_\tau\left[\int_\tau^T|Z^j_s|^2{\rm d}s\right]\right)^\delta (T-t)^{1-\delta}\leq \big(1+\|Z\|^2_{{\rm BMO}[t,T]}\big)(T-t)^{1-\delta},
$$
it follows from \cite[(i) of Proposition A.2]{FanHuTang2025EJP} with $r=T$ that for each $t\in [T-h,T]$,
\begin{equation}\label{eq:6.10}
\begin{array}{l}
\Dis\|Y^i\|_{\s^\infty_{[t,T]}}+\|Z^i\|^2_{{\rm BMO}_{[t,T]}}\\[2mm]
\Dis \hspace{0.2cm} \leq  \Dis {4(\gamma +1) \over \gamma^2}\exp\left\{4\gamma \exp(\beta T)\left(\|\xi^i\|_{\infty}+\|\tilde\alpha\|_{\lcal^{\infty}}\right)\right\} \Big(1+5n\lambda^2T^{1-\delta}\\[4mm]
\Dis \hspace{0.6cm} +\|\tilde\alpha\|_{\lcal^{\infty}}+\beta T \exp(\beta T)\left(\|\xi^i\|_{\infty}+\|\tilde\alpha\|_{\lcal^{\infty}}\right)
+5\|v\|^2_{{\rm BMO}_{[t,T]}}+5\beta^2\|Y\|_{\s^\infty_{[t,T]}} (T-t)\\[2mm]
\Dis \hspace{0.6cm} +5nc^2 \sum_{j=1}^{i-1} \|Z^j\|^2_{{\rm BMO}_{[t,T]}}
+5n\theta^2 \|Z\|^2_{{\rm BMO}_{[t,T]}}
+5n\lambda^2 \|Z\|^2_{{\rm BMO}_{[t,T]}} (T-t)^{1-\delta}\Big) .
\end{array}
\end{equation}
Then, it follows from \eqref{eq:6.10} along with the definitions of constants $\bar C_1$, $\bar C_2$, $\bar C_3$, $\bar C_4$ and $\bar C_5$ that the desired inequality \eqref{eq:6.5} is fulfilled in this case when $T-t\leq 1$ and $\theta\leq 1$.\vspace{0.2cm}

(3) $g^i$ satisfies \ref{A:H4}(iii). For this case, it follows from \eqref{eq:6.3} that almost everywhere on $\Omega\times [T-h,T]$, for each $(y,z)\in \R\times \R^{1\times d}$, we have
\begin{equation}\label{eq:6.11}
f^i(t,y,z)\, {\rm sgn}(y) \leq  \Dis \breve\alpha_t+\beta |y|+\bar\lambda |z|,
\end{equation}
where
$$
\breve\alpha_t:=\bar\alpha_t+\beta \sqrt{|Y_t|}+\bar\lambda \sum_{j\neq i} |Z^j_t|.
$$
Applying \cite[(iii) of Lemma 3.1]{FanHuTang2025EJP} with $\theta=0$, $\bar\delta=0$, $\bar t=T-h$, $V=\bar V=Z$ and $\hat\alpha=\bar\alpha$ we can deduce that for each $t\in [T-h,T]$,
\begin{equation}\label{eq:6.12}
\Dis \left\|\breve\alpha\right\|_{\mcal^{\infty}_{[t,T]}} \leq \Dis \left\|\bar\alpha\right\|_{\mcal^{\infty}_{[t,T]}}+\beta \sqrt{\|Y\|_{\s^\infty_{[t,T]}}} (T-t)+\bar\lambda \sqrt{n} \|Z\|_{{\rm BMO}_{[t,T]}} (T-t)^{1\over 2}<+\infty.
\end{equation}
Combining \eqref{eq:6.11} and \eqref{eq:6.12} implies that on $\Omega\times [T-h,T]$, the generator $f^i$ satisfies the first inequality for the case of $y<0$ and the second inequality for the case of $y>0$ in
\cite[assumption (A4)]{FanHuTang2025EJP} with $\dot\alpha=\breve\alpha$ and $\bar\beta=\beta$. It follows from
\cite[(i) of Proposition A.3]{FanHuTang2025EJP} that for each $t\in [T-h,T]$,
\begin{equation}\label{eq:6.13}
\begin{array}{l}
\Dis \|Y^i\|_{\s^\infty_{[t,T]}}+\left\|Z^i\right\|_{{\rm BMO}_{[t,T]}}^2 \leq \Dis 1+c_0 \exp\left(2\beta T+ 2\bar\lambda^2 T\right)\left(\|\xi^i\|_\infty^2
+6\|\bar\alpha\|_{\mcal^\infty}^2\right)\vspace{0.1cm}\\
\hspace{4.6cm}\Dis +6\beta^2 c_0 \exp\left(2\beta T+2\bar\lambda^2 T\right) \|Y\|_{\s^\infty_{[t,T]}}(T-t)^2\vspace{0.1cm}\\
\hspace{4.6cm}\Dis +6n\bar\lambda^2c_0 \exp\left(2\beta T+2\bar\lambda^2 T\right) \|Z\|_{{\rm BMO}_{[t,T]}}^2(T-t).
\end{array}
\end{equation}
Then, it follows from \eqref{eq:6.13} along with the definitions of constants $\bar C_1$, $\bar C_2$, $\bar C_3$, $\bar C_4$ and $\bar C_5$ that the desired inequality \eqref{eq:6.5} is fulfilled in this case when $T-t\leq 1$.\vspace{0.2cm}

Furthermore, in light of \eqref{eq:6.5}, by induction for $i$ it is not very hard to deduce that for each $i=1,\cdots,n$ and each $t\in [T-h,T]$, if $T-t\leq 1$ and $\theta\leq 1$, then
\begin{equation}\label{eq:6.14}
\begin{array}{l}
\hspace*{-0.4cm}\Dis \sum_{j=1}^{i}\left(\|Y^j\|_{\s^{\infty}_{[t,T]}}+\|Z^j\|_{{\rm BMO}_{[t,T]}}^2\right)\\
\hspace*{-0.4cm}\ \ \leq \Dis \left(\bar C_1+\bar C_2 \|Y\|_{\s^{\infty}_{[T-t,T]}} (T-t)+\bar C_3 \|Z\|_{{\rm BMO}_{[T-t,T]}}^2(T-t)^{1-\delta\over 2}+\theta\bar C_5 \|Z\|_{{\rm BMO}_{[T-t,T]}}^2\right)\bar C_6^i,
\end{array}
\end{equation}
where
\begin{equation}\label{eq:6.15}
\bar C_6^1:=1\ \ {\rm and}\ \
\bar C_6^i:=\bar C_6^{i-1}+1+\bar C_4\bar C_6^{i-1}\ {\rm for}\ i=2,\cdots,n,\ {\rm recursively}.
\end{equation}
In particular, letting $i=n$ in \eqref{eq:6.14} implies that
for each $t\in [T-h,T]$, if $T-t\leq 1$ and $\theta\leq 1$, then
$$
\begin{array}{lll}
\Dis \|Y\|_{\s^{\infty}_{[t,T]}}+\|Z\|_{{\rm BMO}_{[t,T]}}^2 &\leq & \Dis \bar C_1\bar C_6^n+\bar C_2 \bar C_6^n\|Y\|_{\s^{\infty}_{[T-t,T]}} (T-t)\\
&&\Dis +\bar C_3 \bar C_6^n\|Z\|_{{\rm BMO}_{[T-t,T]}}^2(T-t)^{1-\delta\over 2}+\theta \bar C_5 \bar C_6^n \|Z\|_{{\rm BMO}_{[T-t,T]}}^2.
\end{array}
$$
Thus, when $\theta\leq 1/(4\bar C_5 \bar C_6^n)$, we have
\begin{equation}\label{eq:6.16}
\|Y\|_{\s^{\infty}_{[T-\eps,T]}}+\|Z\|_{{\rm BMO}_{[T-\eps,T]}}^2 \leq \Dis 2\bar C_1 \bar C_6^n,
\end{equation}
where
$$
\eps:=\min\Bigg\{h,\  1, \ {1\over 2\bar C_2 \bar C_6^n}, \ \Big({1\over 4\bar C_3 \bar C_6^n}\Big)^{2\over 1-\delta}\Bigg\}>0.\vspace{0.3cm}
$$

Finally, for $m\geq 1$, we define successively the following constants:
\begin{equation}\label{eq:6.17}
\bar C_7^1:= 2\bar C_1\bar C_6^n\ \ \ {\rm and}\ \ \ \bar C_7^{m+1}:=\bar C_7^m+2\Phi(\bar C_7^m) \bar C_6^n,
\end{equation}
where the constant $\bar C_6^n$ is defined in \eqref{eq:6.15}. And, we let $m_0$ be the unique positive integer satisfying
$T-h \in [T-m_0\eps,T-(m_0-1)\eps)$, or equivalently,
\begin{equation}\label{eq:6.18}
{h\over \eps}\leq m_0 <{h\over \eps}+1.
\end{equation}
If $m_0=1$, it then follows from \eqref{eq:6.16} and \eqref{eq:6.17} that
$$
\|Y\|_{\s^{\infty}_{[T-h,T]}}+\|Z\|_{{\rm BMO}_{[T-h,T]}}^2 \leq \Dis \bar C_7^1=\bar C_7^{m_0}.
$$
If $m_0=2$, it then follows from \eqref{eq:6.16}  and \eqref{eq:6.17} that
\begin{equation}\label{eq:6.19}
\|Y_{T-\eps}\|_{\infty}\leq \|Y\|_{\s^{\infty}_{[T-\eps,T]}}\leq \bar C_7^1.
\end{equation}
Now, consider the following system of BSDEs
$$
  Y_t=Y_{T-\eps}+\int_t^{T-\eps} g(s,Y_s,Z_s){\rm d}s-\int_t^{T-\eps} Z_s {\rm d}B_s, \ \ t\in [T-h,T-\eps].
$$
In light of \eqref{eq:6.19} and the definition of $\Phi(\cdot)$, by \cite[Propositions A.1, A.2 and A.3]{FanHuTang2025EJP} we can use a similar argument as that obtaining \eqref{eq:6.5} to get that for each $i=1,\cdots,n$, and each $\theta\leq 1$,
$$
\begin{array}{l}
\Dis \|Y^i\|_{\s^{\infty}_{[T-h,T-\eps]}}+\|Z^i\|_{{\rm BMO}_{[T-h,T-\eps]}}^2\\[1mm]
\Dis \ \ \leq \Dis \Phi(\bar C_1^1)+\bar C_2 \|Y\|_{\s^{\infty}_{[T-h,T-\eps]}} \eps +\bar C_3 \|Z\|_{{\rm BMO}_{[T-h,T-\eps]}}^2 \eps^{1-\delta\over 2}\\
\Dis \ \ \ \ \ +\bar C_4 \sum_{j=1}^{i-1} \|Z^j\|_{{\rm BMO}_{[T-h,T-\eps]}}^{2}
+\theta\bar C_5 \|Z\|_{{\rm BMO}_{[T-h,T-\eps]}}^2.
\end{array}
$$
We would like to emphasize that for the proof of the last assertion, we need to use $T-\eps$ instead of $T$ and $Y_{T-\eps}$ instead of $\eta$ in \cite[Propositions A.1 and A.3]{FanHuTang2025EJP}, and use $T-\eps$ instead of $r$ in \cite[Propsotion A.2]{FanHuTang2025EJP}. Moreover, in light of the definition of $\eps$, by induction and a similar argument as that from \eqref{eq:6.14} to \eqref{eq:6.16} we can further get that when $\theta\leq 1/(4\bar C_5 \bar C_6^n)$,
$$
\|Y\|_{\s^{\infty}_{[T-h,T-\eps]}}+\|Z\|_{{\rm BMO}_{[T-h,T-\eps]}}^2 \leq \Dis 2\Phi(\bar C_7^1) \bar C_6^n.
$$
Combining  the last inequality, \eqref{eq:6.16} and \eqref{eq:6.17} implies that for each $\theta\leq 1/(4\bar C_5 \bar C_6^n)$, we have
$$
\|Y\|_{\s^{\infty}_{[T-h,T]}}+\|Z\|_{{\rm BMO}_{[T-h,T]}}^2 \leq \Dis \bar C_7^1+2\Phi(\bar C_7^1) \bar C_6^n = \bar C_7^2=\bar C_7^{m_0},
$$
Proceeding the above computation gives that if $\theta\leq 1/(4\bar C_5 \bar C_6^n)$ and $m_0$ satisfies \eqref{eq:6.18}, then
\begin{equation}\label{eq:6.20}
\|Y\|_{\s^{\infty}_{[T-h,T]}}+\|Z\|_{{\rm BMO}_{[T-h,T]}}^2 \leq \bar C_7^{m_0}=\bar C_7^{\left\lceil {h\over \eps}\right\rceil}\leq \bar C_7^{1+\left\lceil {T\over \eps_0}\right\rceil}=:\widetilde K,
\end{equation}
where the constant $\bar C_7^{m_0}$ is defined in \eqref{eq:6.17}, $\lceil x\rceil$ represents the minimum of integers which is equal to or bigger than $x\in \R$, and
$$
\eps_0:=\min\Bigg\{{1\over 2\bar C_2 \bar C_6^n}, \ \Big({1\over 4\bar C_3 \bar C_6^n}\Big)^{2\over 1-\delta}\Bigg\}\in (0,1).
$$
Then, the desired conclusion follows immediately from \eqref{eq:6.20}. The proof is complete.\vspace{0.2cm}
\end{proof}

Now, we can give the proof of \cref{thm:2.11}.

\begin{proof}[Proof of \cref{thm:2.11}]
With \citet[Theorem 2.3]{FanHuTang2025EJP} and \cref{pro:3.1} in hands, we can closely follow the proof of \citet[Theorem 4.1]{cheridito2015Stochastics} to prove our \cref{thm:2.11}. All the details are omitted here.\vspace{0.2cm}
\end{proof}

\subsection{The proof of \cref{thm:2.13}\vspace{0.2cm}}
\label{subsec3.5-proof-Th2.13}

We would like to mention that the whole idea of \cref{thm:2.13} is mainly enlightened by the authors' \cite[Theorem 2.27]{FanHuTang2025EJP}, but some mistakes in the proof of the latter one are found and corrected here.\vspace{0.2cm}

First of all, we define
\begin{equation}\label{eq:7.21e}
\bar g(t,y,z):=g\Big(t,\ y+\int_0^t H_s {\rm d}B_s,\ z+H_t\Big),\ \ (\omega,t,y,z)\in \Omega\times\T\times \R^n\times \R^{n\times d}.
\end{equation}
Consider the following system of BSDEs
\begin{equation}\label{eq:7.21}
\bar Y_t=\bar \xi+\int_t^T \bar g(s,\bar Y_s,\bar Z_s){\rm d}s-\int_t^T \bar Z_s {\rm d}B_s, \ \ t\in\T,
\end{equation}
or, equivalently,
$$
\bar Y_t^i=\bar\xi^i+\int_t^T \bar g^i(s,\bar Y_s,\bar Z_s) {\rm d}s-\int_t^T \bar Z_s^i {\rm d}B_s, \ \ t\in\T,\ \ i=1,\cdots,n.
$$
It is easy to verify that BSDE \eqref{eq:7.21} admits a unique solution $(\bar Y,\bar Z)\in \s^\infty(\R^n)\times {\rm BMO}(\R^{n\times d})$ if and only if BSDE \eqref{eq:1.1} admits a unique solution $(Y,Z):=(\bar Y+\int_0^\cdot H_s {\rm d}B_s,\ \bar Z+H)$ on $\T$ such that
$$
\Big(Y-\int_0^\cdot H_s {\rm d}B_s,\,  Z\Big)\in \s^\infty(\R^n)\times {\rm BMO}(\R^{n\times d}).
$$
Thus, in light of \cref{thm:2.11}, for completing the proof of \cref{thm:2.13} it suffices to prove that the generator $\bar g$ defined in \eqref{eq:7.21e} also satisfies assumptions \ref{A:H1"} and \ref{A:H4'} with some other appropriate parameters when the constant $\theta$ is smaller than a given constant $\theta_0$ depending only on $(p,\bar p,\gamma)$. It is evident that $\bar g$ satisfies \ref{A:H1"} since $g$ satisfies it. \vspace{0.2cm}

In what follows, we prove that the generator $\bar g$ also satisfies assumption \ref{A:H4'} with some other appropriate parameters under the given conditions. Assume that $\alpha\in \ecal^{\infty}(p\gamma)$ for some real $p>1$, and $|H|^2\in \ecal^{\infty}(2\bar p (q\gamma)^2)$ for some $\bar p>1$ with $q=p/(p-1)$ such that $1/p+1/q=1$. Define
$$
\eps:={p-1\over p+1},\ \ \hat p:={(p+1)(\bar p+1)\over 2(p+\bar p)}>1,\ \ \tilde p:={p+1\over 2\hat p}={p+\bar p\over \bar p+1}>1\ \ {\rm and}\ \ \tilde q:={p+1\over p-2\hat p+1}={p+\bar p\over p-1}>1.
$$
It is evident that $1/\tilde p+1/\tilde q=1$. Since $g$ satisfies assumption \ref{A:H4'}, by an indentical argument to \cite[(iv) of Remark 2.2]{FanHuTang2025EJP}, in assumption \ref{A:H4'} we can without loss of generality suppose that for all $i=1,\cdots,n$, $f=g^i$ satisfies any one of the three conditions \ref{A:H4'}(i), \ref{A:H4'}(ii) and \ref{A:H4'}(iii). Thus, for $i=1,\cdots,n$, we need to address the following three cases.\vspace{0.2cm}

(1) $g^i$ satisfies \ref{A:H4'}(i). For this case, it follows from the definition of $\bar g$ that $\as$, for each $(y,z)\in\R^n\times \R^{n\times d}$, we have
$$
\Dis \bar g^i(\omega,t,y,z) \leq  \Dis \alpha_t(\omega)+\lambda\sum_{j\neq i} |z^j+H_t^j(\omega)|+\frac{\gamma}{2} |z^i+H_t^i(\omega)|^2
$$
and
$$
\begin{array}{l}
\Dis \bar g^i(\omega,t,y,z)\geq \Dis \frac{\bar \gamma}{2} |z^i+H_t^i(\omega)|^2 -\bar\alpha_t(\omega)-\sum_{j=i+1}^{n} \left(\bar\lambda |z^j+H_t^j(\omega)|^{1+\delta}+\theta |z^j+H_t^j(\omega)|^2\right)\\
\hspace{2.4cm} \Dis -\bar c\sum_{j=1}^{i-1}|z^j+H_t^j(\omega)|^2.
\end{array}
$$
Note that for each $a,b\in\R$, we have\vspace{0.1cm}
$$
\begin{array}{l}
\Dis (a+b)^2\leq 2a^2+2b^2,\ \ (a+b)^{1+\delta}\leq 2a^{1+\delta}+2b^{1+\delta},\vspace{0.1cm}\\
\Dis (a+b)^2\leq (1+\eps)a^2+\left(1+{1\over \eps}\right)b^2\ \ {\rm and}\ \ (a+b)^2\geq {1\over 2}a^2-b^2.
\end{array}
$$
It follows that $\as$, for each $(y,z)\in\R^n\times \R^{n\times d}$,
\begin{equation}\label{eq:7.23}
\begin{array}{l}
\Dis \frac{\bar \gamma}{4} |z^i|^2 -\hat\alpha_t(\omega)-\sum_{j=i+1}^{n} \left(2\bar\lambda |z^j|^{1+\delta}+2\theta |z^j|^2\right)-2\bar c\sum_{j=1}^{i-1}|z^j|^2\leq \Dis \bar g^i(\omega,t,y,z)\\
\quad\quad\quad \leq  \Dis \breve\alpha_t(\omega)+\lambda\sum_{j\neq i} |z^j|+\frac{\gamma(1+\eps)}{2} |z^i|^2
\end{array}\vspace{0.1cm}
\end{equation}
with
$$
\breve\alpha_t(\omega):=\alpha_t(\omega)+\lambda\sum_{j\neq i}|H_t^j(\omega)|+\frac{\gamma(1+\eps)}{2\eps} |H_t^i(\omega)|^2
$$
and
$$
\hat\alpha_t(\omega):= \bar\alpha_t(\omega)+2\bar c\sum_{j=1}^{i-1}|H_t^j(\omega)|^2+\frac{\bar \gamma}{2} |H_t^i(\omega)|^2+\sum_{j=i+1}^{n} \left(2\bar\lambda |H_t^j(\omega)|^{1+\delta}+2\theta |H_t^j(\omega)|^2\right).
$$
By virtue of  $\bar\alpha\in \mcal^\infty$ and $|H|^2\in \ecal^{\infty}(2\bar p (q\gamma)^2)$ along with the fact that
$$
\hat\alpha\leq 2n\bar\lambda+\bar\alpha+\left(2\bar c+{\bar\gamma\over 2}+2\bar\lambda+2\theta\right) |H|^2,\vspace{-0.2cm}
$$
by \cite[(i) of Remark 2.1]{FanHuTang2025EJP} we can deduce that\vspace{-0.1cm}
\begin{equation}\label{eq:7.24}
\hat\alpha \in \mcal^\infty.
\end{equation}
On the other hand, from the definitions of constants $q, \hat p, \tilde p, \tilde q,\eps$ and $\breve\alpha$ it is not very hard to verify that $p=\tilde p \hat p (1+\eps)$, $\tilde q\hat p(1+\eps)=(\bar p+1)q<2\bar pq$ and there exists a positive constant $\bar\eps>0$ such that
$$
\breve\alpha\leq \alpha+{n\lambda^2\over 4\bar\eps}+(\bar \eps+q\gamma)|H|^2\quad {\rm and} \quad (\bar\eps+q\gamma)\tilde q \hat p\gamma (1+\eps)=2\bar p(q\gamma)^2.
$$
Then, it follows from the integrability condition of $\alpha$ and $|H|^2$ that $\alpha\in \ecal^{\infty}(\tilde p \hat p\gamma (1+\eps))$ and $(\bar \eps+q\gamma) |H|^2\in \ecal^{\infty}(\tilde q \hat p\gamma (1+\eps))$. Since $1/\tilde p+1/\tilde q=1$, by \cite[(iii) of Remark 2.1]{FanHuTang2025EJP} we have
\begin{equation}\label{eq:7.25}
\breve\alpha\in \ecal^{\infty}(\hat p\gamma (1+\eps)).
\end{equation}
Combining \eqref{eq:7.23}, \eqref{eq:7.24} and \eqref{eq:7.25} implies that the generator $\bar g^i$ satisfies \ref{A:H4'}(i) with parameters $(\breve\alpha, \hat\alpha,\gamma(1+\eps),\bar\gamma/2,2\bar c,\lambda,2\bar\lambda,2\theta)$ instead of $(\alpha, \bar\alpha,\gamma,\bar\gamma,\bar c,\lambda,\bar\lambda,\theta)$.\vspace{0.2cm}

(2) $g^i$ satisfies \ref{A:H4'}(ii) and violates \ref{A:H4'} (i) and (iii), and $-g^i$ violates \ref{A:H4'} (i). Since $H^i\equiv 0$ in this case, it follows from the definition of $\bar g$ that $\as$, for each $(y,z)\in\R^n\times \R^{n\times d}$, we have\vspace{0.1cm}
$$
|\bar g^i(\omega,t,y,z)|\leq \tilde\alpha_t(\omega)+\frac{\gamma}{2} |z^i|^2+\Big\{\bar v_t(\omega)+c\sum_{j=1}^{i-1}|z^j| +\sum_{j=i+1}^{n} \big(\lambda |z^j|^\delta+\theta |z^j|\big)\Big\}|z^i|
$$
with
$$
\bar v_t(\omega):=v_t(\omega)+n\lambda +(c+\lambda+\theta)\sum_{j\neq i}|H_t^j(\omega)|
\leq v_t(\omega)+n\lambda+(c+\lambda+\theta)\sqrt{n}|H_t(\omega)|.
$$
And, it follows from $v\in {\rm BMO}(\R)$ and $H\in {\rm BMO}(\R^{n\times d})$ that $\bar v\in {\rm BMO}(\R)$, which indicates that $\bar g^i$ satisfies \ref{A:H4'}(ii) with
parameters $(\tilde\alpha,\bar v,c,\gamma,\lambda,\theta)$ instead of $(\tilde\alpha, v,c,\gamma,\lambda,\theta)$.\vspace{0.2cm}

(3) $g^i$ satisfies \ref{A:H4'}(iii). In this case, it follows from the definition of $\bar g$ that $\as$, for each $(y,z)\in\R^n\times \R^{n\times d}$, we have\vspace{-0.1cm}
$$
|g^i(\omega,t,y,z)|\leq \check\alpha_t(\omega)+\bar\lambda |z|\vspace{-0.1cm}
$$
with\vspace{-0.2cm}
$$
\check\alpha_t(\omega):=\bar\alpha_t(\omega)+\bar\lambda |H_t(\omega)|\leq \bar\lambda+\bar\alpha_t(\omega)+\bar\lambda |H_t(\omega)|^2.
$$
And, it follows from $\bar\alpha\in \mcal^\infty$ and $H\in {\rm BMO}(\R^{n\times d})$ that $\check\alpha\in \mcal^\infty$, which indicates that $\bar g^i$ satisfies \ref{A:H4'}(iii) with
parameters $(\check\alpha,\bar\lambda)$ instead of $(\bar\alpha,\bar\lambda)$.\vspace{0.2cm}

Up to now, we have proved the desired conclusion. The proof of \cref{thm:2.13} is then complete.\vspace{0.2cm}

\appendix
\section{Auxiliary results on one-dimensional quadratic BSDEs}
\renewcommand{\appendixname}{}

We consider the following one-dimensional BSDE:
\begin{equation}\label{eq:A.1}
  Y_t=\eta+\int_t^T f(s,Y_s,Z_s){\rm d}s-\int_t^T Z_s {\rm d}B_s, \ \ t\in\T,
\end{equation}
where the terminal value $\eta\in \F_T$ is a real-valued random variable, the generator $f(\omega, t,y,z):\Omega\times\T\times\R\times\R^{1\times d}\To \R$
is $(\F_t)$-progressively measurable for each $(y,z)\in \R\times\R^{1\times d}$, and the solution $(Y_t,Z_t)_{t\in\T}$ is a pair of $(\F_t)$-progressively measurable processes with values in $\R\times\R^{1\times d}$ which satisfies \eqref{eq:A.1}.\vspace{0.1cm}

We always assume that $(\hat\alpha_t)_{t\in\T}$ is an $(\F_t)$-progressively measurable nonnegative real-valued process, and $\varphi_t(\omega,r):\Omega\times \T\times [0,+\infty)\mapsto  [0,+\infty)$ is a nonnegative real-valued function satisfying the following two conditions:
\begin{itemize}
\item For each $r\in [0,+\infty)$, $\varphi_t(\omega,r)$ is $(\F_t)$-progressively measurable and
$$\ps,\ \ \int_0^T\varphi_t(\omega,r){\rm d}t<+\infty;\vspace{-0.1cm}$$
\item Almost everywhere in $\Omega\times \T$, the function $r\mapsto \varphi_t(\omega,r)$ is increasing,
\end{itemize}
and that for each $p\geq 1$,
\begin{equation}\label{eq:A.2}
\E\left[\exp\left(p|\eta|+p\int_0^T \hat\alpha_s{\rm d}s\right)\right]<+\infty.\vspace{0.2cm}
\end{equation}

Let us further introduce the following assumptions on the generator $f$.

\begin{enumerate}
\renewcommand{\theenumi}{(A\arabic{enumi})}
\renewcommand{\labelenumi}{\theenumi}

\item\label{A:A1} Almost everywhere in $\Omega\times\T$, $f(\omega,t,\cdot,\cdot)$ is continuous.

\item\label{A:A2} Almost everywhere in $\Omega\times\T$, for each $(y,z)\in \R\times\R^{1\times d}$, we have
\begin{equation}\label{eq:A.3}
    {\rm sgn}(y)f(\omega,t,y,z)\leq \hat\alpha_t(\omega)+\beta|y|+\frac{\gamma}{2} |z|^2
\end{equation}
and
\begin{equation}\label{eq:A.4}
|f(\omega,t,y,z)|\leq \hat\alpha_t(\omega)+\varphi_t(\omega,|y|)+\frac{\gamma}{2} |z|^2.\vspace{0.1cm}
\end{equation}
Furthermore, almost everywhere in $\Omega\times\T$, $f$ satisfies the following inequality
\begin{equation}\label{eq:A.5}
  \begin{array}{l}
   \hspace*{-0.4cm}\Dis {\bf 1}_{\{y_1-\theta y_2>0\}}\left(f(\omega,t,y_1,z_1)-\theta f(\omega,t,y_2,z_2)\right)\vspace{0.1cm}\\
   \hspace*{-0.4cm}\ \ \leq \Dis (1-\theta)\left(\hat\alpha_t(\omega)+ \beta \left|\frac{y_1-\theta y_2}{1-\theta}\right|+\frac{\gamma}{2} \left|\frac{z_1-\theta z_2}{1-\theta}\right|^2+k(|y_1|+|y_2|)+\bar k (|z_1|^{1+\delta}+|z_2|^{1+\delta})\right)\\
   \Dis \hspace{2.2cm}\RE (y_1,y_2,z_1,z_2)\in \R\times\R\times\R^{1\times d}\times\R^{1\times d}, \ \RE\theta\in (0,1)
   \end{array}
\end{equation}
with $\bar k =0$ or along with the following additional inequality
\begin{equation}\label{eq:A.6}
   \Dis f(\omega,t,y,z)\geq \frac{\bar\gamma}{2}|z|^2-\beta|y|-\hat\alpha_t(\omega), \ \ \RE (y,z)\in \R\times\R^{1\times d}.
\end{equation}
\end{enumerate}

\begin{enumerate}
\renewcommand{\theenumi}{(A2')}
\renewcommand{\labelenumi}{\theenumi}

\item\label{A:A2'} Let \ref{A:A2} hold with \eqref{eq:A.5} and \eqref{eq:A.6} being replaced, respectively, with the following two inequalities
\begin{equation}\label{eq:A.7}
  \begin{array}{l}
   \hspace*{-0.4cm}\Dis -{\bf 1}_{\{y_1-\theta y_2<0\}}\left(f(\omega,t,y_1,z_1)-\theta f(\omega,t,y_2,z_2)\right)\vspace{0.1cm}\\
   \hspace*{-0.4cm}\ \ \leq \Dis (1-\theta)\left(\hat\alpha_t(\omega)+ \beta \left|\frac{y_1-\theta y_2}{1-\theta}\right|+\frac{\gamma}{2} \left|\frac{z_1-\theta z_2}{1-\theta}\right|^2+k(|y_1|+|y_2|)+\bar k (|z_1|^{1+\delta}+|z_2|^{1+\delta}\right)\\
   \Dis \hspace{2.2cm}\RE (y_1,y_2,z_1,z_2)\in \R\times\R\times\R^{1\times d}\times\R^{1\times d},\ \RE\theta\in (0,1)
   \end{array}
\end{equation}
and
\begin{equation}\label{eq:A.8}
   \Dis f(\omega,t,y,z)\leq -\frac{\bar\gamma}{2}|z|^2+\beta|y|+\hat\alpha_t(\omega), \ \ \RE (y,z)\in \R\times\R^{1\times d}.
\end{equation}
\end{enumerate}

\cref{lem:A.1} below is a slight generalization of Theorem 5 in \citet{FanHuTang2020CRM}.  The proof is almost the same, and is thus  omitted here.

\begin{lem}\label{lem:A.1}
Assume that \eqref{eq:A.2} holds and the generator $f$ satisfies assumptions \ref{A:A1} and \ref{A:A2}. Then BSDE \eqref{eq:A.1} admits a unique solution $(Y,Z)$ such that for each $p\geq 1$,
\begin{equation}\label{eq:A.9}
\E\left[\exp\left(p\sup_{t\in \T}|Y_t|\right)\right] <+\infty
\end{equation}
and
\begin{equation}\label{eq:A.10}
\E\left[\left(\int_0^T |Z_s|^2{\rm d}s \right)^{p\over 2} \right] <+\infty.\vspace{0.1cm}
\end{equation}
Moreover, if $f$ satisfies assumption \eqref{eq:A.6}, then for each $\eps\in (0,\bar\gamma/9]$, we also have
\begin{equation}\label{eq:A.11}
\E\left[\exp\left(\frac{\bar \gamma}{2}\eps\int_0^T |Z_s|^2{\rm d}s \right)\right] \leq \E\left[\exp\left(3\eps(2+\beta T)\sup_{t\in [0,T]}|Y_t|+3\eps\int_0^T \hat\alpha_s {\rm d}s\right)\right]<+\infty. \vspace{0.4cm}
\end{equation}
\end{lem}

Concerning the previous assumptions, we have the following remark.

\begin{rmk}\label{rmk:A.2}
(i) Let
$$\tilde f(\omega,t,y,z):=-f(\omega,t,-y,-z).$$
Then, both \ref{A:A1} and \ref{A:A2} hold for $f$ if and only if both \ref{A:A1} and \ref{A:A2'} hold for $\tilde f$. Consequently, the assertions of Lemma \ref{lem:A.1} still hold when \ref{A:A2} is replaced with \ref{A:A2'}, and \eqref{eq:A.6} is replaced with \eqref{eq:A.8}. \vspace{0.2cm}

(ii) In the same way as that in \citet[Proposition 3.5]{FanHu2021SPA}, it is not difficult to verify that \eqref{eq:A.5} holds with $\bar k=0$ for $f$ but with different parameters provided that it satisfies \eqref{eq:A.3} and any one of the following three conditions:
\begin{enumerate}
\item [a)] Almost everywhere in $\Omega\times\T$, $f(\omega,t,\cdot,\cdot)$ is convex;
\item [b)] Almost everywhere in $\Omega\times\T$, for each $(y,z)\in \R\times\R^{1\times d}$, $f(\omega,t,\cdot,z)$ is Lipschitz continuous and $f(\omega,t, y,\cdot)$ is convex;
\item [c)] $f(\omega,t,y,z):=l(y)q(z)$, where both functions $l:\R\To\R$ and $q:\R^{1\times d}\To\R$ are bounded Lipschitz continuous, and the function $q(z)$ has a bounded support.
\end{enumerate}
Furthermore, it has been shown in \citet[Proposition 3]{FanHuTang2020CRM} that \eqref{eq:A.5} holds for $f$ but with different parameters provided that it satisfies \eqref{eq:A.3} and either of the following two conditions:
\begin{enumerate}
\item [d)] Almost everywhere in $\Omega\times\T$, for each $(y,z)\in \R\times\R^{1\times d}$, $f(\omega,t,\cdot,z)$ is Lipschitz continuous, and $f(\omega,t, y,\cdot)$ is $\delta$-locally Lipschitz continuous,  i.e., $\as$, for each $(y_i,z_i)\in \R\times\R^{1\times d}$, $i=1,2$, we have
\begin{equation}\label{eq:A.12}
|f(\omega,t,y_1,z_1)-f(\omega,t,y_2,z_2)|\leq \beta |y_1-y_2|+\gamma (1+|z_1|^\delta+|z_2|^\delta)|z_1-z_2|.\vspace{-0.1cm}
\end{equation}
\item [e)] Almost everywhere in $\Omega\times\T$, for each $(y,z)\in \R\times\R^{1\times d}$, $f(\omega,t,\cdot,z)$ is convex, and $f(\omega,t, y,\cdot)$ is $\delta$-locally Lipschitz continuous, i.e., \eqref{eq:A.12} holds with $y_1=y_2=y$.\vspace{0.1cm}
\end{enumerate}
In addition, we have \eqref{eq:A.5} with $\bar k=0$ (resp. \eqref{eq:A.5} and \eqref{eq:A.6})  for the sum of finite numbers of generators satisfying
\eqref{eq:A.5} with $\bar k=0$ (resp. \eqref{eq:A.5} and \eqref{eq:A.6}) but with different parameters.\vspace{0.2cm}

(iii) Letting $y_1=y_2=y$ and $z_1=z_2=z$ in  \eqref{eq:A.5} and \eqref{eq:A.7}, respectively, yields that
$$
{\bf 1}_{\{y>0\}} f(\omega,t,y,z)\leq 2\bar k+\hat\alpha_t(\omega)+(\beta+2k)|y|+\left(\frac{\gamma}{2}+2\bar k\right)|z|^{2}\vspace{-0.1cm}
$$
and\vspace{-0.1cm}
$$
-{\bf 1}_{\{y<0\}} f(\omega,t,y,z)\leq 2\bar k+\hat\alpha_t(\omega)+(\beta+2k)|y|+\left(\frac{\gamma}{2}+2\bar k\right)|z|^{2},
$$
whose combination implies that $f$ satisfies \eqref{eq:A.3} but with different parameters. In addition, letting first $z_1=z_2=z$ in  \eqref{eq:A.5} (or \eqref{eq:A.7}) and then letting $\theta\To 1$ yields that
$$
{\bf 1}_{\{y_1-y_2>0\}} \left(f(\omega,t,y_1,z)-f(\omega,t,y_2,z)\right)\leq \beta |y_1-y_2|,
$$
which indicates that $f$ satisfies the usual monotonicity condition in $y$.\vspace{0.2cm}
\end{rmk}

\cref{lem:A.3} below provides some a priori estimates on the solution of one-dimensional quadratic BSDEs with unbounded terminal values.

\begin{lem}\label{lem:A.3}
Assume that \eqref{eq:A.2} holds and $(Y_t,Z_t)_{t\in\T}$ is an adapted solution of BSDE \eqref{eq:A.1} such that \eqref{eq:A.9} is satisfied. \vspace{0.2cm}

(i) If it holds that
\begin{equation}\label{eq:A.13}
\as,\ \ {\rm sgn}(Y_t)f(t,Y_t,Z_t)\leq \hat\alpha_t+\beta|Y_t|+\frac{\gamma}{2} |Z_t|^2,
\end{equation}
then we have
\begin{equation}\label{eq:A.15}
\exp\left(\gamma |Y_t|\right)\leq \E_t\left[\exp\left(\gamma e^{\beta T}|\eta|+\gamma e^{\beta T} \int_t^T \hat\alpha_s{\rm d}s\right)\right],\ \ \ t\in\T
\end{equation}
and for each $p>1$,
\begin{equation}\label{eq:A.16}
\E\left[\left(\int_0^T |Z_s|^2{\rm d}s\right)^{p\over 2}\right]\leq C_p \E\left[\left(\int_0^T \hat\alpha_s{\rm d}s\right)^p+\exp\left(p\gamma \sup_{t\in\T}|Y_t|\right)\right],\vspace{0.1cm}
\end{equation}
where $C_p>0$ is a positive constant depending only on $(p,\beta,\gamma,T)$.\vspace{0.2cm}

(ii) If it holds that
\begin{equation}\label{eq:A.17}
 \as,\ \ {\bf 1}_{\{Y_t>0\}}f(t,Y_t,Z_t)\leq \hat\alpha_t+\beta|Y_t|+\frac{\gamma}{2} |Z_t|^2+ \check\alpha_t,
\end{equation}
where $(\check\alpha_t)_{t\in\T}$ is an $(\F_t)$-progressively measurable nonnegative real-valued process satisfying that for some real $q>1$,
\begin{equation}\label{eq:A.18}
\E\left[\left(\int_0^T \check\alpha_t{\rm d}t\right)^{q}\right]<+\infty,
\end{equation}
then for each $t\in\T$ and $p\geq 1$, we have
\begin{equation}\label{eq:A.19}
\begin{array}{l}
\Dis \exp\left(p\gamma Y_t^+\right)\leq \Dis \E_t\left[\exp\left\{p\gamma e^{\beta T}\left(\eta^+ +\int_t^T \hat\alpha_s{\rm d}s\right)\right\}\right.\vspace{0.2cm}\\
\Dis \hspace*{3cm} +\left. p\gamma e^{\beta T}
\exp\left\{p\gamma e^{\beta T} \left(\sup_{t\in\T} Y_t^+ +\int_0^T\hat\alpha_s{\rm d}s\right)\right\}\int_0^T \check\alpha_t {\rm d}t\right]. \vspace{0.2cm}
\end{array}
\end{equation}
\end{lem}

\begin{proof}
In light of \eqref{eq:A.2}, \eqref{eq:A.9}, \eqref{eq:A.13}, \eqref{eq:A.17} and \eqref{eq:A.18}, inequalities \eqref{eq:A.15} and \eqref{eq:A.19} can be deduced by applying It\^{o}-Tanaka's formula to
$$\exp\left(\gamma e^{\beta t}|Y_t|+\gamma\int_0^t e^{\beta s}\hat\alpha_s {\rm d}s\right)\ \ {\rm and} \ \ \exp\left(p\gamma e^{\beta t} Y_t^+ +p\gamma\int_0^t e^{\beta s}\hat\alpha_s {\rm d}s\right),\ {\rm respectively},$$
and using an identical argument as that in the proof of \citet[(i) of Proposition 1]{FanHuTang2020CRM}. Furthermore, in light of \eqref{eq:A.2}, \eqref{eq:A.9} and \eqref{eq:A.13}, using It\^{o}-Tanaka's formula to compute
$${1\over \gamma^2}\big(\exp(\gamma |Y_t|)-\gamma |Y_t|-1\big)$$
and then applying the BDG inequality and H\"{o}lder's inequality, we can derive \eqref{eq:A.16}; see more details in \citet[Corollary 4]{BriandHu2008PTRF}.\vspace{0.2cm}
\end{proof}

\cref{lem:A.4} below establishes a stability result on one-dimensional quadratic BSDEs with unbounded terminal values and general assumptions on the generator, which strengthens in \citet[Proposition 7]{BriandHu2008PTRF} and inspires the proofs of \cref{thm:2.2,thm:2.5}.

\begin{lem}\label{lem:A.4}
Let $f$ satisfy assumptions \ref{A:A1} and \ref{A:A2} (resp. \ref{A:A2'}), and $(Y_t,Z_t)_{t\in\T}$ be the unique solution to BSDE \eqref{eq:A.1} such that \eqref{eq:A.9} holds. And, for each $m\geq 1$, let $\eta^{(m)}$ be a real-valued $\F_T$-measurable random variable, $f^{(m)}$ be a generator which satisfies \ref{A:A1} and \ref{A:A2} (resp. \ref{A:A2'}) with an $(\F_t)$-progressively measurable nonnegative real-valued  process $\hat\alpha^{(m)}$ instead of $\hat\alpha$, and $(Y^{(m)}_t,Z^{(m)}_t)_{t\in\T}$ be the unique solution to BSDE associated with $(\eta^{(m)},f^{(m)})$ such that \eqref{eq:A.9} with $Y^{(m)}$ instead of $Y$ holds. We further assume that for each $p\geq 1$,\vspace{0.1cm}
\begin{equation}\label{eq:A.20}
\E\left[\exp\left\{p|\eta|+p\int_0^T \hat\alpha_t {\rm d}t \right\}\right]+\sup_{m\geq 1}\E\left[\exp\left\{p|\eta^{(m)}|+p\int_0^T \hat\alpha^{(m)}_t {\rm d}t \right\}\right]<+\infty.\vspace{0.1cm}
\end{equation}
If $\ps$, $\eta^{(m)}\To \eta$ as $m\To \infty$ and there exists a real $q>1$ such that
\begin{equation}\label{eq:A.21}
\lim\limits_{m\To \infty} \E\left[\left(\int_0^T \left| f^{(m)}(s,Y_s, Z_s)-f(s,Y_s, Z_s)\right| {\rm d}s\right)^q\right]=0,
\end{equation}
then for each $p\geq 1$, we have
\begin{equation}\label{eq:A.22}
\lim\limits_{m\To \infty} \E\left[\exp\left(p \sup_{t\in\T}|Y^{(m)}_t-Y_t|\right)\right]=1.\vspace{0.2cm}
\end{equation}
And, if the function $\varphi_t(\omega,r)$ in \eqref{eq:A.4} also satisfies that for some constant $k>0$,
\begin{equation}\label{eq:A.23}
\as,\ \ \varphi_t(\omega,x)\leq k\exp (kx),\ \ x\in [0,+\infty),
\end{equation}
then for each $p\geq 1$,
\begin{equation}\label{eq:A.24}
\lim_{m\To \infty} \E\left[\left(\int_0^T|Z^{(m)}_s-Z_s|^2{\rm d}s\right)^{p\over 2}
\right]=0.\vspace{0.3cm}
\end{equation}
\end{lem}

\begin{proof}
We only prove the case that all $f^{(m)}$ satisfy assumption \ref{A:A1} and  inequalities~\eqref{eq:A.3}-\eqref{eq:A.6} of assumption \ref{A:A2}. The other case can be proved in an identical way. For simplicity, we assume without loss of generality that $\gamma>1$. Otherwise, we replace  $\gamma$ with $1+\gamma$.\vspace{0.2cm}

First of all, in light of \eqref{eq:A.3}, from (i) of \cref{lem:A.3} and \eqref{eq:A.20} along with Doob's maximal inequality for martingale we can conclude that the sequence $(Y^{(m)}_t,Z^{(m)}_t)_{t\in\T}$ satisfies
\begin{equation}\label{eq:A.25}
\RE p\geq 1,\ \ \sup_{m\geq 1} \E\left[\exp\left(p \sup_{t\in\T}|Y^{(m)}_t|\right)
+\left(\int_0^T|Z^{(m)}_s|^2{\rm d}s\right)^{p\over 2}
\right]<+\infty.
\end{equation}
And, due to the inequality~\eqref{eq:A.6}, from \eqref{eq:A.11} in \cref{lem:A.1} along with \eqref{eq:A.20}, \eqref{eq:A.25} and H\"{o}lder's inequality we know that
\begin{equation}\label{eq:A.26}
\E\left[\exp\left(\frac{\bar \gamma^2}{18}\int_0^T |Z_s|^2{\rm d}s \right)\right]+\sup_{m\geq 1}\E\left[\exp\left(\frac{\bar \gamma^2}{18}\int_0^T |Z_s^{(m)}|^2{\rm d}s \right)\right]<+\infty. \vspace{0.1cm}
\end{equation}
Furthermore, it follows from Young's inequality that for each $a,b,\eps>0$ and $\delta\in [0,1)$,
\begin{equation}\label{eq:A.27}
ab^{1+\delta}=\left(\left(\frac{1+\delta}{2\eps}\right)^{\frac{1+\delta}{1-\delta}} a^{\frac{2}{1-\delta}}\right)^{\frac{1-\delta}{2}}
\left(\frac{2\eps}{1+\delta}b^2\right)^{\frac{1+\delta}{2}}\leq \eps b^2+\frac{1-\delta}{2}\left(\frac{1+\delta}{2\eps}\right)^{\frac{1+\delta}{1-\delta}} a^{\frac{2}{1-\delta}}.
\end{equation}
By virtue of the last inequality and H\"{o}lder's inequality, from \eqref{eq:A.9}, \eqref{eq:A.20}, \eqref{eq:A.25} and \eqref{eq:A.26} we obtain the existence of a constant $C(p)\geq 1$ depending only on $(p,\beta,\gamma,\bar\gamma, \delta, k,\bar k, T)$ such that for each $p>1$,
\begin{equation}\label{eq:A.28}
\begin{array}{l}
\Dis\left(\frac{p}{p-1}\right)^{p} \sup_{m\geq 1}\E\Bigg[\exp\Bigg(p\int_0^T \bigg(\hat\alpha^{(m)}_t+k|Y_t|+\bar k|Z_t|^{1+\delta}\vspace{0.1cm}\\
\hspace{4.9cm}\Dis +\hat\alpha_t+k|Y^{(m)}_t|+\bar k|Z^{(m)}_t|^{1+\delta}\bigg){\rm d}t \Bigg)\Bigg]\leq C(p). \vspace{0.2cm}
\end{array}
\end{equation}

Next, for each fixed $\theta\in (0,1)$ and $m\geq 1$, define
$$
\delta_\theta^m  U:=\frac{Y^{(m)}-\theta Y}{1-\theta}\ \  {\rm and} \ \ \delta_\theta^m  V:=\frac{Z^{(m)}-\theta Z}{1-\theta}.
$$
Then the pair $(\delta_\theta^m  U,\delta_\theta^m  V)$ verifies the following BSDE:\vspace{0.1cm}
$$
  \delta_\theta^m  U_t=\delta_\theta^m  U_T +\int_t^T \delta_\theta^m f (s) {\rm d}s-\int_t^T \delta_\theta^m  V_s \cdot {\rm d}B_s, \ \ \ \ t\in\T,
$$
where
\begin{equation}\label{eq:A.29}
\as,\ \ \Dis \delta_\theta^m f(t):=\Dis \frac{f^{(m)}(t, Y^{(m)}_t, Z^{(m)}_t)-\theta f^{(m)}(t, Y_t, Z_t)}{1-\theta}+\ D^m_\theta(t)
\end{equation}
with
\begin{equation}\label{eq:A.30}
D^m_\theta(t):=\frac{\theta}{1-\theta}\left( f^{(m)}(t,Y_t, Z_t)-f(t,Y_t, Z_t) \right).
\vspace{0.1cm}
\end{equation}
Since the inequality~\eqref{eq:A.5} is satisfied for the generator $f^{(m)}$, it follows from \eqref{eq:A.29} that $\as$, for each $(y,z)\in \R\times \R^{1\times d}$,
\begin{equation}\label{eq:A.31}
\begin{array}{l}
\Dis{\bf 1}_{\{\delta_\theta^m  U_t>0\}} \delta_\theta^m f(t)\leq \hat\alpha^{(m)}_t+\beta |\delta_\theta^m  U_t|+{\gamma\over 2} |\delta_\theta^m  V_t|^2+k(|Y_t|+|Y^{(m)}_t|)\vspace{0.1cm}\\
\hspace{3.1cm}\Dis +\bar k (|Z_t|^{1+\delta}+|Z^{(m)}_t|^{1+\delta})+|D^m_\theta(t)|.
\end{array}
\end{equation}
In light of \eqref{eq:A.31}, \eqref{eq:A.28}, \eqref{eq:A.20}, \eqref{eq:A.21} and \eqref{eq:A.25}, by using (ii) of \cref{lem:A.3}, Doob's maximal inequality and H\"{o}lder's inequality we can deduce that for each $\theta\in (0,1)$ and $\bar q\in (1,q)$, noting that $C(p)\geq 1$ for $p\geq 1$ and $e^x\geq 1$ for $x\geq 0$,
\begin{equation}\label{eq:A.32}
\begin{array}{l}
\Dis \E\left[\exp\left(\bar q\gamma \sup_{t\in\T}\left(\delta_\theta^m  U_t\right)^+\right)\right] \vspace{0.2cm}\\
\ \ \leq \Dis 2^{\bar q} C\left(2\bar q\gamma e^{\beta T}\right)\E\left[\exp\left\{2\bar q\gamma e^{\beta T}\left(\delta_\theta^m  U_T\right)^+ \right\}\right]+(2\gamma)^{\bar q} e^{\bar q\beta T}C\left({2q\bar q\gamma\over q-\bar q} e^{\beta T}\right)\vspace{0.2cm}\\
\Dis \ \ \ \ \cdot \sup_{m\geq 1}\E\left[\exp\left\{{2q\bar q\gamma\over q-\bar q}e^{\beta T}\sup_{t\in\T} \left(\delta_\theta^m  U_t\right)^+\right\}\right]\cdot \left(\E\left[\left(\int_0^T |D^m_\theta(t)| {\rm d}t\right)^{q}\right]\right)^{\bar q\over q}.
\end{array}
\end{equation}
Furthermore, in light of \eqref{eq:A.20}, \eqref{eq:A.30}, \eqref{eq:A.21} and the assumption that $\ps$, $\eta^{(m)}\To \eta$ along with the facts that $x\leq e^x$ for $x\geq 0$ and
$$
\left(\frac{\eta^{(m)}-\theta \eta}{1-\theta}\right)^+\leq \frac{|\eta^{(m)}- \eta|}{1-\theta}+|\eta|,
$$
letting $m\To \infty$ in \eqref{eq:A.32} yields that for each $\theta\in (0,1)$ and $\bar q\in (1,q)$,
\begin{equation}\label{eq:A.33}
\limsup_{m\To \infty}\E\left[\sup_{t\in\T}\left(Y_t^{(m)}-\theta Y_t\right)^+\right] \leq {1-\theta\over \bar q\gamma}2^{\bar q} C\left(2\bar q\gamma e^{\beta T}\right)\E\left[\exp\left\{2\bar q\gamma e^{\beta T}|\eta| \right\}\right].\vspace{0.2cm}
\end{equation}

On the other hand, in the last paragraph by interchanging the position between $(Y,Z)$ and $(Y^{(m)},Z^{(m)})$ and replacing \eqref{eq:A.29} and\eqref{eq:A.30} with, respectively,
$$
\Dis \delta_\theta^m f(t):=\Dis \ D^m_\theta(t)+\frac{f^{(m)}(t, Y_t, Z_t)-\theta f^{(m)}(t, Y^{(m)}_t, Z^{(m)}_t)}{1-\theta}
$$
and
$$
D^m_\theta(t):=\frac{f(t,Y_t, Z_t)-f^{(m)}(t,Y_t, Z_t)}{1-\theta},\vspace{0.1cm}
$$
we can also conclude that for each $\theta\in (0,1)$ and $\bar q\in (1,q)$,
\begin{equation}\label{eq:A.34}
\limsup_{m\To \infty}\E\left[\sup_{t\in\T}\left(Y_t-\theta Y_t^{(m)}\right)^+\right] \leq {1-\theta\over \bar q\gamma}2^{\bar q} C\left(2\bar q\gamma e^{\beta T}\right)\E\left[\exp\left\{2\bar q\gamma e^{\beta T}|\eta| \right\}\right].
\end{equation}

Observe that
$$
Y^{(m)}_t-Y_t \leq (Y^{(m)}_t-\theta Y_t)^+ +(1-\theta)|Y_t|\ \ {\rm and}\ \ Y_t-Y^{(m)}_t \leq (Y_t-\theta Y^{(m)}_t)^+ +(1-\theta)|Y^{(m)}_t|.
$$
It follows from \eqref{eq:A.33} and \eqref{eq:A.34} that for each $\theta\in (0,1)$ and $\bar q\in (1,q)$,
$$
\begin{array}{l}
\Dis \limsup_{m\To \infty}\E\left[\sup_{t\in\T}|Y_t^{(m)}-Y_t|\right]\vspace{0.2cm} \\
\ \ \leq \Dis (1-\theta) \left({2^{\bar q+1}\over \bar q\gamma} C\left(2\bar q\gamma e^{\beta T}\right)\E\left[\exp\left\{2\bar q\gamma e^{\beta T}|\eta| \right\}\right]+ \sup_{m\geq 1}\E\left[\sup_{t\in\T}\left(|Y_t|+|Y_t^{(m)}|\right)\right] \right).
\end{array}
$$
In light of \eqref{eq:A.25}, sending $\theta\To 1$ in the last inequality yields that
\begin{equation}\label{eq:A.36}
\lim\limits_{m\To \infty} \E\left[\sup_{t\in\T}|Y_t^{(m)}-Y_t|\right]=0,
\end{equation}
which along with \eqref{eq:A.25} leads to the desired conclusion \eqref{eq:A.22}.\vspace{0.2cm}

Finally, let \eqref{eq:A.23} be further satisfied. We prove that \eqref{eq:A.24} holds. In fact, it follows from It\^{o}'s formula that for each $m\geq 1$,
\begin{equation}\label{eq:A.37}
\begin{array}{l}
\Dis \E\left[\int_0^T |Z^{(m)}_s-Z_s|^2 {\rm d}s\right]\vspace{0.1cm}\\
\ \ \leq \Dis \E\left[|\eta^{(m)}-\eta|^2+2\sup_{t\in\T}|Y^{(m)}_t-Y_t| \int_0^T \left|f^{(m)}(s,Y^{(m)}_s, Z^{(m)}_s)-f(s,Y_s, Z_s)\right| {\rm d}s \right].
\end{array}
\end{equation}
On the other hand, by virtue of \eqref{eq:A.4} with \eqref{eq:A.23} for $f^{(m)}$ and $f$ as well as \eqref{eq:A.20} and \eqref{eq:A.25} we can conclude that
\begin{equation}\label{eq:A.38}
\sup_{m\geq 1}\E\left[\left(\int_0^T \left|f^{(m)}(s,Y^{(m)}_s, Z^{(m)}_s)-f(s,Y_s, Z_s)\right| {\rm d}s\right)^2\right]<+\infty.
\end{equation}
Then, by using H\"{o}lder's inequality, the desired conclusion \eqref{eq:A.24} for the case of $p=2$ follows immediately from \eqref{eq:A.37}, \eqref{eq:A.38}, \eqref{eq:A.22} and \eqref{eq:A.25}. Furthermore, by virtue of \eqref{eq:A.25} we can verify the general case of $p\geq 1$. The proof is then complete.
\end{proof}



\setlength{\bibsep}{2pt}
\bibliographystyle{elsarticle-harv}

\begin{thebibliography}{26}
\expandafter\ifx\csname natexlab\endcsname\relax\def\natexlab#1{#1}\fi
\expandafter\ifx\csname url\endcsname\relax
  \def\url#1{\texttt{#1}}\fi
\expandafter\ifx\csname urlprefix\endcsname\relax\def\urlprefix{URL }\fi

\bibitem[{Barrieu and El~Karoui(2013)}]{BarrieuElKaroui2013AoP}
Barrieu P., El~Karoui N., 2013. Monotone stability of quadratic
  semimartingales with applications to unbounded general quadratic {BSDE}s.
  Ann. Probab. 41~(3B), 1831--1863.

\bibitem[{Bismut(1973)}]{Bismut1973JMAA}
Bismut J., 1973. Conjugate convex functions in optimal stochastic control. J.
  Math. Anal. Appl. 44~(2), 384--404.

\bibitem[{Briand, Delyon, Hu, Pardoux and Stoica(2003)}]{Briand-etal2003SPA}
Briand P., Delyon B., Hu Y., Pardoux E., Stoica L, 2003. $L^p$ solutions of backward stochastic differential equations. Stochastic Process. Appl. 108~(1), 109--29.

\bibitem[{Briand and Elie(2013)}]{BriandElie2013SPA}
Briand P., Elie R., 2013. A simple constructive approach to quadratic {BSDE}s
  with or without delay. Stochastic Process. Appl. 123, 2921--2939.

\bibitem[{Briand and Hu(2006)}]{BriandHu2006PTRF}
Briand P., Hu Y., 2006. {BSDE} with quadratic growth and unbounded terminal
  value. Probab. Theory Related Fields 136~(4), 604--618.

\bibitem[{Briand and Hu(2008)}]{BriandHu2008PTRF}
Briand P., Hu Y., 2008. {Q}uadratic {BSDE}s with convex generators and
  unbounded terminal conditions. Probab. Theory Related Fields 141~(3),
  543--567.

\bibitem[{Cheridito and Nam(2015)}]{cheridito2015Stochastics}
Cheridito P., Nam K., 2015. Multidimensional quadratic and subquadratic
  {BSDE}s with special structure. Stochastics An International Journal of
  Probability and Stochastic Processes 87~(5), 871--884.

\bibitem[{Delbaen, Hu and Richou(2011)}]{DelbaenHuRichou2011AIHPPS}
Delbaen F., Hu Y., Richou A., 2011. On the uniqueness of solutions to
  quadratic {BSDE}s with convex generators and unbounded terminal conditions.
  Ann. Inst. Henri Poincar\'{e} Probab. Stat. 47~(2), 559--574.

\bibitem[{Delbaen, Hu and Richou(2015)}]{DelbaenHuRichou2015DCD}
Delbaen F., Hu Y., Richou A., 2015. On the uniqueness of solutions to
  quadratic {BSDE}s with convex generators and unbounded terminal conditions:
  the critical case. Discrete Contin. Dyn. Syst. 35~(11), 5273--5283.

\bibitem[El Karoui, Peng and Quenez(1997)]{ElKarouiPengQuenez1997MF}
El Karoui N., Peng S., Quenez M.-C., 1997. Backward stochastic differential equations in finance. Math. Finance 7~(1), 1--71. 

\bibitem[{Fan(2016)}]{Fan2016SPA}
Fan S., 2016. Bounded solutions, ${L}^{p}\ (p>1)$ solutions and ${L}^1$
  solutions for one-dimensional {BSDE}s under general assumptions. Stochastic
  Process. Appl. 126, 1511--1552.

\bibitem[{Fan and Hu(2021)}]{FanHu2021SPA}
Fan S., Hu Y., 2021. Well-posedness of scalar BSDEs with sub-quadratic generators and related PDEs. Stochastic Process. Appl., 131, 21--50.

\bibitem[{Fan, Hu and Tang(2020)}]{FanHuTang2020CRM}
Fan S., Hu Y., Tang S., 2020. On the uniqueness of solutions to quadratic
  {BSDE}s with non-convex generators and unbounded terminal conditions. C. R. Math. Acad. Sci. Paris  358~(2),  227--235.

\bibitem[{Fan, Hu and Tang(2023)}]{FanHuTang2023JDE}
Fan S., Hu Y., Tang S., 2023. Multi-dimensional backward
  stochastic differential equations of diagonally quadratic generators: the
  general result. J. Differ. Equ., 368, 105--140.

\bibitem[{Fan, Hu and Tang(2025)}]{FanHuTang2025EJP}
Fan S., Hu Y., Tang S., 2025. Multi-dimensional non-Markovian backward stochastic differential equations of interactively quadratic generators. Electron. J. Probab. 30~(169), 1--57.

\bibitem[{Fan, Hu and Tang(2026)}]{FanHuTang2026MF}
Fan S., Hu Y., Tang S., 2026. One-dimensional nonlinear backward stochastic differential equations: A unified theory and applications. Math. Finance, forthecoming, https://doi.org/10.1111/mafi.70046.

\bibitem[{Fan, Wang and Yong(2025)}]{FanWangYong2025JDE}
Fan S., Wang T., Yong J., 2025. Multi-dimensional super-linear backward stochastic volterra integral equations. J. Differ. Equ. 437, Paper No. 113350, 62pp.

\bibitem[{Frei(2014)}]{Frei2014SPA}
Frei C., 2014. Splitting multidimensional {BSDE}s and finding local
  equilibria. Stochastic Process. Appl. 124~(8), 2654--2671.

\bibitem[{Frei and Dos~Reis(2011)}]{FreiDosReis2011MFE}
Frei C., Dos~Reis G., 2011. A financial market with interacting investors:
  does an equilibrium exist? Mathematics and Financical economics 4~(3),
  161--182.

\bibitem[Hao, Hu, Tang and Wen(2025)]{HaoHuWen2025AAP}
Hao T., Hu Y., Tang S., Wen J., 2025. Mean-field backward stochastic differential equations and nonlocal PDEs with quadratic growth. Ann. Appl. Probab. 35~(3), 2128--2174.

\bibitem[{Harter and Richou(2019)}]{HarterRichou2019EJP}
Harter J., Richou A., 2019. A stability approach for solving multidimensional
  quadratic {BSDE}s. Electron. J. Probab. 24~(4), 1--51.

\bibitem[Hu, Imkeller and M\"{u}ller(2005)]{HuImkellerMuller2005AAP}
Hu Y., Imkeller P., M\"{u}ller M., 2005. Utility maximization in incomplete markets. Ann. Appl. Probab. 15, 1691--1712. 

\bibitem[{Hu and Tang(2016)}]{HuTang2016SPA}
Hu Y., Tang S., 2016. Multi-dimensional backward stochastic differential
  equations of diagonally quadratic generators. Stochastic Process. Appl. 126~(4), 1066--1086.

\bibitem[Hu, Tang and Wang(2022)]{HuTangWang2022SPA}
Hu Y., Tang S., Wang F., 2022. Quadratic {G}-BSDEs with convex generators and unbounded terminal conditions. Stochastic Process. Appl. 153, 363--390.

\bibitem[{Jamneshan, Kupper and Luo(2017)}]{JamneshanKupperLuo2017ECP}
Jamneshan A., Kupper M., Luo P., 2017. Multidimensional quadratic {BSDE}s
  with separated generators. Electronic Communications in Probability 22~(58),
  1--10.

\bibitem[{Kardaras, Xing and {\v{Z}}itkovi{\'c}(2022)}]{Kardaras2022Book}
Kardaras C., Xing H., {\v{Z}}itkovi{\'c} G., 2022. Incomplete stochastic equilibria with exponential utilities: close to {P}areto optimality. In: Yin, G., Zariphopoulou, T. (eds), Stochastic Analysis, Filtering, and Stochastic Optimization. Springer, Cham, pp.267--292.

\bibitem[{Kobylanski(2000)}]{Kobylanski2000AP}
Kobylanski M., 2000. Backward stochastic differential equations and partial
  differential equations with quadratic growth. Ann. Probab. 28~(2), 558--602.

\bibitem[{Kramkov and Pulido(2016)}]{KramkovPulido2016SIAM}
Kramkov D., Pulido S., 2016. Stability and analytic expansions
  of local solutions of systems of quadratic bsdes with applications to a price
  impact model. SIAM Journal on Financial Mathematics 7~(1), 567--587.

\bibitem[{Kramkov and Pulido(2016)}]{KramkovPulido2016AAP}
Kramkov D., Pulido S., 2016. A system of quadratic {BSDE}s
  arising in a price impact model. Ann. Appl. Probab. 26~(2),
  794--817.

\bibitem[{Luo and Fan(2018)}]{LuoFan2018SD}
Luo H., Fan S., 2018. Bounded solutions for general time interval {BSDE}s
  with quadratic growth coefficients and stochastic conditions. Stoch. Dynam.
  18~(5), Paper No. 1850034, 24pp.

\bibitem[{Luo(2020)}]{Luo2020EJP}
Luo P., 2020. A type of globally solvable {BSDE}s with triangularly quadratic
  generators. Electron. J. Probab. 25~(112), 1--23.

\bibitem[{Pardoux and Peng(1990)}]{PardouxPeng1990SCL}
Pardoux E., Peng S., 1990. Adapted solution of a backward stochastic
  differential equation. Syst. Control Lett. 14~(1), 55--61.

\bibitem[{Peng(1999)}]{Peng1999Book}
Peng S., 1999. Open problems of backward stochastic differential equations, in: Chen S., et al. (Eds.), Control of Distributed Parameter and Stochastic Systems, Hangzhou, 1998, Kluwer Academic Publishers, Boston, 1999, pp. 265--273.

\bibitem[{Tevzadze(2008)}]{Tevzadze2008SPA}
Tevzadze R., 2008. Solvability of backward stochastic differential equations
  with quadradic growth. Stochastic Process. Appl. 118~(3), 503--515.

\bibitem[{Xing and {\v{Z}}itkovi{\'c}(2018)}]{XingZitkovic2018AoP}
Xing H., {\v{Z}}itkovi{\'c} G., 2018. A class of globally solvable
  {M}arkovian quadratic {BSDE} systems and applications. Ann. Probab. 46~(1), 491--550.

\bibitem[Yong and Zhou(1999)]{YongZhou1999}
Yong J., Zhou X., 1999. Stochastic controls. Hamiltonian systems and HJB equations. Applications of Mathematics (New York), 43, Springer-Verlag, New York.

\end{thebibliography}

\def\cprime{$'$}

\end{document}